\documentclass[reqno]{amsart}
\title[Resolvent expansions for three-dimensional discrete Schr\"odinger operators]{
Resolvent Expansions at All Spectral Thresholds for Three-Dimensional Discrete Schr\"odinger Operators}
\author{Yanhong Chen, Sisi Huang and Xiaohua Yao}

\address{Yanhong Chen, Department of Mathematics, Central China Normal University, Wuhan, 430079, P.R. China}
\email{cyh23@mails.ccnu.edu.cn}

\address{Sisi Huang, Department of Mathematics, Central China Normal University, Wuhan, 430079, P.R. China}
\email{hss@mails.ccnu.edu.cn}

\address{Xiaohua Yao, Department of Mathematics, Key Laboratory of Nonlinear Analysis and Applications(Ministry of Education), and Hubei Key Laboratory of Mathematical Sciences, Central China Normal University, Wuhan, 430079, P.R. China}
\email{yaoxiaohua@ccnu.edu.cn}

\date{\today}
\usepackage{color}
\usepackage{amsmath}
\usepackage{amssymb}
\usepackage[mathscr]{euscript}
\usepackage{paralist}
\usepackage{amsthm}
\usepackage{ulem}
\usepackage{bm}
\usepackage{extarrows}
\usepackage{diagbox}
\usepackage{multirow}
\usepackage{makecell}
\usepackage{threeparttable}
\usepackage{tikz}
\usetikzlibrary{patterns}
\usetikzlibrary{decorations.markings}
\usepackage{caption}
\usepackage{verbatim}
\usepackage{pifont}
\usepackage{esint}

\usepackage{cite}
\usepackage{hyperref}
\newtheorem{definition}{Definition}[section]
\newtheorem{theorem}[definition]{Theorem}
\newtheorem{lemma}[definition]{Lemma}
\newtheorem{remark}[definition]{Remark}
\newtheorem{proposition}[definition]{Proposition}

\newcommand\R{\mathbb{R}}
\newcommand\Z{\mathbb{Z}}
\newcommand\B{\mathbb{B}}
\newcommand\C{\mathbb{C}}
\newcommand\N{\mathbb{N}}
\newcommand\T{\mathbb{T}}
\newcommand\mscH{\mathcal{H}}
\newcommand\mcaF{\mathcal{F}}
\newcommand\mcaI{\mathcal{I}}
\newcommand\mcaR{\mathcal{R}}
\newcommand\mcaG{\mathcal{G}}

\newcommand\mcaA{\mathcal{A}}
\newcommand\mcaB{\mathcal{B}}

\newcommand\mcaD{\mathcal{D}}
\newcommand\mcaQ{\mathcal{Q}}
\newcommand\mcaS{\mathcal{S}}

\newcommand\mcaM{\mathcal{M}}
\newcommand\mcaU{\mathcal{U}}
\newcommand\mcaV{\mathcal{V}}
\newcommand\mcaH{\mathcal{H}}

\newcommand\mcaP{\mathcal{P}}
\newcommand\mcaE{\mathcal{E}}
\numberwithin{equation}{section}

\makeatletter
\def\section{\@startsection{section}{1}%
  \z@{.7\linespacing\@plus\linespacing}{.5\linespacing}%
  {\normalfont\bfseries}}
\makeatother
\begin{document}

\begin{abstract}
This paper studies resolvent expansions at all four spectral thresholds
for the discrete Schr\"odinger operator $H=-\Delta+V $ on $\Z^3$, where $\Delta$ denotes the discrete Laplacian on $\mathbb{Z}^3$ and $V(n)$ is a real-valued decaying potential of $\mathbb{Z}^3$. At the thresholds 0, 4, 8, and 12, we first give a complete classification of the threshold types and derive the corresponding resolvent expansions in weighted $\ell^2$-spaces, including the non-generic cases of resonances and eigenvalues. 

As a result, we establish the Jensen-Kato type time-decay estimates of the  evolution $e^{itH} P_{ac}(H)$ to discrete Schr\"odinger equations with thresholds singularities on $\Z^3$, with decay $O(t^{-3/2})$ when all thresholds are regular and $O(t^{-1/2})$ otherwise.
These results provide a three-dimensional lattice analogue of the Jensen--Kato classical  analysis at zero energy threshold on $\mathbb{R}^3$.
\end{abstract}
\maketitle
\tableofcontents

\section{Introduction and Main Results}
\subsection{Introduction}
\ \

\vskip0.1cm
The analysis of threshold energies for Schr\"odinger operators plays a fundamental role in spectral and scattering theory. Near a threshold,
the resolvent may develop singularities in the presence of resonant states or eigenfunctions, and these singularities  directly influence the long-time behavior of the evolution through Stone's formula
\begin{equation}\label{eq:Stone-formula-intro}
e^{-itH}P_{\mathrm{ac}}(H)
=
\frac{1}{2\pi i}
\int_{\sigma_{\mathrm{ac}}(H)}
e^{-it\lambda}
\bigl(
\mcaR_V^+(\lambda)-\mcaR_V^-(\lambda)
\bigr)\,d\lambda.
\end{equation}
Consequently, threshold resolvent expansions provide a necessary  tool for
the study of dispersive estimates and
the $L^p$-boundedness of wave operators.

In the case of the continuous Schr\"odinger operator $
H_c=-\Delta_{\mathbb{R}^d}+V(x),$
the only spectral threshold is zero. In dimension $d=3$,
Jensen and Kato \cite{JK79} first gave the famous analysis of the perturbed
resolvent near zero and classified the threshold as a regular point or
an exceptional point of the first, second, or third kind. Under suitable
assumptions on the potential, the resolvent
$
\mathscr{R}(\omega)
=
(-\Delta_{\mathbb{R}^3}+V-\omega)^{-1}
$
admits the following expansions in appropriate weighted Sobolev spaces:
\begin{align}\label{resol expan of JK79}
\begin{split}
\mathscr{R}(\omega)&=B^0_0+O(\omega^{\frac{1}{2}}),\qquad\qquad\qquad\quad \ \ \text{if 0 is a regular point},\\
\mathscr{R}(\omega)&=-iB^1_{-1}\omega^{-\frac{1}{2}}+O(1),\qquad \qquad\ \ \   \text{if 0 is an exceptional point of the first kind},\\
\mathscr{R}(\omega)&=-\mathscr{P}_0\omega^{-1}-iB^2_{-1}\omega^{-\frac{1}{2}}+O(1),\quad \text{if 0 is an exceptional point of the second kind},\\
\mathscr{R}(\omega)&=-\mathscr{P}_0\omega^{-1}-iB^3_{-1}\omega^{-\frac{1}{2}}+O(1),\quad \text{if 0 is an exceptional point of the third kind},
\end{split}
\end{align}
where $\mathscr{P}_0$ is the orthogonal projection onto the zero
eigenspace. These expansions show explicitly how the threshold spectral
structure governs the large-time dynamics of $e^{itH_c}$: in the regular case one
obtains the classical Jensen--Kato weighted decay of order
$t^{-3/2}$, whereas a zero-energy resonance produces a slower
$t^{-1/2}$ contribution.

The threshold theory was subsequently developed in other dimensions
and became an important ingredient in dispersive and scattering
analysis. Jensen \cite{Jen80,Jen84} treated higher-dimensional cases,
while Jensen and Nenciu \cite{JN01} developed a unified reduction
scheme and obtained complete threshold expansions in dimensions one and
two. The effects of threshold singularities on dispersive estimates were further studied in
\cite{Mur82,JSS91,GS04,Sch05,ES04,ES06,EG13,EG13w,EGG14,GG15,GG17,RS04,Sch07},
while threshold behavior also plays an important role in the  $L^p$ theory of wave operators
\cite{Yaj06}. More recently, threshold analysis has been extended to higher-order Schr\"odinger operators; see, for example,
\cite{FSWY20,CHHZ24} and the references therein.


The purpose of the present paper is to develop an analogous threshold
theory for discrete Schr\"odinger operators on the lattice. Let
$H_0=-\Delta_{\mathbb{Z}^d}$ be 
the discrete Laplacian by
\begin{equation}\label{eq:discrete-laplacian}
(\Delta_{\mathbb{Z}^d}\phi)(n)
:=
\sum_{j=1}^d
\big(
\phi(n+e_j)+\phi(n-e_j)-2\phi(n)
\big),
\qquad
n\in\mathbb{Z}^d,
\end{equation}
with
\[
e_j=(0,\ldots,0,\underset{j}{1},0,\ldots,0).
\]
Under the discrete Fourier transform, $H_0$ is unitarily equivalent to
the multiplication operator with symbol
\begin{equation}\label{eq:discrete-symbol}
\mathcal{M}_d(x)
=
\sum_{j=1}^{d}
\bigl(2-2\cos x_j\bigr),
\qquad
x=(x_1,\ldots,x_d)\in(-\pi,\pi]^d.
\end{equation}
Consequently, $
\sigma(H_0)=[0,4d]$,
and the critical values of the dispersion relation are
\[
\Lambda_d
=
\{0,4,8,\ldots,4d\}.
\]
These critical values constitute the spectral thresholds of the
discrete Laplacian.

This is one of the fundamental differences between the Euclidean and
lattice problems. For the standard Euclidean Laplacian there is a
single spectral threshold at zero, whereas the lattice dispersion
relation possesses several critical values, including thresholds in
the interior of the continuous spectrum. Moreover, these thresholds
have different local geometries. The endpoints $0$ and $4d$ arise
from elliptic critical points of $\mathcal{M}_d$, whereas the interior
thresholds are generated by saddle points and therefore have a
hyperbolic character. 

The distinctive features of these thresholds are already reflected in the study of the free resolvent kernel of $-\Delta_{\Z^d}$, which takes the form
\begin{equation}\label{eq:free-lattice-resolvent}
R_0(z,n)
=
\frac{1}{(2\pi)^d}
\int_{(-\pi,\pi]^d}
\frac{e^{in\cdot x}}
{\mathcal{M}_d(x)-z}\,dx.
\end{equation}
Unlike the Euclidean setting, where the rotational invariance of the symbol $|x|^2$
allows a radial reduction and explicit free resolvent formulas in terms of elementary or Bessel functions. The lattice integral in \eqref{eq:free-lattice-resolvent} generally lacks such a straightforward radial reduction, except in the one-dimensional case. Over recent years, however, substantial progress has been achieved, most notably by Ito and Jensen \cite{IJ19,IJ21}, who analyzed the resolvent behavior at both endpoint and embedded thresholds for multidimensional discrete Laplacians.  
Their results reveal how the local elliptic or ultrahyperbolic geometry of the dispersion relation enters the threshold asymptotics, thus providing a key free resolvent foundation for the study of perturbed lattice operators.

For perturbed discrete Schr\"odinger operators \(H=-\Delta_{\mathbb{Z}^d}+V\), the explicit free resolvent kernel facilitates a detailed threshold analysis in dimension \(d=1\). Resolvent expansions for generic potentials were derived in \cite{KKK06,PS08}, while Ito and Jensen \cite{IJ15} later gave a complete classification of the thresholds \(0\) and \(4\), including non-generic cases. In contrast, the higher-dimensional setting \(d\geq2\) is considerably more intricate.

For perturbed discrete case in dimension \(d=2\), threshold behavior in the regular case and the associated long-time dynamics were studied by Komech, Kopylova, and Vainberg in \cite{KKV08}. Subsequently, in dimension \(d=3\), Kopylova \cite{Kop09} derived Puiseux-type expansions of the perturbed resolvent near the four critical energies \(\Lambda_3=\{0,4,8,12\}\) under the assumption that all thresholds are regular. These expansions yield the weighted long-time decay estimate
\begin{equation}\label{Regular-Kato-Jensen}
\bigl\|
e^{-itH}P_{\mathrm{ac}}(H)
\bigr\|_{\ell^{2,\sigma}\to\ell^{2,{-\sigma}}}
=
O(t^{-3/2}),\qquad t\rightarrow+\infty,
\end{equation}
which thereby demonstrates that the classical Jensen--Kato mechanism persists on the lattice in the generic setting.
By contrast, {\it the non-generic threshold theory in dimensions \(d \geq 2\) remains considerably less understood}. 

Motivated by these questions, in the present paper we focus on a
three-dimensional lattice analogue of the Jensen--Kato threshold
theory at all spectral thresholds. More precisely, we consider
$$H=-\Delta+V,\qquad \Delta=\Delta_{\mathbb{Z}^3},$$
where $V$ is a real-valued
sequence $\{V(n)\}_{n\in\mathbb{Z}^3}$ satisfying $|V(n)|
\lesssim
\langle n\rangle^{-\beta}$
for a suitable $\beta>0$. We classify the possible threshold types at all four thresholds $0,4,8,12$, and derive the corresponding expansions of the perturbed resolvent
$
R_V(z)=(H-z)^{-1}$
at each threshold as follows.  

At the elliptic thresholds $\lambda=0$ and $12$, three cases arise: regular, resonance, and eigenvalue. As $|\omega|\to0$,
the corresponding leading resolvent behaviors take the form
\begin{align}\label{eq:elliptic-threshold-expansions}
\begin{aligned}
(\text{regular case})\qquad
R_V(\omega)
&=
\mathbb{A}_0+O(\omega^{1/2}),
&
R_V(12-\omega)
&=
\mathscr{A}_0+O(\omega^{1/2}),
\\
(\text{resonance case})\qquad
R_V(\omega)
&=
\mathbb{A}_1\omega^{-1/2}+O(1),
&
R_V(12-\omega)
&=
\mathscr{A}_1\omega^{-1/2}+O(1),
\\
(\text{eigenvalue case})\qquad
R_V(\omega)
&=
-\mathcal{P}_0\omega^{-1}
+O(\omega^{-1/2}),
&
R_V(12-\omega)
&=
\mathcal{P}_{12}\omega^{-1}
+O(\omega^{-1/2}),
\end{aligned}
\end{align}
where $\mathcal{P}_0$ and $\mathcal{P}_{12}$ denote the corresponding
threshold eigenprojections.

At the hyperbolic thresholds $\lambda=4$ and $8$, the absence of
threshold resonances reduces the classification to the regular and
eigenvalue cases. The corresponding leading behaviors are
\begin{align}\label{eq:hyperbolic-threshold-expansions}
\begin{aligned}
(\text{regular case})\qquad
R_V(4+\omega)
&=
\mathbb{B}_0
+\mathbb{B}_1\omega^{1/2}
+O(\omega),
&
R_V(8-\omega)
&=
\mathscr{B}_0
+\mathscr{B}_1\omega^{1/2}
+O(\omega),
\\
(\text{eigenvalue case})\qquad
R_V(4+\omega)
&=
-\mathcal{P}_4\omega^{-1}
+O(\omega^{-1/2}),
&
R_V(8-\omega)
&=
\mathcal{P}_8\omega^{-1}
+O(\omega^{-1/2}),
\end{aligned}
\end{align}
where $\mathcal{P}_4$ and $\mathcal{P}_8$ are the corresponding
threshold eigenprojections. The precise operator topologies, decay
assumptions on the potential, and remainder estimates are stated in the main resolvent expansion Theorems \ref{Asymptotic expansion theorem_0} and \ref{Asymptotic expansion theorem_4} below.

 Combining the resolvent expansions near all four
thresholds with the spectral representation
\eqref{eq:Stone-formula-intro}, we establish Jensen--Kato type
time-decay estimates for the discrete Schr\"odinger evolution $e^{itH}P_{ac}(H)$ on
$\mathbb{Z}^3$. When all four thresholds are regular, we recover the $t^{-3/2}$ decay estimate \eqref{Regular-Kato-Jensen}. 
In the presence of a threshold obstruction, the corresponding singular behavior of the boundary values of the resolvent leads instead to the slower decay $
O(t^{-1/2}).$
Thus, as in the classical Jensen--Kato theory on $\mathbb{R}^3$, the
long-time dynamics are governed by the precise spectral structure at the threshold. In the lattice setting, however, this principle must be accounted for simultaneously at four distinct critical energies.

\subsection{Main Results}
\ \

\vskip0.1cm
To state our main results, we first introduce the classification of the threshold types.
Recall that
\begin{equation*}
   (\Delta \phi)(n)=(\Delta_{\Z^3} \phi)(n)=\sum\limits_{j=1}^{3}(\phi(n+e_j)+\phi(n-e_j)-2\phi(n)).
\end{equation*}
For any $x=(x_1,x_2,x_3)\in\T^3:=(-\pi,\pi]^3$, denote by $\mcaF$ the following Fourier transform on $\ell^2(\Z^3)$:
\begin{equation}\label{fourier transform}
(\mcaF\phi)(x):=(2\pi)^{-3/2}\sum_{n\in\Z^3}^{}e^{-in\cdot x}\phi(n), \quad \phi\in\ell^2(\Z^3).
\end{equation}
Obviously, $\mcaF$ is a unitary operator from $\ell^2(\Z^3)$ to $L^2(\T^3)$, and its inverse is given by
\begin{equation*}
    (\mcaF^{-1}g)(n):=(2\pi)^{-3/2}\int_{\T^3}e^{ix\cdot n}g(x)dx,\quad g\in L^2(\T^3).
\end{equation*}
It thus yields that $-\Delta$ is unitarily equivalent to the following multiplication operator $\mcaM(x)$:
\begin{equation}\label{mcaM}
    -(\mcaF^{-1}\Delta\mcaF \phi)(x)=\mcaM(x)\phi(x),\quad \mcaM(x)=\sum\limits_{j=1}^{3}(2-2{\rm cos}x_j).
\end{equation}
This implies that spectrum of $-\Delta$ is purely absolutely continuous, i.e.,
$$\sigma(-\Delta)=\sigma_{ac}(-\Delta)=[0,12].$$
The critical values 0, 4, 8  and 12 of $\mcaM(x)$ are called the {\it thresholds} of $-\Delta$. Moreover, we also say that $0$ and $12$ are \textbf{elliptic thresholds} and $4$ and $8$
are \textbf{hyperbolic thresholds}. 

Note that unlike the continuous analogue $-\Delta_{\R^3}$ on $\R^3$, where 0 is the only threshold, our discrete setting involves three additional thresholds: 4, 8 and 12, which complicate the classification of resonances. To illustrate this, let $s\in\R$ and denote by
$\mcaB(s,-s)$ the set of all bounded linear operators from $\ell^{2,s}(\Z^3)$ to $\ell^{2,-s}(\Z^3)$, where 
$$\ell^{2,s}(\Z^3):=\big\{\{\phi(n)\}_{n\in\Z^3}\mid\sum\limits_{n\in\Z^3}\left<n\right>^{2s}|\phi(n)|^2<\infty\big\}.$$
\begin{definition}\label{defin of resonance types 0,12}{\rm{\bf (Resonant types at elliptic thresholds).}}
 Let $H=-\Delta+V$ on $\Z^3$ and $|V(n)|\lesssim \left<n\right>^{-\beta}$ for some $\beta>2$. Then, for $\lambda \in \{0,12\}$, we have 
 \vskip0.2cm
\begin{itemize}
\item[{\rm(i)}] $\lambda$ is \textbf{a regular point of $H$} if $H\phi=\lambda\phi$ has only the zero solution in $\bigcap\limits_{s>1/2} \ell^{2, -s}(\Z^3)$.
\vskip0.15cm
\item [{\rm(ii)}] $\lambda$ is \textbf{a resonance of $H$} if  $H\phi=\lambda\phi$ has a non-zero solution in $\bigcap\limits_{s>1/2} \ell^{2, -s}(\Z^3)$ but only the zero solution in $\ell^{2}(\Z^3)$.
\vskip0.15cm
\item [{\rm(iii)}] $\lambda$ is \textbf{an eigenvalue of $H$} if  $H\phi=\lambda\phi$ has a non-zero solution in $\ell^2(\Z^3)$.
\end{itemize}
\end{definition}
\vskip0.1cm
Before giving the resonance types at the hyperbolic thresholds, we recall from Tadano and Taira \cite[Theorem 1.8]{TT19} that, for any $\lambda\in [0,12]$ and $s>1$,  
\begin{equation*}
    (-\Delta-\lambda\mp i0)^{-1}:=\lim\limits_{\varepsilon>0,\ \varepsilon\rightarrow0}(-\Delta-\lambda\mp i\varepsilon)^{-1}\ \ \text{in the operator norm of}\ \ \mcaB(s,-s).
\end{equation*}
Combining this with the weighted resolvent estimate (cf. \cite[Proposition 3.4]{TT19}) 
\begin{equation*}
\sup\limits_{z\in\C\setminus\R}\|(-\Delta-z)^{-1}\|_{\mcaB(1,-1)}<\infty, 
\end{equation*}
and using a standard density argument, we can conclude that the limiting operator $(-\Delta-\lambda\mp i0)^{-1}$ extends uniquely as an element of $\mcaB(1,-1)$ and satisfies
\begin{equation}\label{proper of birmann operator}
    (-\Delta-\lambda)(-\Delta-\lambda\mp i0)^{-1}f=f,\qquad f\in\ell^{2,1}(\Z^3).
\end{equation}
\begin{definition}\label{defin of resonance types 4,8}{\rm{\bf (Resonant types at hyperbolic thresholds).}}
 Let $H=-\Delta+V$ on $\Z^3$ and $|V(n)|\lesssim \left<n\right>^{-\beta}$ for some $\beta>2$. Then, for $\lambda \in \{4,8\}$, we have 
 \vskip0.2cm
\begin{itemize}
\item[{\rm(i)}] $\lambda$ is \textbf{a regular point of $H$} if $\phi+(-\Delta-\lambda\mp i0)^{-1}V\phi=0$ has only the zero solution in $\ell^{2,-1}(\Z^3)$.
\vskip0.15cm
\item [{\rm(ii)}] $\lambda$ is \textbf{a resonance of $H$} if $\phi+(-\Delta-\lambda\mp i0)^{-1}V\phi=0$ has a non-zero solution in $\ell^{2,-1}(\Z^3)$ but only the zero solution in $\ell^{2}(\Z^3)$.
\vskip0.15cm
\item [{\rm(iii)}] $\lambda$ is \textbf{an eigenvalue of $H$} if $\phi+(-\Delta-\lambda\mp i0)^{-1}V\phi=0$ has a non-zero solution in $\ell^2(\Z^3)$.
\end{itemize}
\end{definition}
\begin{remark}\label{remark on defini of resonances}
\ \ 

{\rm (1) It is clear that when $V\equiv0$, all thresholds $0,4,8,12$ are regular points of $H$. 

(2) Note that Nomura and Taira \cite[Theorem 1.7]{NT20} proved that, under the assumption $\beta>\frac{7}{2}$, any $\phi\in\ell^{2,-1}(\Z^3)$ satisfying $\phi+(-\Delta-\lambda\mp i0)^{-1}V\phi=0$ for $\lambda=4$ and $8$, necessarily belongs to $\phi\in\ell^2(\Z^3)$. Hence, no resonance occurs at the hyperbolic thresholds $4$ and $8$. 

(3) Moreover, by \eqref{proper of birmann operator}, for $\lambda\in\{4,8\}$, we have 
$$
\{\phi\in\ell^{2,-1}(\Z^3)\mid\phi+(-\Delta-\lambda\mp i0)^{-1}V\phi=0\}\subseteq\{\phi\in\ell^{2,-1}(\Z^3)\mid H\phi=\lambda\phi\}.
$$
However, the reverse inclusion does not hold in general.  Nevertheless, under suitable decay assumptions on $V$, they are equivalent when $\phi\in\ell^2(\mathbb Z^3)$. In this case, the eigenvalue condition in Definition \ref{defin of resonance types 4,8} can equivalently be replaced by the usual eigenvalue equation. We refer to Appendix \ref{example} for more details.
}
\end{remark}

Now we present the expansions of $R_V(\omega):=(H-\omega)^{-1}$. In this paper, unless otherwise noted, we choose the branch of
$\sqrt{\omega}$ for $\omega \in \mathbb{C}\setminus [0,\infty)$ with
$\operatorname{Im}\sqrt{\omega}>0$.

\begin{theorem}[{\bf Resolvent expansions at elliptic thresholds}]\label{Asymptotic expansion theorem_0}
 Let $H=-\Delta+V$ with $|V(n)|\lesssim \left<n\right>^{-\beta}$ for some $\beta>0$. 
Then the following asymptotic expansions hold in $\mathcal{B}(\sigma,-\sigma)$ as $|\omega|\rightarrow 0$ with $\arg\omega\in(0,2\pi)$. 
\begin{itemize}
\item[{\rm(i)}] If $0$ {\rm(}resp. $12${\rm)} is a regular point of $H$ and $\beta>3$, then for any $\sigma>\frac{3}{2}$, 
\begin{equation}\label{(t)asy expan of regular 0}
R_V(\omega)=\mathbb{A}_0+O(\omega^{\frac{1}{2}}),\qquad
R_V(12-\omega)=\mathscr{A}_0+O(\omega^{\frac{1}{2}}).
\end{equation}
\item [{\rm(ii)}] If $0$ {\rm(}resp. $12${\rm)} is a resonance of $H$ and $\beta>5$, then for any $\sigma>\frac{5}{2}$, 
\begin{equation}\label{(t)asy expan of resonance 0}
R_V(\omega)=\mathbb{A}_1\omega^{-\frac{1}{2}}+O(1),\qquad 
R_V(12-\omega)=\mathscr{A}_1\omega^{-\frac{1}{2}}+O(1).
\end{equation}
\item [{\rm(iii)}] If $0$ {\rm(}resp. $12${\rm)} is an eigenvalue of $H$ and $\beta>7$, then for any $\sigma>\frac{7}{2}$, 
\begin{equation}\label{(t)asy expan of eigenvalue 0}
R_V(\omega)=-\mcaP_0\omega^{-1}+O(\omega^{-\frac{1}{2}}),\qquad R_V(12-\omega)=\mcaP_{12}\omega^{-1}+O(\omega^{-\frac{1}{2}}).
\end{equation}
\end{itemize}
Here $\mcaP_0$ and $\mcaP_{12}$ are the orthogonal projections onto the eigenspaces of $H$ at $0$ and $12$, respectively. Moreover, these expansions are twice differentiable in $\mcaB(\sigma,-\sigma)$ with $\sigma>\frac{5}{2}$.
\end{theorem}
For convenience, we write $R_V^{\pm}(4+\omega):=R_V(4+\omega)$ and $R_V^{\pm}(8-\omega):=R_V(8-\omega)$ when $\pm\operatorname{Im}\omega>0.$
\begin{theorem}[{\bf Resolvent expansions at hyperbolic thresholds}]\label{Asymptotic expansion theorem_4}
 Let $H=-\Delta+V$ with $|V(n)|\lesssim \left<n\right>^{-\beta}$ for some $\beta>0$. 
Then the following asymptotic expansions hold in $\mathcal{B}(\sigma,-\sigma)$  as $|\omega|\rightarrow0$ with $\pm\operatorname{Im}\omega>0$.  
\begin{itemize}
\item[{\rm(i)}] If $4$ {\rm(}resp. $8${\rm)} is a regular point of $H$ and $\beta>7$, then for any $\sigma>\frac{7}{2}$,
\begin{equation}\label{(t)asy expan of regular 4}
R_V^{\pm}(4+\omega)=\mathbb{B}^{\pm}_0+\mathbb{B}^{\pm}_1\omega^{\frac{1}{2}}+O(\omega),\qquad
R_V^{\pm}(8-\omega)=\mathscr{B}^{\pm}_0+\mathscr{B}^{\pm}_1\omega^{\frac{1}{2}}+O(\omega).
\end{equation}
\item [{\rm(ii)}] If $4$ {\rm(}resp. $8${\rm)} is an eigenvalue of $H$ and $\beta>11$, then for any $\sigma>\frac{11}{2}$, 
\begin{equation}\label{(t)asy expan eigenvalue 4}
R_V^{\pm}(4+\omega)=-\mcaP_4\omega^{-1}+O(\omega^{-\frac{1}{2}}),\qquad R_V^{\pm}(8-\omega)=\mcaP_8\omega^{-1}+O(\omega^{-\frac{1}{2}}).
\end{equation}
\end{itemize}
Here $\mcaP_4$ and $\mcaP_{8}$ are are the orthogonal projections onto the eigenspaces of $H$ at $4$ and $8$, respectively. Moreover, these expansions are twice differentiable in $\mcaB(\sigma,-\sigma)$ with $\sigma>\frac{11}{2}$.
\end{theorem}
The above results extend Kopylova's result \cite{Kop09}~(which treats only regular thresholds), to nonregular threshold cases, including resonances and eigenvalues. The detailed proofs of Theorems \ref{Asymptotic expansion theorem_0}$\sim$\ref{Asymptotic expansion theorem_4} above will be provided in Section \ref{proof of Asy theorem}.

As an application of the resolvent expansions, we  establish the weighted  $\ell^2$ time-decay estimates for solutions to the following Schr\"{o}dinger equation on $\Z^3$:
\begin{equation*}
\left\{
\begin{aligned}
&i\partial_t u=(-\Delta+V)u, \\
&u(0, n) = u_0(n),
\end{aligned}
\right. \quad (t, n) \in {\R} \times {\Z^3}.
\end{equation*}
\begin{theorem}\label{decay theorem}
Let $H=-\Delta+V$ with $|V(n)|\lesssim \left<n\right>^{-\beta}$ for $\beta>11$. Suppose that $H$ has no positive eigenvalues in the interval $\mcaI=(0,4)\cup(4,8)\cup(8,12)$, and let $P_{ac}(H)$ denote the spectral projection onto the absolutely continuous spectrum space of $H$. Then the following decay estimates hold in $\mcaB(\sigma,-\sigma)$ for any $\sigma>\frac{11}{2}$:
\begin{equation*}
    e^{-itH}P_{ac}(H)=\begin{cases}
        O( t^{-\frac{3}{2}}),& all\ 0,4,8,12 \  are\ regular\ points \ of \ H\\
        O( t^{-\frac{1}{2}}),&otherwise
    \end{cases},\quad t\rightarrow+\infty.
\end{equation*}
\end{theorem}
This result extends Kopylova's result \cite{Kop09} (where only the case in which all thresholds are regular is considered) to nonregular threshold cases, with the corresponding decay rate $O(t^{-1/2})$.

\begin{remark}{\rm (Embedded eigenvalue problem)}
    {\rm The absence of embedded eigenvalues is a natural spectral assumption in the above time-decay estimates, since such eigenvalues lead to singularities of the resolvent inside the continuous spectrum. 
    
 For discrete Schrödinger operators, this issue has drawn considerable interest and has been studied under stronger decay assumptions on the potential (see e.g. \cite{HHNO16,HMK22,HSSS12,IM14,Liu22,LMT25} and references therein). In particular, Isozaki and Morioka \cite{IM14} proved the absence of embedded eigenvalues for compactly supported potentials on $\Z^d$. More recently, Liu, Matos and Treuer \cite{LMT25} extended this conclusion to potentials decaying faster than any exponential, namely,
    $$
\limsup_{|n|\to\infty}\frac{\ln |V(n)|}{|n|}<-\gamma
\qquad\text{for every }\gamma>0.
$$
These results, however, do not directly cover the polynomial decay assumption imposed in this theorem. Thus, the absence of embedded eigenvalues remains an interesting problem for further investigation. 
 }
\end{remark}

\subsection{Organization of the Paper}
\ \

\vskip0.1cm
The remainder of this paper is organized as follows. In Section \ref{proof of Asy theorem}, we prove the resolvent expansions stated in Theorems \ref{Asymptotic expansion theorem_0} and \ref{Asymptotic expansion theorem_4}, building on a refined analysis of the kernel coefficients in Kopylova's free resolvent expansion given in Appendix \ref{further infor of coeffi}. In Section \ref{sec of application}, we apply these expansions to establish Jensen--Kato type time-decay estimates for the discrete Schr\"odinger evolution on \(\mathbb Z^3\). The argument relies on a limiting absorption principle with spectral derivatives and on the spectral representation of $e^{-itH}P_{ac}(H)$, both of which are established in Appendix \ref{LAP and repre of eitHpacH}. Finally, Appendix \ref{example} clarifies the relation between the threshold equation and the Birman--Schwinger equation at the hyperbolic thresholds.

\section{Proof of Resolvent Expansions} \label{proof of Asy theorem}
In this section, we prove the asymptotic expansions of the perturbed resolvent $R_V(\omega)$ near the spectral thresholds $\omega= 0, 4, 8, 12$, namely, Theorems \ref{Asymptotic expansion theorem_0} and \ref{Asymptotic expansion theorem_4}. 

Following the reduction scheme of Jensen and Nenciu \cite{JN01}, we analyze $R_V(\omega)$ through an auxiliary operator constructed from the free resolvent. More precisely, let $R_0(\omega)=(-\Delta-\omega)^{-1}$ denote the free resolvent, and define \begin{equation}\label{def of Mw}
M(\omega):=U+vR_0(\omega)v,\quad U(n)=\left\{\begin{aligned}
&1,\quad\ \  V(n)\geq0,\\
&-1, \ V(n)<0,
\end{aligned}\right.\qquad v(n)=|V(n)|^{\frac{1}{2}}.
 \end{equation} 
Using $U^2=I$ (the identity operator) and the resolvent identity 
\begin{equation}\label{R0-Rv}
    R_0(\omega)-R_V(\omega)=R_0(\omega)VR_V(\omega)=R_V(\omega)VR_0(\omega),
\end{equation}
one readily verifies that $M(\omega)$ is invertible on $\ell^2(\Z^3)$ and that
$$M^{-1}(\omega)=U-UvR_V(\omega)vU.$$
Moreover, by \eqref{R0-Rv} and $V=vUv$, we have
\begin{equation}\label{rela betw RV and inver of Mw}
R_V(\omega)=R_0(\omega)-R_0(\omega)vM^{-1}(\omega)vR_0(\omega).
\end{equation}
Therefore, to obtain the asymptotic expansions of $R_V(\omega)$ near the spectral thresholds, it suffices to analyze $R_0(\omega)$ and $M^{-1}(\omega)$.
\subsection{Free Resolvent Expansions }\label{FREE}
\ \

\vskip0.1cm
We begin with the free resolvent $R_0(\omega)$. Although its operator-valued expansions near the thresholds $0$ and $4$ were established in \cite[Propositions 2.4 and 2.6]{Kop09}, the analysis of $M^{-1}(\omega)$ requires finer information on the structure and form of certain coefficients. We therefore supplement these expansions by computing the relevant coefficients explicitly.
We denote by $K$ the operator with kernel $K(n,m)$:
 $$(Kf)(n):=\sum\limits_{m\in\Z^3}^{}K(n,m)f(m).$$ 
\begin{lemma}\label{lemma of R0 0}
Let $ N = 0, 1, 2\cdots$ and $\sigma > N +\frac{3}{2}.$ Then the following expansion holds in $\mathcal{B}(\sigma, -\sigma):$
\begin{equation}\label{expan of R0}
\ R_{0}(\omega) =\sum\limits_{k=0}^{N}A_k \omega^{k/2}+ O(\omega^{{(N+1)}/{2}}),\quad \left| \omega\right|\rightarrow0,\quad \arg\omega\in(0,2\pi).
\end{equation}
Here $A_{k}\in\ \mathcal{B}(s, -s)$ with $ s>k+\frac{1}{2}.$ Moreover, if $k$ is odd, the kernel $A_k(n-m)$ of $A_k$ is given by
\begin{equation*}
 A_k(n)=\frac{i(\frac{k-1}{2})!}{4\pi k!}\sum\limits_{\substack{k_1+k_2+k_3=\frac{k-1}{2}\\  k_1,k_2,k_3\in \N}} \prod\limits_{j=1}^{3}\frac{(\frac{1}{2}-n_j)_{k_j}(\frac{1}{2}+n_j)_{k_j}}{k_j!},\quad n=(n_1,n_2,n_3)\in\Z^3
\end{equation*}
with $(a)_{\ell}:=a(a+1)\cdots(a+\ell-1)$ for any integer $\ell\geq1$ and $ (a)_0=1$. In particular, $A_1(n)=\frac{i}{4\pi}$.
\end{lemma}
For convenience, we write $R_0^{\pm}(4+\omega):=R_0(4+\omega)$ and $R_0^{\pm}(8-\omega):=R_0(8-\omega)$ when $\pm\operatorname{Im}\omega>0.$
\begin{lemma}\label{lemma of free reso expan at 4}
{\rm(\cite[Proposition 2.6]{Kop09})}
Let $ N =-1, 0, 1,\cdots$ and $\sigma >2 N +\frac{7}{2}.$ Then the following expansion holds in $\mathcal{B}(\sigma, -\sigma):$
\begin{equation}\label{expan at 4}
R^{\pm}_{0}(4+\omega) =\sum\limits_{k=0}^{N}E^{\pm}_k\omega^{k}+ \sum\limits_{k=0}^{N}B^{\pm}_k\omega^{k+1/2}+O(\omega^{N+1}),\quad \left| \omega \right|\rightarrow0,\quad \pm\operatorname{Im} \omega>0,
\end{equation}
where the operators $E^{\pm}_{k},B^{\pm}_{k}\in \mathcal{B}(s, -s)$ with $ s>2k+\frac{3}{2}.$ Moreover, $E^-_k=(E^+_k)^*$ and $B^-_k=-(B^+_k)^*$. 
\end{lemma}

\begin{remark}
{\rm 
As for the other two thresholds $8$ and $12$, by means of the unitary transform $$(\mcaU\phi)(n)=(-1)^{n_1+n_2+n_3}\phi(n), \quad n=(n_1,n_2,n_3),$$ we have $R_0(12-4k-\omega)=-\mcaU R_0(4k+\omega)\mcaU$ for $k=0,1$, and therefore, 
\begin{align}
R_{0}(12-\omega) &=\sum\limits_{k=0}^{N}\mcaA_k \omega^{k/2}+ O(\omega^{{(N+1)}/{2}}),\quad \left| \omega\right|\rightarrow0,\quad \arg\omega\in(0,2\pi),\label{expan of R012}\\
R^{\pm}_{0}(8-\omega) &=\sum\limits_{k=0}^{N}\mcaE^{\pm}_k\omega^{k}+ \sum\limits_{k=0}^{N}\mcaB^{\pm}_k\omega^{k+1/2}+O(\omega^{N+1}),\quad \left| \omega \right|\rightarrow0,\quad \pm \operatorname{Im} \omega>0,\label{expan of R08}
\end{align}
where $\mcaA_{k}=-\mcaU A_k\mcaU,\ \mcaE^{\pm}_{k}=-\mcaU E^{\pm}_k\mcaU$ and $ \mcaB^{\pm}_{k}=-\mcaU B^{\pm}_k\mcaU.$
}
\end{remark}

\begin{remark}\label{remark of free expan}\ \

{\rm {\bf(1)} At the threshold $0$, Kopylova \cite[Proposition 2.4]{Kop09} established the operator-valued expansion \eqref{expan of R0}. In her argument, the discrete symbol is locally replaced by its quadratic approximation $|\theta|^2$, which captures the structure of the expansion but does not yield explicit kernel formulas for the coefficients associated with the exact discrete symbol $$\mcaM(\theta)=\sum\limits_{j=1}^{3}(2-2\cos\theta_j),\qquad \theta=(\theta_1,\theta_2,\theta_3)\in (-\pi,\pi]^3.$$ Since such kernel information is required in our subsequent analysis, we rederive the expansion using the exact discrete symbol and compute the odd-order kernels $A_{2k+1}(n)$ explicitly; see Proposition \ref{prop of Ak calcul}.
\vskip0.1cm
{\bf(2)} At the threshold $4$, we use the expansion in \cite[Proposition 2.6]{Kop09}, which is derived from the exact discrete symbol. The expansion there is stated for $\operatorname{Im}\omega>0$, while the corresponding expansion for $\operatorname{Im}\omega<0$ follows from the identity $R_0(4+\omega)=R_0^*(4+\overline{\omega}).$
The argument in \cite{Kop09} analyzes the contributions from the critical points, while the contribution from the noncritical region is not made explicit. For our purposes, we verify that the latter contributes only integer powers of $\omega$. It follows that the half-integer coefficients $B^{\pm}_k(n)$ are determined entirely by the three critical points. In particular, we need the explicit form of $B^{\pm}_0(n)$, which we compute as follows (see Proposition \ref{prop of B0 calcu}):
\begin{equation}\label{kernel B0pmn}
 B^{\pm}_0(n)=\frac{\mp}{4\pi}\big((-1)^{n_1}+(-1)^{n_2}+(-1)^{n_3}\big),\qquad n=(n_1,n_2,n_3)\in\Z^3.
\end{equation}

{\bf(3)} Moreover, the coefficients $A_0$ and $E_0^\pm$ can be identified with fundamental solutions of $-\Delta$ and $-\Delta-4$, respectively. More precisely, $(-\Delta )A_0=\delta_0$ and $(-\Delta-4)E^{\pm}_0=\delta_0$, where $\delta_0$ denotes the Kronecker delta. Their kernels are commonly referred to as the corresponding {\it Green functions}. Formally, the Green function at energy $\lambda$ admits the Fourier representation
\begin{equation*}
    G(\lambda,n)=(2\pi)^{-3}\int_{(-\pi,\pi]^3}\frac{e^{in\cdot\theta}}{\mcaM(\theta)-\lambda}d\theta.
\end{equation*}

 In the Euclidean setting, such Green function information plays an important role in relating the inversion scheme for $M^{-1}(\omega)$ to the characterization of threshold eigenfunctions and resonances (see, e.g., \cite{JN01}). In higher-dimensional lattice settings, however, no comparably simple general kernel formulas are available, except in the one-dimensional case (cf. \cite{KKK06,IJ21}). This lack of explicit kernel information makes it difficult to relate the threshold subspaces obtained from the inversion scheme to actual threshold eigenfunctions and resonances.

In our case, the situation at the threshold $\lambda=0$ is more favorable. Since $(\mcaM(\theta))^{-1}\sim|\theta|^{-2}$ near $\theta=0$ is integrable on $(-\pi,\pi]^3$, the Green function $G(0,n)$ is well defined. Moreover, Martinsson and Rodin \cite[Example 1.1]{MR02} obtained the asymptotic relation
\begin{equation}\label{relation of phi0 and g0}
    |G(0,n)-G_d(n)|\lesssim \frac{1}{|n|^3},\quad G_d(n)=\frac{1}{4\pi|n|},\quad n\neq0,
\end{equation}
where $G_d(n)$ is the continuum Green function of $-\Delta_{\R^3}$ evaluated at lattice points. Since $A_0\in\mcaB(s,-s)$ for all $s>1/2$, we have $A_0(\cdot)\in\ell^{2,-s}(\Z^3)$ for all $s>1/2$. By \eqref{relation of phi0 and g0}, the same holds for $G(0,\cdot)$. Moreover, both are fundamental solutions of $-\Delta$. Therefore, we conclude that $A_0(n)=G(0,n)
$. We remark that this identity, together with \eqref{relation of phi0 and g0}, plays an important role in the subsequent analysis of $M^{-1}(\omega)$ and in the characterization of threshold resonances at $\lambda=0$ (see Lemma \ref{chara of solution and space}).

The situation at the threshold $\lambda=4$ is different. No analogous asymptotic description of the corresponding Green function is available for our purposes. However, the coefficients $E_0^\pm$ can be identified with the boundary values of the free resolvent:
\begin{equation*}
    E^{\pm}_0=(-\Delta-4\mp i0)^{-1}\quad \text {on}\ \ \ell^{2,s}(\Z^3)\ \ \text {for}\ \  s>7/2.
\end{equation*}
This identification allows us to make use of the threshold analysis of Nomura and Taira \cite{NT20}, in particular their absence-of-resonance result at hyperbolic thresholds. Together with the structure obtained from the inversion scheme, their result provides a key ingredient in the characterization of the threshold space at $\lambda=4$ (see Lemma \ref{chara of solution and space 4}).
}    
\end{remark}
\subsection{Resolvent Expansions at Elliptic Thresholds \texorpdfstring{$0$}{0} and \texorpdfstring{$12$}{12}}\label{subsec of expan at 0}
\ \

\vskip0.1cm
To state our result, we first introduce some notation and auxiliary operators. 
Let $A_k$ and $\mcaA_k$ be defined as in \eqref{expan of R0} and \eqref{expan of R012}, respectively, and set $$S:=U+vA_0v,\qquad \mcaS:=U+v\mcaA_0v.$$
Throughout this paper, we use the inner product
$$\left<f,g\right>:=\sum\limits_{m\in\Z^3}^{}f(m)\overline{g(m)},\quad f,g\in\ell^2(\Z^3).$$
Let $S_1$ be the orthogonal projection onto $\operatorname{Ker}S$, and let $S_2$ be the orthogonal projection onto the kernel of $S_1vA_1vS_1$ restricted to $S_1\ell^2(\Z^3)$. Since $A_1(n)=\frac{i}{4\pi}$, we have
\begin{equation*}
    S_1\ell^2(\Z^3)=\{f\in \ell^2(\Z^3)~|~ Sf=0\},\qquad S_2\ell^2(\Z^3)=\{f\in S_1\ell^2(\Z^3)~|\left<f,v\right>=0\}.
\end{equation*}
Similarly, let $\mcaS_1$ be the orthogonal projection onto $\operatorname{Ker}\mcaS$, and let $\mcaS_2$ be the orthogonal projection onto the kernel of $\mcaS_1v\mcaA_1v\mcaS_1$ restricted to $\mcaS_1\ell^2(\Z^3)$. Then 
\begin{equation*}
\mcaS_1\ell^2(\Z^3)=\{f\in \ell^2(\Z^3)~|~ \mcaS f=0\},\qquad \mcaS_2\ell^2(\Z^3)=\{f\in \mcaS_1\ell^2(\Z^3)~|\left<f,\mcaU v\right>=0\}.
\end{equation*}

\begin{proposition}\label{Asymptotic expansion theorem 0}
 Let $H=-\Delta+V$ with $|V(n)|\lesssim \left<n\right>^{-\beta}$ for some $\beta>0$. Then the following asymptotic expansions hold in $\mcaB(\sigma,-\sigma)$ as $|\omega|\rightarrow 0$ with $\arg\omega\in(0,2\pi)$.  
\begin{itemize}
\item[{\rm(1)}] If $0$ {\rm(}resp. $12${\rm)} is a regular point of $H$ and $\beta>3$, then for any $\sigma>\frac{3}{2}$, 
\begin{equation*}
R_V(\omega)=A_0-A_0vD_{0}vA_0+O(\omega^{\frac{1}{2}}),\quad
R_V(12-\omega)=\mcaA_0-\mcaA_0v\mcaD_{0}v\mcaA_0+O(\omega^{\frac{1}{2}}).
\end{equation*}
\item [{\rm(2)}] If $0$ {\rm(}resp. $12${\rm)} is a resonance of $H$ and $\beta>5$, then for any $\sigma>\frac{5}{2}$, 
\begin{equation*}
R_V(\omega)=-A_0vS_1D_{1}S_1vA_0\omega^{-\frac{1}{2}}+O(1),\quad 
R_V(12-\omega)=-\mcaA_0v\mcaS_1\mcaD_{1}\mcaS_1v\mcaA_0\omega^{-\frac{1}{2}}+O(1).
\end{equation*}
\item [{\rm(3)}] If $0$ {\rm(}resp. $12${\rm)} is an eigenvalue of $H$ and $\beta>7$, then for any $\sigma>\frac{7}{2}$, 
\begin{equation*}
R_V(\omega)=-A_0vS_2D_{2}S_2vA_0\omega^{-1}+O(\omega^{-\frac{1}{2}}),\quad R_V(12-\omega)=-\mcaA_0v\mcaS_2\mcaD_{2}\mcaS_2v\mcaA_0\omega^{-1}+O(\omega^{-\frac{1}{2}}).
\end{equation*}
\end{itemize}
Here
\begin{itemize}
\item $D_0=(S+S_1)^{-1},\quad \mcaD_0=(\mcaS+\mcaS_1)^{-1};$
    \vskip0.1cm
    \item $D_{1}=(S_1vA_1vS_1+S_2)^{-1}\ {\rm on} \ S_1\ell^2(\Z^3),\quad \mcaD_{1}=(\mcaS_1v\mcaA_1v\mcaS_1+\mcaS_2)^{-1}\ {\rm on} \ \mcaS_1\ell^2(\Z^3);$
    \vskip0.1cm
    \item $ D_{2}=(S_2vA_2vS_2)^{-1}\ {\rm on} \ S_2\ell^2(\Z^3),\quad\qquad\  \mcaD_{2}=(\mcaS_2v\mcaA_2v\mcaS_2)^{-1}\ {\rm on} \ \mcaS_2\ell^2(\Z^3),$
\end{itemize}
and $A_0vS_2D_{2}S_2vA_0$ and $-\mcaA_0v\mcaS_2\mcaD_{2}\mcaS_2v\mcaA_0$  
 are the orthogonal projections onto the eigenspaces of $H$ at $0$ and $12$, respectively. Moreover, these expansions are twice differentiable in $\mcaB(\sigma,-\sigma)$ with $\sigma>\frac{5}{2}$.
\end{proposition}
Having established the required free resolvent expansions, it remains to analyze $M^{-1}(\omega)$. Since the free resolvent expansions at the thresholds 0 and 12 have the same structure, we only treat the threshold 0; the argument at 12 is analogous. The expansion of $M^{-1}(\omega)$ is obtained by successive inversion. We first present two lemmas. The first connects the inversion scheme with the spectral characterization of the threshold.

\begin{lemma}\label{chara of solution and space}
{\rm Let $H=-\Delta+V$ with $|V(n)|\lesssim\left<n\right>^{-\beta}$ for $\beta>7$. Then
\begin{itemize}
\item [{\rm(i)}]$f\in S_1\ell^2(\Z^3)\Longleftrightarrow \exists\ \phi\in \bigcap\limits_{s>1/2} \ell^{2, -s}(\Z^3)$ such that $H\phi=0$. Moreover, $f=Uv\phi$ and $\phi=-A_0vf$.
\vskip0.15cm
\item [{\rm(ii)}]$f\in S_2\ell^2(\Z^3)\Longleftrightarrow \exists\ \phi\in \ell^2(\Z^3)$ such that $H\phi=0$. Moreover, $f=Uv\phi$ and $\phi=-A_0vf$.
\end{itemize}
}
\end{lemma}
\begin{remark}
{\rm Note that when $\beta>1$, both $UvA_0v$ and $Uv\mcaA_0v$ are compact operators on $\ell^2(\Z^3)$. Since  $U$ is invertible, it follows that 
    $US=I+UvA_0v$ and  $U\mcaS=I+Uv\mcaA_0v$
are Fredholm of index zero. Hence $$\dim(S_2\ell^2(\Z^3))\leq\dim(S_1\ell^2(\Z^3))<+\infty,\qquad\dim(\mcaS_2\ell^2(\Z^3))\leq\dim(\mcaS_1\ell^2(\Z^3))<+\infty.$$
Combining Definition \ref{defin of resonance types 0,12} with this lemma, we obtain the following characterization:
\begin{itemize}
\item $0$ is a regular point of $H$ $\Longleftrightarrow$
    $S_1\ell^2(\Z^3)=\{0\}$$\Longleftrightarrow$$S$ is invertible on $\ell^2(\Z^3)$.
\vskip0.1cm
\item 
    $0$ is a resonance of $H$ $\Longleftrightarrow$ $S_1\ell^2(\Z^3)\neq\{0\}$ but $S_2\ell^2(\Z^3)=\{0\}$$\Longleftrightarrow$ $S$ is not invertible on $\ell^2(\Z^3)$ but $S_1vA_1vS_1$ is invertible on $S_1\ell^2(\Z^3)$.
\vskip0.1cm
\item $0$ is an eigenvalue of $H$ $\Longleftrightarrow$ $S_2\ell^2(\Z^3)\neq\{0\}$$\Longleftrightarrow$ $S_1vA_1vS_1$ is not invertible on $S_1\ell^2(\Z^3)$. 
\end{itemize}
Thus, the successive invertibility conditions arising in the inversion
scheme precisely distinguish the regular, resonant, and eigenvalue cases
at the threshold.
}
\end{remark}
The proof of this lemma will be postponed to the end of this subsection. Moreover, the following lemma will be frequently used.
\begin{lemma}\label{lemm of expa}
{\rm (\cite[Lemma 2.1]{JN01})} Let $\mscH$ be a complex Hilbert space. Let $A$ be a closed operator and $P$ a projection. Suppose $A+P$ has a bounded inverse. Then $A$ has a bounded inverse if and only if
$$a\equiv P-P(A+P)^{-1}P$$
has a bounded inverse in $P\mscH$, and in this case
$$A^{-1}=(A+P)^{-1}+(A+P)^{-1}Pa^{-1}P(A+P)^{-1}.$$
\end{lemma}
\begin{remark}\label{remark of invertibility}
{\rm In particular, if $A$ is bounded on $\mcaH$ and has closed range, $P$ is the orthogonal projection onto $\operatorname{Ker}A$, and $\operatorname{Ker}A^*=\operatorname{Ker}A$, then $A+P$ is invertible on $\mcaH$. Moreover, we can obtain
$$(A+P)^{-1}P=P=P(A+P)^{-1}.$$
This criterion will be used repeatedly in the expansion of $M^{-1}(\omega)$. 
    }
\end{remark}
\vskip0.3cm
\begin{proof}[Proof of Proposition {\rm\ref{Asymptotic expansion theorem 0}}]
We first show that all the operators $D_j$ are well-defined. Since
$A_0(n)=G(0,n)$ is real-valued and $A_1(n)=\frac{i}{4\pi}$, we have $S^*=S$ and $(S_1vA_1vS_1)^*=-S_1vA_1vS_1$. Moreover, $S$ is Fredholm, and $S_1$ is of finite rank. Hence
$\operatorname{Ran}S$ is closed, while
$\operatorname{Ran}(S_1vA_1vS_1)$ is finite-dimensional and therefore
closed. It then follows from Remark \ref{remark of invertibility} that
$D_0$ and $D_1$ are well-defined and satisfy
\begin{equation*}
D_0S_1=S_1=S_1D_0,\qquad D_1S_2=S_2=S_2D_1.
\end{equation*}
It remains to consider $D_2$. Since $S_2\ell^2(\Z^3)$ is
finite-dimensional, it suffices to prove that
$S_2vA_2vS_2$ is injective on $S_2\ell^2(\Z^3)$. Suppose that
$f\in S_2\ell^2(\Z^3)$ satisfies $S_2vA_2vS_2f=0$. Then
$$0=\left<S_2vA_2vS_2f,f\right>=\left<A_2vf,vf\right>.$$
Since $f\in S_2\ell^2(\Z^3)$, we have
$\langle f,v\rangle=0$, and hence $A_1vf=0$. Taking $N=2$ and
$\omega=i\varepsilon$, $\varepsilon>0$, in Lemma
\ref{lemma of R0 0}, we therefore obtain
\begin{align*}
\left<A_2vf,vf\right> &=\Big<\lim\limits_{\varepsilon\rightarrow0}\frac{R_0(i\varepsilon)-A_0-A_1(i\varepsilon)^{\frac{1}{2}}}{i\varepsilon}vf,vf\Big>=\lim\limits_{\varepsilon\rightarrow0} \Big<\frac{R_0(i\varepsilon)-A_0}{i\varepsilon}vf,vf\Big>\\
&=\lim\limits_{\varepsilon\rightarrow0}\int_{(-\pi,\pi]^3}\frac{1}{(\mathcal{M}(x)-i\varepsilon)\mathcal{M}(x)}|(\mathcal{F}(vf))(x)|^2dx\\
&=\big<\mcaF(A_0vf),\mcaF (A_0vf)\big>\\
&=\big<A_0vf,A_0vf\big>.
\end{align*}
Indeed, since $f\in S_2\ell^2(\mathbb Z^3)$, $\langle f,v\rangle=0$ and Lemma \ref{chara of solution and space} (2) gives $A_0vf\in\ell^{2}(\Z^3)$. By the Fourier representation of $A_0(n)$ and Plancherel's theorem,
$\mathcal{F}(vf)(x)/\mathcal{M}(x)=\mcaF(A_0vf)(x)\in L^2((-\pi,\pi]^3)$. Moreover, 
$$\left| \frac{1} {(\mcaM(x)-i\varepsilon)\mcaM(x)} \right| \leq \frac{1}{|\mcaM(x)|^2},$$
so the last limit follows from the dominated convergence theorem. Consequently, $A_0vf=0$. This combined with $Sf=0$ yields $f=0$. Thus
$S_2vA_2vS_2$ is injective on $S_2\ell^2(\Z^3)$ and hence invertible
there. Therefore $D_2$ is well-defined.

Having established the well-definedness of all the operators $D_j$, we now compute the asymptotic expansion of $M^{-1}(\omega)$ in each threshold type.
\vskip0.2cm
{\bf{(1)}} Suppose that {\bf$\bm 0$ is a regular point of $\bm H$}, that is, $S$ is invertible on $\ell^2(\Z^3)$, and assume that {${\beta>3}$}.
Taking $N=0$ in \eqref{expan of R0}, we obtain the following expansion on $\ell^{2}(\Z^3)$:
 \begin{equation*}
M(\omega)=U+vR_0(\omega)v=S+O(\omega^{\frac{1}{2}}),\quad |\omega|\rightarrow 0,\  \arg\omega\in(0,2\pi).
 \end{equation*}
Combining $D_0=S^{-1}$ and Neumann series expansion, we derive 
\begin{equation}\label{inver of M 0 regular}
\bm {M^{-1}(\omega)=D_{0}+D_0O(\omega^{\frac{1}{2}})D_0.}
\end{equation}
\vskip0.1cm
{\bf{(2)}} Suppose that {\bf$\bm 0$ is a resonance of $\bm H$}, i.e., $S$ is not invertible on $\ell^2(\Z^3)$, while $S_1vA_1vS_1$ is invertible on $S_1\ell^2(\Z^3)$, and assume that $\beta>5$. Taking $N=1$ in \eqref{expan of R0}, we obtain
\begin{equation*}
M(\omega)=S+vA_1v\,\omega^{\frac{1}{2}}+O(\omega),\quad |\omega| \rightarrow0,
\end{equation*}
and hence,
\begin{equation}\label{expre M+S_1}
\left(M(\omega)+S_1\right)^{-1}=D_{0}-D_{0}vA_1vD_{0}\,\omega^{\frac{1}{2}}+D_0O(\omega)D_0,\quad |\omega|\rightarrow0.
\end{equation}
Using $S_1D_{0}=S_1=D_{0}S_1$, we have
 \begin{equation*}
    M_1(\omega):=S_1-S_1\big(M(\omega)+S_1\big)^{-1}S_1=\omega^{\frac{1}{2}}\big(S_1vA_1vS_1+S_1O(\omega^{\frac{1}{2}})S_1\big),\quad |\omega|\rightarrow0.
\end{equation*}
By assumption, $M_1(\omega)$ is invertible on $S_1\ell^2(\Z^3)$ and its inverse admits the following expansion:
\begin{equation}\label{expan of M1 1st}
    M^{-1}_1(\omega)=D_{1}\omega^{-\frac{1}{2}}+D_1O(1)D_1,\quad |\omega|\rightarrow0.
\end{equation}
Applying Lemma \ref{lemm of expa} to
$M(\omega)$, we obtain
\begin{equation}\label{rela inver of M and M1}
    M^{-1}(\omega)=(M(\omega)+S_1)^{-1}\big(I+S_1M^{-1}_1(\omega)S_1(M(\omega)+S_1)^{-1}\big).
\end{equation}
Substituting \eqref{expre M+S_1}  and \eqref{expan of M1 1st} into \eqref{rela inver of M and M1}, we derive
\begin{equation}\label{inver of M 0 resonan}
\bm{M^{-1}(\omega)=S_1D_{1}S_1\omega^{-\frac{1}{2}}+O(1).}
\end{equation}
\vskip0.2cm
{\bf{(3)}} Let {\bf$\bm 0$ be an eigenvalue of $\bm H$} and $\beta>7$. Taking $N=2$ in \eqref{expan of R0} and repeating the process as in (2), we further obtain
\begin{equation}\label{M+S_1 eigenv}
    \left({M}(\omega)+S_1\right)^{-1}=D_{0}-D_{0}vA_1vD_{0}\,\omega^{\frac{1}{2}}+\big((D_{0}vA_1v)^2D_{0}-D_{0}vA_2vD_{0}\big)\omega+D_0O(\omega^{\frac{3}{2}})D_0
\end{equation}
and 
\begin{equation*}
    M_1(\omega)=\omega^{\frac{1}{2}}\Big(S_1vA_1vS_1+\big(S_1vA_2vS_1-S_1vA_1vD_{0}vA_1vS_1\big)\omega^{\frac{1}{2}}+S_1O(\omega)S_1\Big)=:\omega^{\frac{1}{2}}\widetilde{M}_1(\omega).
\end{equation*}
Since $S_2\ell^2(\Z^3)={\rm Ker}S_1vA_1vS_1|_{S_1\ell^2(\Z^3)}$, the operator $\widetilde{M}_1(\omega)+S_2$ is invertible on $S_1\ell^2(\Z^3)$ and 
\begin{equation}\label{expre M_1+S_2}
\left(\widetilde{M}_1(\omega)+S_2\right)^{-1}=D_{1}+D_1(S_1vA_1vD_{0}vA_1vS_1-S_1vA_2vS_1)D_1\omega^{\frac{1}{2}}+D_1O(\omega)D_1,\quad |\omega|\rightarrow0.
\end{equation}
Using $S_2D_{j}=S_2=D_{j}S_2\ (j=0,1),\ S_2S_1=S_2=S_1S_2$ together with $  A_1vS_2=0$, we obtain $$M_2(\omega):=S_2-S_2\big(\widetilde{M}_1(\omega)+S_2\big)^{-1}S_2=\omega^{\frac{1}{2}}\Big(S_2vA_2vS_2+S_2O(\omega^\frac{1}{2})S_2\Big).$$ 
Since $S_2vA_2vS_2$ is invertible on $S_2\ell^2(\Z^3)$. It follows that ${M}_2(\omega)$ is invertible on $S_2\ell^2(\Z^3)$ with
\begin{equation}\label{inverse of M2}
M^{-1}_2(\omega)=S_2D_{2}S_2\,\omega^{-\frac{1}{2}}+S_2O(1)S_2,\quad |\omega|\rightarrow0.  
\end{equation}
Applying Lemma \ref{lemm of expa} to
$\widetilde{M}_1(\omega)$, and using \eqref{expre M_1+S_2} and \eqref{inverse of M2}, we get
\begin{equation}\label{rela inver of tutaM1 and M2}
   \widetilde{ M}^{-1}_1(\omega)=(\widetilde{M}_1(\omega)+S_2)^{-1}\Big(I+S_2M^{-1}_2(\omega)S_2(\widetilde{M}_1(\omega)+S_2)^{-1}\Big)=S_2D_2S_2\,\omega^{-\frac{1}{2}}+D_1O(1).
\end{equation}
Therefore,
\begin{equation*}
M^{-1}_1(\omega)=S_2D_{2}S_2\,\omega^{-1}+D_1O(\omega^{-\frac{1}{2}}),\quad |\omega|\rightarrow0.
\end{equation*}
Combining this with \eqref{M+S_1 eigenv} and \eqref{rela inver of M and M1}, gives 
\begin{equation}\label{inver of M 0 eigen}
\bm{M^{-1}(\omega)=S_2D_{2}S_2\,\omega^{-1}+O(\omega^{-\frac{1}{2}}).}
\end{equation}
Hence, substituting \eqref{inver of M 0 regular}, \eqref{inver of M 0 resonan}, \eqref{inver of M 0 eigen} into \eqref{rela betw RV and inver of Mw}, respectively, we obtain the expansions in Proposition \ref{Asymptotic expansion theorem 0}. Moreover, combining the derivative argument for $R_0(\omega)$ in \cite[Remark 2.5]{Kop09} with the corresponding expansions of $M^{-1}(\omega)$ and the resolvent identity \eqref{rela betw RV and inver of Mw}, the above perturbed expansions can be differentiated twice in $\mcaB(\sigma,-\sigma)$ for $\sigma>\frac52$.
\vskip0.1cm
{\bf(4)} Finally, we prove that $\mcaP_0:=A_0vS_2D_{2}S_2vA_0$ is the projection onto the eigenspace of the eigenvalue $\lambda=0$. 
Denote by $V_0$ the eigenspace of the eigenvalue $0$ of $H$, that is,
$V_0:=\{\phi\in\ell^2(\Z^3):H\phi=0\}.
$ Let $\{f_j\}^{N}_{j=1}$ be an orthonormal basis  of $S_2\ell^2(\Z^3)$. For each $j\in\{1,2,\cdots, N\}$, define $\phi_j:=-A_0vf_j$. By Lemma \ref{chara of solution and space} (ii), $\phi_j\in V_0$ and $f_j=Uv\phi_j$. Moreover, since $\{f_j\}^{N}_{j=1}$ is linearly independent, $\{\phi_j\}^{N}_{j=1}$ is a basis of the set $V_0$. We claim that 
\begin{equation}\label{projection p0}
\mcaP^2_0=\mcaP_0,\quad \mcaP^*_0=\mcaP_0,\quad \mcaP_0\ell^2(\Z^3)=V_0,
\end{equation}
which yields the desired result. 

To this end, let $A:=(a_{mj})_{N\times N}$ be the matrix of $S_2vA_2vS_2$ with respect to the orthonormal basis $\{f_j\}^{N}_{j=1}$ of $S_2\ell^2(\Z^3)$, where
\begin{equation*}
a_{mj}=\left<S_2vA_2vS_2f_m,f_j\right>=\left<A_2vf_m,vf_j\right>=\left<A_0vf_m,A_0vf_j\right>=\left<\phi_m,\phi_j\right>.
\end{equation*}
For any $g\in\ell^2(\Z^3)$, we  compute
\begin{equation*}
\mcaP_0g=A_0vS_2D_{2}S_2vA_0g=-\sum\limits_{j=1}^{N}\left<g,\phi_j\right>A_0vS_2D_{2}f_j=\sum\limits_{k,j=1}^{N}\left<g,\phi_j\right>b_{jk}\phi_k,
\end{equation*}
where $b_{jk}$ is the $(j,k)$ element of the inverse matrix $A^{-1}$ of $A$. Particularly, note that 
$$\mcaP_0\phi_m=\sum\limits_{k=1}^{N}\sum\limits_{j=1}^{N}a_{mj}b_{jk}\phi_k=\sum\limits_{k=1}^{N}(AA^{-1})_{mk}\phi_k=\phi_m,$$
which establishes \eqref{projection p0} and thus, we complete the proof. 
\end{proof}
Now we complement the proof of Lemma \ref{chara of solution and space}. 
\begin{proof}[Proof of Lemma {\rm\ref{chara of solution and space}}]

{\bf\underline{(i)}} ``$\bm{\Longrightarrow}$" Given $f\in S_1{\ell^2(\Z^3)}$, then
$(U+vA_0v)f=Sf=0.$
Multiplying both sides by $U$, it then yields that 
$f=-UvA_0vf$. Define $\phi:=-A_0vf$. Then $f=Uv\phi$ and hence $vf=V\phi$. Since $-\Delta A_0=\delta_0$, it follows that
$H\phi=(-\Delta+V)\phi=\Delta A_0vf+V\phi=-vf+vf=0.$
Let 
$$W:=\bigcap\limits_{s>1/2} \ell^{2, -s}(\Z^3).$$
To show that $\phi\in W$, it suffices to study the behavior of $\phi(n)$ for large $|n|$. Recall from Remark \ref{remark of free expan} (3) that  $A_0(n)=G(0,n)$. Assume that $|n|\geqslant5$. We decompose $(A_0vf)(n)$ as
\begin{align}
-\phi(n)&=(A_0vf)(n)=\sum\limits_{m}^{}G(0,n-m)v(m)f(m)\notag\\
=&\sum\limits_{m\neq n}^{}(G(0,n-m)-G_d(n-m))v(m)f(m)+\sum\limits_{m\neq n}^{}G_d(n-m)v(m)f(m)+G(0,0)v(n)f(n)\notag\\
=&:\Phi_1(n)+\Phi_2(n)+\Phi_3(n).\label{A0vf in S1}
\end{align}
Clearly, 
$|\Phi_3(n)|\lesssim\left<n\right>^{-\frac{\beta}{2}}\in W$ if $\beta>2$. Moreover, let $\mathbf{1}_{\{n\neq0\}}(n)$ be the characteristic function on the set $\{n\neq0\}$. By \eqref{relation of phi0 and g0}, we have
\begin{equation*}
    |\Phi_1(n)|\lesssim (\varphi*\left|vf\right|)(n),\quad \varphi(n)=\mathbf{1}_{\{n\neq0\}}(n)|n|^{-3}.
\end{equation*}
Since $vf\in\ell^1(\Z^3)$ for $\beta>3$ and $\varphi(n)\in\ell^2(\Z^3)$, Young's inequality gives
\begin{equation}\label{A}
\Phi_1(n)\in\ell^2(\Z^3)\subseteq W.
\end{equation}
For $\Phi_2(n)$, we divide it into the following three parts
\begin{align*}
\Phi_2(n)=\frac{1}{4\pi}\sum\limits_{m\neq n}^{}\frac{1}{|n-m|} v(m)f(m)=\frac{1}{4\pi}\sum\limits_{j=1}^{3}\sum\limits_{m\in R_j}\frac{1}{|n-m|} v(m)f(m)=:\frac{1}{4\pi}\sum\limits_{j=1}^{3}\Phi_{2j}(n),
\end{align*}
where $R_1=\{m\mid|m-n|=1\}$ and 
\begin{equation*}
   R_2=\{m\mid |m-n|>1,|m|\geqslant\frac{|n|}{2}\},\quad R_3=\{m\mid|m-n|>1,|m|<\frac{|n|}{2}\}.
\end{equation*}
For $\Phi_{21}(n)$, since $\left<m\right>\approx\left<n\right>$ for any $m\in R_1$, we obtain, for $\beta>5$,  
\begin{equation}\label{Phi21}
|\Phi_{21}(n)|\lesssim\left<n\right>^{-1}\sum\limits_{m\in R_1}\left<m\right>^{-\frac{\beta}{2}+1}|f(m)|\lesssim\left<n\right>^{-1}.
\end{equation}
For $\Phi_{22}(n)$, the Cauchy--Schwarz inequality and $\left<m\right>\gtrsim\left<n\right>$ for any $m\in R_2$ yield
 \begin{align}
    |\Phi_{22}(n)|^2&\lesssim \sum\limits_{m\in R_2}|n-m|^{-2}\left<m\right>^{-\beta}\lesssim\sum\limits_{1<|n-m|\leq |n|/4}\left<n\right>^{-\beta}+\sum\limits_{\substack{m\in R_2\\ |n-m|>|n|/4}}\left<n\right>^{-2}|m|^{-\beta}\notag\\
    &\lesssim\left<n\right>^{3-\beta}+\left<n\right>^{-2}\lesssim \left<n\right>^{-2}~({\rm since}\ \beta>5).\label{Phi22}
 \end{align} 
Finally, if $m \in R_3$, then
$|m-n| \approx |n|$. Hence, 
\begin{equation*}
|\Phi_{23}(n)|\lesssim \left<n\right>^{-1}\sum\limits_{m\in R_3}^{}|v(m)||f(m)|\lesssim \left<n\right>^{-1}.
\end{equation*}
Combining this with \eqref{Phi21} and \eqref{Phi22}, we conclude that $|\Phi_2(n)|\lesssim\left<n\right>^{-1}\in W $ and $\phi\in W$ follows.
\vskip0.2cm
``$\bm{\Longleftarrow}$" Conversely, suppose that there exists $\phi\in W$ satisfying the equation $H\phi=0$.
 Let $f=Uv\phi$. We will show that $\phi(n)=-(A_0vf)(n)$. Once this identity is established, it follows that
$$Sf=(U+vA_0v)f=v\phi-v\phi=0,$$
and hence $f\in S_1\ell^2(\Z^3)$. Indeed, since $H\phi=0$, we have $\Delta\phi=V\phi=vf=\Delta(-A_0vf).$
Moreover, by the argument in {\bf\underline{(i)~``$\bm{\Longrightarrow}$"}}, we have $A_0vf\in W$. In other words, the equation $\Delta(\phi+A_0vf)=0$ has solution in 
 $W$. By the uniqueness of solutions in $W$, it follows that $\phi=-A_0vf$.
\vskip0.2cm
{\bf{\underline{(ii)}}}~``$\bm{\Longrightarrow}$" Assume that $f\in S_{2}\ell^2(\Z^3)$, i.e., $f\in S_{1}\ell^2(\Z^3) $ and $\left<f,v\right>=0$. Then by part {\bf{\underline{(i)}}}, $\exists~\phi(n)=-(A_0vf)(n)\in W$ such that $H\phi=0$. We claim that $\phi \in\ell^2(\Z^3).$
In fact, for $|n|\geq5$, using $\left<f,v\right>=0$, we have 
\begin{equation}
\begin{split}
-\phi(n)&=(A_0vf)(n)=\sum\limits_{m}^{}\big(G(0,n-m)-(4\pi|n|)^{-1})v(m)f(m)\\
&=\Phi_1(n)+\sum\limits_{m\neq n}^{}\frac{|n|-|n-m|}{4\pi|n-m|\cdot |n|}v(m)f(m)+(G(0,0)-(4\pi|n|)^{-1})v(n)f(n)\\
&=:\Phi_1(n)+\widetilde{\Phi}_2(n)+\widetilde{\Phi}_3(n),
\end{split}
\end{equation}
where $\Phi_1(n)\in \ell^2(\Z^3)$ is defined in \eqref{A0vf in S1}. It is clear that $|\widetilde{\Phi}_3(n)|\lesssim \left<n\right>^{-\frac{\beta}{2}}\in \ell^2(\Z^3)$ if $\beta>3$.
As for $\widetilde{\Phi}_2(n)$, we note that
\begin{equation*}
    |\widetilde{\Phi}_2(n)|\lesssim \left<n\right>^{-1}\sum\limits_{m\neq n}\frac{1}{|n-m|}\left<m\right>^{-\frac{\beta-2}{2}}|f(m)|.
\end{equation*}
Since $\beta>7$, applying the same argument used for $\Phi_2(n)$ to the latter sum, we obtain $|\widetilde{\Phi}_2(n)|\lesssim\left<n\right>^{-2}\in \ell^2(\Z^3)$. Therefore, $\phi\in\ell^2(\Z^3)$.
\vskip0.3cm
``$\bm{\Longleftarrow}$" Suppose there exists $\phi\in \ell^2(\Z^3)$ such that $H\phi=0$. Let $f=Uv\phi$. Since $\ell^2(\Z^3)\subseteq W,$ it follows from part {\bf{\underline{(i)``$\bm{\Longleftarrow}$"}}} that $f\in S_1\ell^2(\Z^3)$ and
$\phi(n)=-(A_0vf)(n)$. It remains to show $\left<f,v\right>=0$, which will yield $f\in S_2\ell^2(\Z^3)$.
To this end, for any $\delta>0$, we introduce
$$F(\delta)=\sum_{n\in\Z^3}V(n)\phi(n)\eta(\delta n),$$
where $\eta(x)\in C^{\infty}_0(\R^3)$ satisfies $\eta(x)=1$ for $|x|\leqslant1$ and $\eta(x)=0$ for $|x|>2$. On the one hand, by the assumptions on $V$ and $\phi\in \ell^2(\Z^3)$, Lebesgue's dominated convergence theorem yields
$$\left<f,v\right>=\left<\phi,V\right>=\lim_{\delta\rightarrow0}F(\delta).$$
On the other hand, for any $\delta>0$, using $V\phi=\Delta \phi$ and the self-adjointness of $\Delta$, we get
\begin{equation*}
    F(\delta)=\sum_{n\in\Z^3}(\Delta \phi)(n)\eta(\delta n)=\sum_{n\in\Z^3}\phi(n)(\Delta G_{\delta})(n),
\end{equation*}
where $G_{\delta}(x)=\eta(\delta x)$. Applying the differential mean value theorem to $\Delta G_{\delta}$, we obtain
\begin{equation*}
 (\Delta G_{\delta})(n)=\delta^2\sum\limits_{j=1} ^{3}(\partial^2_j\eta)\big(\delta(n+(\Theta_j-1)e_j)\big), \quad \Theta_j\in[0,2].
 \end{equation*}
Restricting $0<\delta<1$, using  supp $\eta\subseteq\{|x|\leq2\}$ and uniform boundedness of $\partial^2_j\eta$, we further derive
 \begin{equation*}
 \|\Delta G_{\delta}\|^2_{\ell^2(\Z^3)}\lesssim \delta^4\sum\limits_{|n|\leqslant\frac{3}{\delta}}^{}\lesssim \delta.
 \end{equation*}
 Using $\phi\in\ell^2(\Z^3)$ and H\"{o}lder's inequality, we conclude $\lim\limits_{\delta\rightarrow0}F(\delta)=0$. This completes the proof.
\end{proof}

\subsection{Resolvent Expansions at Hyperbolic Thresholds \texorpdfstring{$4$}{4} and \texorpdfstring{$8$}{8}} \label{subsec of expan at 4}
\ \

\vskip0.1cm
Let $E^{\pm}_k,B^{\pm}_k,\mcaE^{\pm}_k,\mcaB^{\pm}_k$ be defined as in \eqref{expan at 4} and \eqref{expan of R08}, respectively. Define 
$$Q^{\pm}=U+vE^{\pm}_0v,\quad \mcaQ^{\pm}=U+v\mcaE^{\pm}_0v.$$
Denote $v_k(n)=(-1)^{n_k}v(n)$. Let $Q^{\pm}_1$ be the orthogonal projection onto $\operatorname{Ker}Q^{\pm}$, and let $Q^{\pm}_2$ be the orthogonal projection onto the kernel of $Q^{\pm}_1vB^{\pm}_0vQ^{\pm}_1$ restricted to $Q^{\pm}_1\ell^2(\Z^3)$. By \eqref{kernel B0pmn}, we have
\begin{equation*}
Q^{\pm}_1\ell^2(\Z^3)=\{f\in \ell^2(\Z^3)~|~ Q^{\pm}f=0\},\qquad Q^{\pm}_2\ell^2(\Z^3)=\{f\in Q^{\pm}_1\ell^2(\Z^3)~|\left<f,v_k\right>=0,\ k=1,2,3\}.
\end{equation*}
Similarly, let $\mcaQ^{\pm}_1$ be the orthogonal projection onto $\operatorname{Ker}\mcaQ^{\pm}$, and let $\mcaQ^{\pm}_2$ be the orthogonal projection onto the kernel of $\mcaQ^{\pm}_1v\mcaB^{\pm}_0v\mcaQ^{\pm}_1$ restricted to $\mcaQ^{\pm}_1\ell^2(\Z^3)$. Then
\begin{equation*}
\mcaQ^{\pm}_1\ell^2(\Z^3)=\{f\in \ell^2(\Z^3)~|~ \mcaQ^{\pm} f=0\},\qquad
\mcaQ^{\pm}_2\ell^2(\Z^3)=\{f\in \mcaQ^{\pm}_1\ell^2(\Z^3)~|\left<f,\mcaU v_k\right>=0,\ k=1,2,3\}.
\end{equation*}
\begin{proposition}\label{Asymptotic expansion theorem 4}
 Let $H=-\Delta+V$ with $|V(n)|\lesssim \left<n\right>^{-\beta}$ for some $\beta>0$. Then the following asymptotic expansions hold in $\mcaB(\sigma,-\sigma)$ as $|\omega|\rightarrow0$ with $\pm\operatorname{Im}\omega>0$.  
\begin{itemize}
\item[{\rm(1)}] If $4$ {\rm(}resp. $8${\rm)} is a regular point of $H$ and $\beta>7$, then for any $\sigma>\frac{7}{2}$, 
\begin{align*}
\begin{split}
R^{\pm}_V(4+\omega)&=E^{\pm}_0-E^{\pm}_0vF^{\pm}_0vE^{\pm}_0+(B^{\pm}_0+E^{\pm}_0vF^{\pm}_0vB^{\pm}_0vF^{\pm}_0vE^{\pm}_0-B^{\pm}_0vF^{\pm}_0vE^{\pm}_0-E^{\pm}_0vF^{\pm}_0vB^{\pm}_0)\omega^{\frac{1}{2}}+O(\omega),\\
    R^{\pm}_V(8-\omega)&=\mcaE^{\pm}_0-\mcaE^{\pm}_0v\mcaF^{\pm}_0v\mcaE^{\pm}_0+(\mcaB^{\pm}_0+\mcaE^{\pm}_0v\mcaF^{\pm}_0v\mcaB^{\pm}_0v\mcaF^{\pm}_0v\mcaE^{\pm}_0-\mcaB^{\pm}_0v\mcaF^{\pm}_0v\mcaE^{\pm}_0-\mcaE^{\pm}_0v\mcaF^{\pm}_0v\mcaB^{\pm}_0)\omega^{\frac{1}{2}}+O(\omega).
    \end{split}
\end{align*}
\item [{\rm(2)}] If $4$ {\rm(}resp. $8${\rm)} is an eigenvalue of $H$ and $\beta>11$, then for any $\sigma>\frac{11}{2}$,
\begin{align*}
R^{\pm}_V(4+\omega)=-E^{\pm}_0vQ^{\pm}_2F^{\pm}_2Q^{\pm}_2vE^{\pm}_0\omega^{-1}+O(\omega^{-\frac{1}{2}}),\quad  R^{\pm}_V(8-\omega)=-\mcaE^{\pm}_0v\mcaQ^{\pm}_2\mcaF^{\pm}_2\mcaQ^{\pm}_2v\mcaE^{\pm}_0\omega^{-1}+O(\omega^{-\frac{1}{2}}).
\end{align*}
\end{itemize}
Here  
\begin{align*}
F^{\pm}_{0}&=(Q^{\pm}+Q^{\pm}_1)^{-1},\quad \mcaF^{\pm}_{0}=(\mcaQ^{\pm}+\mcaQ^{\pm}_1)^{-1},\\
F^{\pm}_2&=(Q^{\pm}_2vE^{\pm}_1vQ^{\pm}_2)^{-1}\ {\rm on}\ Q^{\pm}_2\ell^2(\Z^3),\quad \mcaF^{\pm}_2=(\mcaQ^{\pm}_2v\mcaE^{\pm}_1v\mcaQ^{\pm}_2)^{-1}\ {\rm on}\ \mcaQ^{\pm}_2\ell^2(\Z^3),
\end{align*}
and  $E^{\pm}_0vQ^{\pm}_2F^{\pm}_2Q^{\pm}_2vE^{\pm}_0$ and $-\mcaE^{\pm}_0v\mcaQ^{\pm}_2\mcaF^{\pm}_2\mcaQ^{\pm}_2v\mcaE^{\pm}_0$ are the orthogonal projections onto the eigenspaces of $H$ at $4$ and $8$, respectively. Moreover, these expansions are twice differentiable in $\mcaB(\sigma,-\sigma)$ with $\sigma>\frac{11}{2}$.
\end{proposition}
 Unlike the elliptic threshold 0, at the hyperbolic thresholds $4$ and $8$, no comparably explicit kernel asymptotics are available for $E^\pm_0$ and $\mcaE^\pm_0$. Thus, the Green-function argument used above to characterize threshold states cannot be applied directly. Instead, we make use of the following results of Nomura and Taira \cite{NT20}, whose analysis of threshold states and the geometry of the Fermi surface provides the substitute for the explicit kernel information available at the elliptic threshold. 
\begin{lemma}\label{Taira lemma}
{\rm (\cite[Theorem 1.7, Lemmas 3.6, 3.7 and Remark 3.8]{NT20})} Let $H=-\Delta+V$ on $\Z^3$ with a real-valued potential $V$ satisfying $|V(n)|\lesssim\langle n\rangle^{-\beta}$ for some $\beta>0$. Let $\lambda\in\{4,8\}$. Then the following statements hold.

{\rm (1)} Let $\beta>7/2$. If $\phi\in\ell^{2,-1}(\Z^3)$ satisfies $\phi+(-\Delta-\lambda\mp i0)^{-1}V\phi=0$, then 
$\phi\in\ell^2(\Z^3)$ and
$\mcaF(V\phi)|_{\mcaM_{\lambda}}=0$, where
\begin{equation*}
    \mcaM_{\lambda}=\{x\in \T^3\mid \mcaM(x)=\lambda\}.
\end{equation*}

{\rm(2)} Let $f\in C^1(\T^3)$ such that $f|_{\mcaM_{\lambda}}=0$. Then, we have $(\mcaM(x)-\lambda\mp i0)^{-1}f=(\mcaM(x)-\lambda)^{-1}f\in L^2(\T^3)$.
\end{lemma}
Using this lemma, we now characterize the above subspaces in terms of threshold states. Since the arguments at the thresholds $4$
and $8$ are analogous, we only treat the threshold $4$. More precisely,
we have the following characterization.
\begin{lemma}\label{chara of solution and space 4}
{\rm Let $H=-\Delta+V$ with $|V(n)|\lesssim \left<n\right>^{-\beta}$ for $\beta>7$. Then
\begin{itemize}
\item [{\rm(i)}]
$f^{\pm}\in Q^{\pm}_1\ell^2(\Z^3)\Longleftrightarrow \exists\ \phi^{\pm}\in \ell^{2,-1}(\Z^3)$ such that $\phi^{\pm}+(-\Delta-4\mp i0)^{-1}V\phi^{\pm}=0$. Moreover, $f^{\pm}=Uv\phi^{\pm}$ and $\phi^{\pm}=-E^{\pm}_0vf^{\pm}$.
\item [{\rm(ii)}] $Q^+_1\ell^2(\Z^3)=Q^-_1\ell^2(\Z^3)$.
\end{itemize}
}
\end{lemma}

\begin{remark}\label{remark of chara 4}
\ \

{\rm

(1)  Under the decay assumption of Lemma \ref{chara of solution and space 4}, there is no genuine resonance at the threshold $4$. Indeed, if
$f^\pm\in Q^\pm_1\ell^2(\Z^3)$, then by Lemma
\ref{chara of solution and space 4} there exists
$\phi^\pm\in\ell^{2,-1}(\Z^3)$ satisfying
$$\phi^{\pm}+(-\Delta-4\mp i0)^{-1}V\phi^{\pm}=0\quad {\rm and}\quad f^\pm=Uv\phi^\pm.$$
Lemma \ref{Taira lemma} (1) then implies that
$\phi^\pm\in\ell^2(\Z^3)$ and 
$\mcaF(V\phi^\pm)|_{\mcaM_4}=0.$
In particular, for $\xi_k:=\pi e_k\in\mcaM_4$, $k=1,2,3$, we have
$$
\left<f^\pm,v_k\right>=(2\pi)^{\frac{3}{2}}\mcaF\big(V\phi^\pm\big)\big(\xi_k\big)=0. 
$$
Hence, by the definition of $Q^\pm_2$, $Q^\pm_1\ell^2(\Z^3) = Q^\pm_2\ell^2(\Z^3).$ In other words, $Q_1^\pm vB_0^\pm vQ_1^\pm=0 $ on $Q_1^\pm\ell^2(\Z^3).$
\vskip0.1cm
(2) Consequently, combining Lemma
\ref{chara of solution and space 4}~(i) with Lemma
\ref{Taira lemma}~(1), the threshold $4$ admits only the regular and
eigenvalue cases. More precisely,
\begin{align}\label{char of Q1pm}
\begin{split}
Q^{\pm}_1\ell^2(\Z^3)&=\{0\} \Longleftrightarrow\phi+(-\Delta-4\mp i0)^{-1}V\phi=0 \ \text{has only the zero solution in }\ell^{2,-1}(\Z^3),\\
Q^{\pm}_1\ell^2(\Z^3)&\neq\{0\} \Longleftrightarrow\phi+(-\Delta-4\mp i0)^{-1}V\phi=0 \ \text{has non-zero solution in }\ell^{2}(\Z^3).
\end{split}
\end{align}
Moreover, Lemma \ref{chara of solution and space 4}~(ii) shows that the
two a priori different classifications arising from the upper and lower
half-plane boundary values coincide.

Arguing as at the threshold $0$, $Q^\pm$ is a Fredholm operator of index zero, and hence
$Q^\pm$ is invertible if and only if $\ker Q^\pm=\{0\}$. Combining this with \eqref{char of Q1pm} and Definition
\ref{defin of resonance types 4,8}, we obtain
\begin{itemize}
    \item 
    $4$ is a regular point of $H$ $\Longleftrightarrow$ $Q^+$ is invertible on $\ell^2(\Z^3)$ $\Longleftrightarrow$ $Q^-$ is invertible on $\ell^2(\Z^3)$.
\vskip0.1cm
\item $4$ is an eigenvalue of $H$ $\Longleftrightarrow$ $Q^+$ is not invertible on $\ell^2(\Z^3)$$\Longleftrightarrow$ $Q^-$ is not invertible on $\ell^2(\Z^3)$. 
\end{itemize}  
}
\end{remark}

\begin{proof} [Proof of Lemma {\rm \ref{chara of solution and space 4}}]
First, taking $N=0$ and letting $\omega=\pm i\varepsilon$ with $\varepsilon>0$ in \eqref{expan at 4}, we have 
\begin{equation*}
E^{\pm}_0=\lim\limits_{\varepsilon\rightarrow0,\ \varepsilon>0}R_0(4\pm i\varepsilon)\quad in\ \mcaB(s,-s),\quad s>7/2.
\end{equation*}
Then it follows that ${E}^{\pm}_0=(-\Delta-4\mp i0)^{-1}$ on $\ell^{2,s}(\Z^3)$ for $s>7/2.$
\vskip0.1cm
{\bf\underline{(i)}} ``$\bm{\Longrightarrow}$" Let $f^{\pm}\in Q^{\pm}_1{\ell^2(\Z^3)}$. Then $(U+vE^{\pm}_0v)f^{\pm}=Q^{\pm}f^{\pm}=0.$ Multiplying both sides by $U$, we obtain
$f^{\pm}=-UvE^{\pm}_0vf^{\pm}$. Since $\beta>7$, we have $vf^{\pm}\in \ell^{2,s}(\Z^3)\subseteq\ell^{2,1}(\Z^3)$ for some $7/2<s\leq\beta/{2}$. Using
$(-\Delta-4\mp i0)^{-1}\in\mcaB(1,-1)$, we may therefore define $$\phi^{\pm}:=-E^{\pm}_0vf^{\pm}=-(-\Delta-4\mp i0)^{-1}vf^{\pm}\in\ell^{2,-1}(\Z^3).$$ Then,
$f^{\pm}=Uv\phi^{\pm}$ and $vf^{\pm}=V\phi^{\pm}$, so that
$\phi^{\pm}+(-\Delta-4\mp i0)^{-1}V\phi^{\pm}=0.$ 
\vskip0.1cm
``$\bm{\Longleftarrow}$" Suppose that there exists $\phi^{\pm}\in \ell^{2,-1}(\Z^3)$ satisfying the equation $$\phi^{\pm}+(-\Delta-4\mp i0)^{-1}V\phi^{\pm}=0.$$
Let $f^{\pm}:=Uv\phi^{\pm}\in\ell^2(\Z^3)$. Then $vf^{\pm}=V\phi^{\pm}\in\ell^{2,s}(\Z^3)$ for some $7/2<s\leq \beta-1$. Therefore, $$\phi^{\pm}=-(-\Delta-4\mp i0)^{-1}V\phi^{\pm}=-E^{\pm}_0vf^{\pm}.$$
Hence 
 $Q^{\pm}f^{\pm}=(U+vE^{\pm}_0v)f^{\pm}=v\phi^{\pm}-v\phi^{\pm}=0,$
which yields $f^{\pm}\in Q^{\pm}_1\ell^2(\Z^3)$.
\vskip0.1cm
{\bf\underline{(ii)}} For any $f\in Q^+_1\ell^2(\Z^3)$, it follows from (i) that there exists $\phi\in\ell^{2,-1}(\Z^3)$ such that $$\phi+(-\Delta-4- i0)^{-1}V\phi=0\ \ \text{and}\ \ f=Uv\phi.$$ Applying Lemma \ref{Taira lemma} (1) to $\phi$, we obtain $\mcaF(V\phi)\big|_{\mcaM_4}=0$. Since $V\phi\in\ell^{2,s}(\Z^3)$ for some $s>\frac{7}{2}$, $\mcaF(V\phi)\in C^1(\T^3)$ by the Sobolev embedding theorem. By Lemma \ref{Taira lemma} (2), we have
\begin{equation*}
    (\mcaM(x)-4\mp i0)^{-1}\mcaF(V\phi)(x)=(\mcaM(x)-4)^{-1}\mcaF(V\phi)(x).
\end{equation*}
This implies that $E^+_0vf=E^+_0V\phi=E^-_0V\phi=E^-_0vf$. Thus, $f\in Q^-_1\ell^2(\Z^3)$ follows from 
$$Q^-f=(U+vE^-_0v)f=Uf+vE^+_0vf=Q^+f=0.$$
Similarly, we can prove $Q^-_1\ell^2(\Z^3)\subseteq Q^+_1\ell^2(\Z^3)$. Hence, $Q^+_1\ell^2(\Z^3)=Q^-_1\ell^2(\Z^3)$ is established. This completes the proof.
\end{proof}

\vskip0.4cm
\begin{proof}[Proof of Proposition {\rm\ref{Asymptotic expansion theorem 4}}]
We first verify that $F_0^\pm$ and $F_2^\pm$ are well-defined. Since $(E^+_0)^*=E^-_0$, we have $(Q^\pm)^*=Q^\mp$. By Lemma
\ref{chara of solution and space 4} (ii), $$\operatorname{Ker}(Q^{\pm})^{*}=\operatorname{Ker}Q^{\mp}=\operatorname{Ran}Q^{\mp}_1=\operatorname{Ran}Q^{\pm}_1=\operatorname{Ker}Q^{\pm}.$$ Moreover, $Q^\pm$ is a Fredholm operator
of index zero and hence has closed range. Therefore, Remark
\ref{remark of invertibility} implies that   $F^{\pm}_0=(Q^{\pm}+Q^\pm_1)^{-1}$ is well-defined.

It remains to consider $F^{\pm}_2= \left(Q_2^\pm vE_1^\pm vQ_2^\pm\right)^{-1}$.  Since $Q^\pm$ is Fredholm, $Q_1^\pm$, and hence $Q_2^\pm$, is of finite
rank. Thus it suffices to prove injectivity. Let $f^\pm\in Q_2^\pm\ell^2(\Z^3)$. By Remark
\ref{remark of chara 4}, $\mcaF(vf^\pm)\big|_{\mcaM_4}=0.$ 
Moreover, the decay of $v$ implies
$vf^\pm\in\ell^{2,s}(\Z^3)$ for some $s>\frac52$, and hence
$\mcaF(vf^\pm)\in C^1(\T^3)$. Lemma \ref{Taira lemma}~(2) therefore gives
$$ (\mcaM(x)-4)^{-1}\mcaF(vf^\pm)\in L^2(\T^3).$$
Proceeding as in the proof of the injectivity of $S_2vA_2vS_2$ at the
threshold $0$, we conclude that
$Q_2^\pm vE_1^\pm vQ_2^\pm$ is injective and hence invertible on
$Q_2^\pm\ell^2(\Z^3)$. Therefore $F_2^\pm$ is well-defined.

\vskip0.1cm
We next turn to the inversion of
$$M^{\pm}(4+\omega):=U+vR^{\pm}_0(4+\omega)v, \qquad \pm\operatorname{Im}\omega>0.$$ 
Using $Q_1^\pm vB_0^\pm vQ_1^\pm=0$, and applying the same iterative inversion scheme as at the threshold $0$, we obtain
\begin{itemize}
    \item If 4 is a regular point of $H$ and $\beta>7$, then 
\begin{equation}\label{asy expan regu 4}
(M^{\pm}(4+\omega))^{-1}=F^\pm_{0}-F^\pm_{0}vB^\pm_0vF^\pm_{0}\omega^{\frac{1}{2}}+F^\pm_0O(\omega)F^\pm_0, \quad|\omega|\rightarrow0,\quad\pm\operatorname{Im} \omega>0.
\end{equation}
\item If 4 is an eigenvalue of $H$ and $\beta>11$, then 
\begin{equation}\label{asy expan 2nd 4}
(M^{\pm}(4+\omega))^{-1}=Q^\pm_2F^\pm_{2}Q^\pm_2\omega^{-1}+O(\omega^{-\frac{1}{2}}),\quad|\omega|\rightarrow0,\quad\pm\operatorname{Im} \omega>0.
\end{equation}
\end{itemize} 
Hence, combining these expansions with \eqref{rela betw RV and inver of Mw}, we obtain the corresponding expansions for $R^\pm_V(4+\omega)$. 
By combining the derivative estimates for the free resolvent in \cite[Corollary 2.10]{Kop09} with the above inversion argument, the expansions in Proposition \ref{Asymptotic expansion theorem 4} can be differentiated twice in $\mcaB(\sigma,-\sigma)$ for $\sigma>\frac{11}{2}$.

Moreover, by Remark \ref{remark on defini of resonances} and the argument used in the proof of Proposition \ref{Asymptotic expansion theorem 0}, we can also show that $E^{\pm}_0vQ^{\pm}_2F^{\pm}_2Q^{\pm}_2vE^{\pm}_0$ is the eigenprojection of $H$ associated with the eigenvalue $4$. This completes the proof.
\end{proof}

\section{Application: Jensen-Kato Type Time-decay Estimates for Schr\"odinger Equations.}\label{sec of application}
In this section, we prove the Jensen--Kato type time-decay estimates for the discrete Schr\"odinger equation on $\Z^3$ stated in Theorem \ref{decay theorem}. We begin with the following limiting absorption principle.
 \begin{lemma}\label{LAP-theorem}
    { Let $H=-\Delta+V$ with $|V(n)|\lesssim \left<n\right>^{-\beta}$ for some $\beta>2$ and $\mcaI=(0,4)\cup(4,8)\cup(8,12)$. Let $[\beta]$ denote the greatest integer not exceeding $\beta$. Then the following statements hold.}
\begin{itemize}
{
\item [{\rm(i)}] The point spectrum $\sigma_p(H)\cap\mcaI$ is discrete, each eigenvalue has finite multiplicity, and the singular continuous spectrum $\sigma_{sc}(H)=\varnothing$.
\vskip0.2cm
\item [{\rm(ii)}] Let $j\in\left\{0,\cdots,[\beta]-2\right\}$ and $\sigma>j+\frac{1}{2}$. Then the norm limits
\begin{equation}\label{boundary value}
   \frac{d^{j}}{d\lambda^j}(\mcaR^{\pm}_V(\lambda))=\lim\limits_{\varepsilon\downarrow0}R^{(j)}_V(\lambda\pm i\varepsilon) \quad {\rm in} \quad \mathcal{B}(\sigma,-\sigma)
    \end{equation}
exist and are norm-continuous on $\mcaI\setminus\sigma_p(H)$, where $R^{(j)}_{V}(\omega)$ denotes the $j$-th derivative of $R_{V}(\omega)$.
    }
\end{itemize}
\end{lemma}
This lemma is a special case of Theorem \ref{general LAP} with $d=3$. It extends the results of \cite{IK12,Kop09}, where the LAP in $\mcaB(s,-s)$ was established for finitely supported potentials. For an overview of the LAP and a proof of this lemma on $\Z^d$, we refer to Appendix \ref{LAP general d}. 

Using the above lemma together with Theorems \ref{Asymptotic expansion theorem_0}$\sim$\ref{Asymptotic expansion theorem_4} and the stationary representation derived in Appendix \ref{Stationary repres of eitHpacH}, we have
\begin{equation}\label{expre of Dt}
   e^{-itH}P_{ac}(H)=\sum\limits_{k=1}^{3}\int_{4(k-1)}^{4k}e^{-it\lambda}\mcaR(\lambda)d\lambda,\quad \mcaR(\lambda):=\frac{1}{2\pi i}(\mcaR_V^{+}-\mcaR_V^{-})(\lambda).
\end{equation}
Based on this representation, we proceed as follows. Fix
$0<\delta\ll1$ and choose a partition of unity on $(0,4)$,
\begin{equation*}
    \xi_1(\lambda)
+\xi_2(\lambda)
+\xi_3(\lambda)=1,\quad \xi_j\in C^{\infty}((0,4)),\quad j=1,2,3,
\end{equation*}
such that supp~$\xi_1\subseteq(0,2\delta\,]$, supp~$\xi_2\subseteq[\delta,4-\delta]$ and supp~$\xi_3\subseteq[4-2\delta,4)$.
We first apply the change of variables $\lambda=4k-u$ to \eqref{expre of Dt} and then use the above partition of unity to obtain
    \begin{equation*}
   e^{-itH}P_{ac}(H)=\sum\limits_{k=1}^{3} \sum\limits_{j=1}^{3} \int_{0}^{4}e^{-it(4k-u)}\xi_j(u)\mcaR(4k-u)du.
\end{equation*}
For the terms involving $\xi_3$, the change of variables $u=4-\lambda$ gives 
\begin{equation*}
    \int_{0}^{4}e^{-it(4k-u)}\xi_3(u)\mcaR(4k-u)du=\int_{0}^{4}e^{-it(4(k-1)+\lambda)}\xi^+_1(\lambda)\mcaR(4(k-1)+\lambda)d\lambda,\quad \xi^+_1(\lambda):=\xi_3(4-\lambda).
\end{equation*}
Combining the above identities, we decompose $e^{-itH}P_{ac}(H)$ into the following five terms:
\begin{equation*}
e^{-itH}P_{ac}(H)=\mcaG(t)+ \sum\limits_{j=0}^{3}\mcaD_j(t),
\end{equation*}
where $\xi^-_1(\lambda):=\xi_1(\lambda)$ and 
\begin{align}\label{gt and Dkt}
\begin{split}
\mcaG(t)&=\sum\limits_{k=1}^3\int_{0}^{4}e^{-it(4k-\lambda)}\xi_2(\lambda)\mcaR(4k-\lambda)d\lambda,\\
\mcaD_0(t)&=\int_{0}^4e^{-it\lambda}\xi^+_1(\lambda)\mcaR(\lambda)d\lambda,\quad \mcaD_3(t)=\int_0^4 e^{-it(12-\lambda)}\xi^-_1(\lambda)\mcaR(12-\lambda)d\lambda, \\
\mcaD_k(t)&=\sum\limits_{\pm}\int_{0}^4e^{-it(4k\pm\lambda)}\xi^{\pm}_1(\lambda)\mcaR(4k\pm\lambda)d\lambda,\quad k=1,2.\\
\end{split}
\end{align}
The proof of Theorem \ref{decay theorem} therefore reduces to establishing the decay estimates for $\mcaG(t)$ and $\{\mcaD_k(t)\}_{k=0}^3$.
\vskip0.3cm
\begin{proposition}\label{theorem of estimate of Ij}
{\rm Under the assumptions on $H$ and $V$ in Theorem \ref{decay theorem}, the following decay estimates hold.
\vskip0.1cm
{\bf(1) Estimates for $\mcaG(t)$.}
\vskip0.1cm
For any $\sigma>\frac{5}{2}$, we have
$\|\mcaG(t)\|_{\mathcal{B}(\sigma, -\sigma)} = O(t^{-2})$ as $t\rightarrow+\infty$.
\vskip0.2cm
{\bf(2) Estimates for $\mcaD_0(t)$ and $\mcaD_3(t)$.} For $k=0,3$, 
\begin{equation*}
    \|\mcaD_{k}(t)\|_{\mcaB(\sigma,-\sigma)}=\begin{cases}
        O(t^{-\frac{3}{2}}),&4k~\text{is a regular point of}~ H\\
        O(t^{-\frac{1}{2}}),&4k~\text{is a resonance or an eigenvalue of}~ H
    \end{cases},\quad t\rightarrow+
\infty,
\end{equation*}
where $\sigma>5/2$ in the regular and resonance cases, and $\sigma>7/2$ in the eigenvalue case. 
\vskip0.2cm
{\bf(3) Estimates for $\mcaD_1(t)$ and $\mcaD_2(t)$.} For $k=1,2$ and $\sigma>11/2$, 
\begin{equation*}
\|\mcaD_k(t)\|_{\mcaB(\sigma,-\sigma)}=\begin{cases}
        O(t^{-\frac{3}{2}}),&4k~\text{is a regular point of}~ H\\
        O(t^{-\frac{1}{2}}),&4k~~\text{is an eigenvalue of}~ H
    \end{cases},\quad t\rightarrow+\infty.
\end{equation*}
}
\end{proposition}
\vskip0.2cm
\begin{lemma}\label{decay rate lemma}
Let $0<\delta\ll1$ and $\alpha\in(-1,1)$. Suppose that $F\in C^{2}((0,\delta],\C)$ satisfies \begin{equation}\label{condi for F}
    F(\delta)=0=F'(\delta),\quad F^{(k)}(\lambda)= O(\lambda^{\alpha-k}),\quad\lambda\in(0,\delta],\quad k=0,1,2.
\end{equation} 
Then
\begin{equation*}
\int_{0}^{\delta}e^{-it\lambda}F(\lambda)d\lambda
=O(t^{-(1+\alpha)}),\quad t\rightarrow+\infty.
\end{equation*}
\end{lemma}
As a consequence, the same conclusion holds for operator-valued functions $F\in C^{2}((0,\delta],\mcaB(\sigma,-\sigma))$ provided that \eqref{condi for F} holds in the operator norm of
$\mcaB(\sigma,-\sigma)$. Indeed, for any $f,g\in\ell^{2,\sigma}(\Z^3)$ with $\|f\|_{\ell^{2,\sigma}(\Z^3)}=1=\|g\|_{\ell^{2,\sigma}(\Z^3)}$, define $F_{f,g}(\lambda)=\left<F(\lambda)f,g\right>$. Then
$F_{f,g}$ satisfies the assumptions of Lemma
\ref{decay rate lemma} uniformly in $f$ and $g$, and therefore
  \begin{equation*}
\Big\|\int_{0}^{\delta}e^{-it\lambda}F(\lambda)d\lambda\Big\|_{\mcaB(\sigma,-\sigma)}=\sup\limits_{f,g}\Big|\Big<\int_{0}^{\delta}e^{-it\lambda}F(\lambda)d\lambda f,g\Big>\Big|=\sup\limits_{f,g}\Big|\int_{0}^{\delta}e^{-it\lambda}F_{f,g}(\lambda)d\lambda\Big|=O(t^{-(1+\alpha)}).
  \end{equation*}
We defer the proof of Lemma \ref{decay rate lemma} to the end of this section.

\begin{proof}[Proof of Proposition {\rm\ref{theorem of estimate of Ij}}]
{\bf (1)} Recall from \eqref{gt and Dkt} and supp $ \xi_2\subseteq[\delta,4-\delta]$ that
\begin{equation*}
\mcaG(t)=\sum\limits_{k=1}^3e^{-4kit}\int_{\delta}^{4-\delta}e^{it\lambda}\xi_2(\lambda)\mcaR(4k-\lambda)d\lambda,\quad \xi_2\in C^{\infty}((0,4)).
\end{equation*}
By Lemma \ref{LAP-theorem}, $\mcaR(4k-\lambda)\in C^{2}([\delta,4-\delta],\mcaB(\sigma,-\sigma))$ for $\sigma>\frac{5}{2}$. Moreover, the support condition on $\xi_2$ implies $\xi_2(\delta)=\xi_2(4-\delta)
=
\xi_2'(\delta)=\xi_2'(4-\delta)=0$. Integrating by parts twice in the inner integral of $\mcaG(t)$ then gives $$\|\mcaG(t)\|_{\mathcal{B}(\sigma, -\sigma)} = O(t^{-2}),\quad t\rightarrow+\infty.$$
\vskip0.2cm
{\bf (2)} From \eqref{gt and Dkt} and supp $\xi^{\pm}_1(\lambda)\subseteq(0,2\delta]$, we have
\begin{equation*}
\mcaD_0(t)=\int_{0}^{2\delta}e^{-it\lambda}\xi^+_1(\lambda)\mcaR(\lambda)d\lambda,\quad \mcaD_3(t)=e^{-12it}\int_0^{2\delta} e^{it\lambda}\xi^-_1(\lambda)\mcaR(12-\lambda)d\lambda.
\end{equation*}
It suffices to consider $\mcaD_0(t)$, since the argument for
$\mcaD_3(t)$ is analogous. Combining Lemma \ref{LAP-theorem} with
Theorem \ref{Asymptotic expansion theorem_0}, we obtain, for
$j=0,1,2$ and $0<\lambda\leq2\delta$,
\begin{equation*}
\|\mathcal{R}^{(j)}(\lambda)\|_{\mathcal{B}(\sigma,-\sigma)}=\begin{cases}
    O(\lambda^{\frac{1}{2}-j}),&0~\  is\ a\ regular\ point\ of\ H\\
    O(\lambda^{-\frac{1}{2}-j}),&0~\ is\ a\ resonance\ or\ an\ eigenvalue\ of\ H
\end{cases}
\end{equation*}
where $\sigma>\frac52$ in the regular and resonance cases and
$\sigma>\frac72$ in the eigenvalue case.
Applying Lemma \ref{decay rate lemma} to $\xi^+_1(\lambda)\mcaR(\lambda)$, we derive the desired estimates.
\vskip0.2cm
{\bf(3)} We only deal with $\mcaD_1(t)$, since the case of $\mcaD_2(t)$ is similar. Suppose first that  4 is a regular point of $H$. Combining Lemma \ref{LAP-theorem} with Theorem \ref{Asymptotic expansion theorem_4}, we obtain, for $0<\lambda\leq 2\delta$, 
\begin{equation*}
    \mathcal{R}(4\pm\lambda)=
        \frac{1}{2\pi i}(\B^+_0-\B^-_0)+\widetilde{\mathcal{R}}_{\pm}(\lambda),
\end{equation*}
where $\widetilde{\mathcal{R}}_{\pm}(\lambda)$ satisfies
\begin{equation*}
\|\widetilde{\mathcal{R}}^{(j)}_{\pm}(\lambda)\|_{\mathcal{B}(\sigma,-\sigma)}=
    O(\lambda^{\frac{1}{2}-j}),\quad j=0,1,2,\quad \sigma>7/2.
\end{equation*}
Then 
\begin{align*}
\mcaD_1(t)&=\sum\limits_{\pm}\int_{0}^4e^{-it(4\pm\lambda)}\xi^{\pm}_1(\lambda)\mcaR(4\pm\lambda)d\lambda\\
&=\frac{e^{-4it}}{2\pi i}(\B^+_0-\B^-_0)\sum\limits_{\pm}\int_{0}^4e^{\mp it\lambda}\xi^{\pm}_1(\lambda)d\lambda+\sum\limits_{\pm}\int_{0}^4e^{-it(4\pm\lambda)}\xi^{\pm}_1(\lambda)\widetilde{\mcaR}_{\pm}(\lambda)d\lambda,\\
&=:\mcaD_{11}(t)+\mcaD_{12}(t).
\end{align*}
Since $\xi^{\pm}_1(4)=0$ and $\xi^{\pm}_1(0)=1$, integration by parts gives 
\begin{equation*}
    \int_{0}^4e^{\mp it\lambda}\xi^{\pm}_1(\lambda)d\lambda=\frac{\mp i}{t}+\frac{\mp i}{t}\int_{0}^{4}e^{\mp it\lambda}\big(\xi^{\pm}_1\big)^{'}(\lambda)d\lambda=\frac{\mp i}{t}+O(t^{-2}),\quad t\rightarrow+\infty.
\end{equation*}
 The leading $t^{-1}$ terms cancel after summing over ``$\pm$", and hence
$$\|\mcaD_{11}(t)\|_{\mcaB(\sigma,-\sigma)}=O(t^{-2}),\quad t\rightarrow +\infty,\quad \sigma>7/2.$$
On the other hand, applying Lemma \ref{decay rate lemma} to $\xi^{\pm}_1(\lambda)\widetilde{\mcaR}_{\pm}(\lambda)$ yields  
$$\|\mcaD_{12}(t)\|_{\mcaB(\sigma,-\sigma)}=O(t^{-3/2}),\quad t\rightarrow +\infty,\quad \sigma >7/2.$$
Therefore, when $4$ is a regular point of $H$ and 
$\sigma>7/2$,
\begin{equation*} \|\mcaD_1(t)\|_{\mcaB(\sigma,-\sigma)}=
O(t^{-3/2}),\quad t\rightarrow+\infty.
\end{equation*}

Suppose now that $4$ is an eigenvalue of $H$. Combining Lemma \ref{LAP-theorem} with Theorem \ref{Asymptotic expansion theorem_4}, we obtain 
\begin{equation*}
\|{\mathcal{R}}^{(j)}(4\pm\lambda)\|_{\mathcal{B}(\sigma,-\sigma)}=
    O(\lambda^{-\frac{1}{2}-j}),\quad j=0,1,2,\quad \sigma>11/2.
\end{equation*}
Using Lemma \ref{decay rate lemma}, we have \begin{equation*} \|\mcaD_1(t)\|_{\mcaB(\sigma,-\sigma)}=
O(t^{-1/2}),\quad t\rightarrow+\infty,\quad \sigma>11/2.
\end{equation*} 
Hence, we complete the proof.
\end{proof}

\begin{proof}[Proof of {\rm Lemma \ref{decay rate lemma}}]
For sufficiently large $t$, we have $t^{-1}<\delta$. We split the integral as
\begin{equation*}
\int_{0}^{\delta} e^{-it\lambda} F(\lambda) d\lambda=\Big(\int_{0}^{\frac{1}{t}} +\int_{\frac{1}{t}}^{\delta} \Big)e^{-it\lambda} F(\lambda) d\lambda =:I_1(t)+I_2(t).
\end{equation*}
 Since $F(\lambda)=O(\lambda^{\alpha})$ on $(0,\delta\,]$ and $\alpha>-1$, we have $|I_1(t)|\lesssim t^{-(\alpha+1)}$. We next estimate $I_2(t)$. Integrating by parts and using $F(\delta)=0$ together with $F(t^{-1})=O(t^{-\alpha})$, we obtain
\begin{equation*}
\begin{split}
    I_2(t)=\frac{e^{-i}}{it}F(t^{-1}) +\frac{1}{it}\int_{\frac{1}{t}}^{\delta}e^{-it\lambda}F^{'}(\lambda)d\lambda
    &=O(t^{-(1+\alpha)})+\frac{1}{it}\int_{\frac{1}{t}}^{\delta}e^{-it\lambda}F^{'}(\lambda)d\lambda.
\end{split}
\end{equation*}
For the second integral,  suppose first that $-1<\alpha<0$. Using $F'(\lambda)=O(\lambda^{\alpha-1})$, we obtain
\begin{equation*}
    \Big|\frac{1}{it}\int_{\frac{1}{t}}^{\delta}e^{-it\lambda}F^{'}(\lambda)d\lambda\Big|\lesssim \frac{1}{t}\int_{\frac{1}{t}}^{\delta}\lambda^{\alpha-1} d\lambda=\frac{1}{\alpha t}(\delta^{\alpha}-t^{-\alpha})\leq \frac{-1}{\alpha}t^{-(1+\alpha)}.
\end{equation*}
When $0\leq\alpha<1$, integrating by parts once more and using $F'(\delta)=0$, we derive 
\begin{equation*}
    \Big|\frac{1}{it}\int_{\frac{1}{t}}^{\delta}e^{-it\lambda}F^{'}(\lambda)d\lambda\Big|=\frac{1}{t^2}\Big|-e^{-i}F'(t^{-1})-\int_{\frac{1}{t}}^{\delta}e^{-it\lambda}F^{''}(\lambda)d\lambda\Big|\lesssim \frac{1}{t^2}\Big(t^{1-\alpha}+\int_{\frac{1}{t}}^{\delta}\lambda^{\alpha-2} d\lambda\Big)\lesssim t^{-1-\alpha}.
\end{equation*}
Thus $|I_2(t)|\lesssim t^{-(1+\alpha)}$, and the desired estimate follows.
\end{proof}

\appendix
\section{Relation Between the Threshold Equation and Birman--Schwinger Equation
}\label{example}
In this appendix, we further clarify the relation between the threshold equation $H\phi=\lambda\phi$ and the Birman--Schwinger equation 
$
    \phi+(-\Delta-\lambda\mp i0)^{-1}V\phi=0
$
at the hyperbolic thresholds $\lambda=4,8$ in the space $\ell^{2,-1}(\Z^3)$. 

As noted in Remark \ref{remark on defini of resonances}, every solution of the Birman--Schwinger equation in $\ell^{2,-1}(\mathbb Z^3)$ also solves the corresponding threshold equation, while the converse may fail in general. We first illustrate this distinction by constructing a nonzero solution $\phi\in\ell^{2,-1}(\mathbb Z^3)$ of the free threshold equation
$(-\Delta-\lambda)\phi=0$ that does not arise from the corresponding Birman--Schwinger equation. We then show that, under a suitable decay assumption on $V$, the two equations are equivalent in $\ell^2(\Z^3)$. 

Since the argument at the threshold 8 is analogous, we consider only the threshold 4. To this end, we first recall the following coarea formula on $\T^3$. 
\begin{lemma}\label{coarea formula lemma}
{\rm (\cite[Lemma 3.1]{NT20})} For $\psi\in C(\T^3)$, we have
\begin{equation*}
    \int_{\T^3}\psi(\xi)d\xi=\int_{0}^{12}\int_{\mcaM_{\lambda}}\psi(\xi)\frac{d\sigma_\lambda(\xi)}{|\nabla\mcaM(\xi)|}d\lambda, 
\end{equation*}
where 
$\mcaM_{\lambda}=\{\xi\in \T^3\mid\mcaM(\xi)=\lambda\}$ and $d\sigma_{\lambda}(\xi)$ denotes the induced surface measure of $\mcaM_{\lambda}$.
\end{lemma}

\begin{proposition}\label{counterexa}
    There exists $ 0\neq\phi\in\ell^{2,-1}(\Z^3)$ such that $(\Delta+4)\phi=0$.
\end{proposition}
\begin{proof}
  Let \begin{equation}\label{Sigma4}
      \Sigma_4:=\{x\in \mcaM_4\mid(\nabla\mcaM)(x)=0\}.
  \end{equation}
Fix $x_0=(x^{(1)}_0,x^{(2)}_0,x^{(3)}_0)\in \mcaM_4\setminus\Sigma_4.$ Without loss of generality, assume that $(\partial^{}_{x_3}\mcaM)(x_0)\neq0$. We introduce the local coordinate $t=\mcaM(x)-4$. Let $y=(x_1,x_2)$, $y_0=(x^{(1)}_0,x^{(2)}_0)$ and define  $$G(t,y,x_3):=\mcaM(y,x_3)-4-t,\quad t\in \R, \quad y\in \R^2,\quad x_3\in\R.$$
Then 
 $$
G(0,y_0,x^{(3)}_0)=0,\qquad (\partial^{}_{x_3}G)(0,y_0,x^{(3)}_0)=(\partial^{}_{x_3}\mcaM)(x_0)\neq0.
$$
By the implicit function theorem, there exist neighborhoods $I$ of $0$ and $D(y_0)$ of $y_0$, together with a unique function $\psi(t,y)\in C^{\infty}(I\times D(y_0))$, such that $\psi(0,y_0)=x^{(3)}_0$ 
 and $x_3=\psi(t,y)$. Then 
 $$x=(y,\psi(t,y))=:\Psi(t,y)\in \Psi(I\times D(y_0)).$$
 
 Since $(\partial^{}_{x_3}\mcaM)(x_0)\neq0$, there exists $(-\delta,\delta)\times D\subseteq I\times D(y_0)$ such that $x_0\in U:=\Psi((-\delta,\delta)\times D)$ and $\nabla\mcaM(x)\neq0$ for $x\in U$. Choose $\rho\in C^{\infty}_0((-\delta,\delta))$ with $\rho(0)=1$ and $0\not\equiv b(y)\in C^{\infty}_0(D)$. Let $a(\Psi(t,y))=\rho(t)b(y).$ Then $a\in C^{\infty}_0(U)$ and  $a|_{\mcaM_4}=b(y)\not\equiv 0$. We define the distribution $\mcaF\phi$ on $\T^3$ by  
\begin{equation}\label{Fphi expre as M4}
   (\mcaF\phi)(\varphi)=\int_{\mcaM_4}a(\xi)\varphi(\xi)\frac{d\sigma_{4}(\xi)}{|\nabla\mcaM(\xi)|},\qquad \varphi\in C^{\infty}(\T^3).  
\end{equation}
This is well defined since $\operatorname{supp}a\subseteq U$ and $\nabla\mcaM\neq0$ on $U$. Moreover,  $\mcaF\phi\neq0$ and $(\Delta+4)\phi=0$ since 
$$(\mcaF\phi)(\overline{a})>0,\qquad(\mcaM(\cdot)-4)(\mcaF\phi)(\varphi)=(\mcaF\phi)((\mcaM(\cdot)-4)\varphi(\cdot))=0.$$ 

It remains to show that $\phi\in \ell^{2,-1}(\Z^3)$, or equivalently, $\mcaF\phi\in H^{-1}(\T^3)$. In the local coordinates $\xi=\Psi(t,y)$, we write $d\xi=J(t,y)dtdy$, where $J$ is the Jacobian. Then for any $\varphi(\xi)\in C^{\infty}_0(U)$, by Lemma \ref{coarea formula lemma},
\begin{equation*}
 (\mcaF\phi)(\varphi)=\int_{D}a(\Psi(0,y))J(0,y)\varphi(\Psi(0,y))dy=:\int_{D}\widetilde{a}(y)\widetilde{\varphi}(0,y)dy,
\end{equation*}
where $\widetilde{a}(y)=a(\Psi(0,y))J(0,y)$ and $\widetilde{\varphi}(t,y)=\varphi(\Psi(t,y))$. Thus, under the smooth coordinate map $\Psi$, the distribution $\mcaF\phi$ corresponds locally to $$W(t,y)=\widetilde{a}(y)\delta_0(t).$$
It suffices to prove $\widetilde{a}(y)\delta_0(t)\in H^{-1}(\R^3)$. Indeed, taking the Fourier transform $\mcaF_c$ on $\R^3$, we obtain 
\begin{equation*}
(\mcaF_c{W})(\tau,\eta)=C(\mcaF_c\widetilde{a})(\eta)\in \mcaS(\R^2),\qquad \tau\in\R,\qquad \eta\in\R^2.
\end{equation*}
Then 
\begin{equation*}
    \|W\|^2_{H^{-1}(\R^3)}=\int_{\R^2}|(\mcaF_c{\widetilde{a}})(\eta)|^2\int_{\R}(1+\tau^2+|\eta|^2)^{-1}d\tau d\eta<\infty.
\end{equation*}
Therefore, $ 0\neq\phi\in\ell^{2,-1}(\Z^3)$ and $(\Delta+4)\phi=0$, as claimed.
\end{proof}
\vskip0.1cm
\begin{proposition}\label{propo of equi betw Hx and Birman eq}
Let $H=-\Delta+V$ on $\Z^3$ and $|V(n)|\lesssim \left<n\right>^{-\beta}$ for $\beta>\frac{5}{2}$. If $\phi\in\ell^2(\Z^3)$ satisfies $H\phi=4\phi$, then $\phi+(-\Delta-4\mp i0)^{-1}V\phi=0.$ 
\end{proposition}
\begin{proof}
Since $\phi\in\ell^2(\Z^3)$ and $H\phi=4\phi$, taking the Fourier transform gives
\begin{equation}\label{m-4 and f}
(\mcaM(\xi)-4)(\mcaF\phi)(\xi)=-\mcaF(V\phi)(\xi).
\end{equation}
Let $f=\mcaF(V\phi)$. Choose $\frac{5}{2}<s\leq\beta$. Then $f\in H^s(\T^3)$, and hence $f\in C^1(\T^3)$ by the Sobolev embedding theorem. We claim that $f|_{\mcaM_4}=0$. Once this holds, then by Lemma \ref{Taira lemma} (2), we have 
\begin{equation*}
    (\mcaM(\xi)-4\mp i0)^{-1}f(\xi)=(\mcaM(\xi)-4)^{-1}f(\xi),
\end{equation*}
which combined with \eqref{m-4 and f} gives $\phi+(-\Delta-4\mp i0)^{-1}V\phi=0.$ 

To prove $f|_{\mcaM_4}=0$, we first consider $\xi_0\in \mcaM_4\setminus\Sigma_4$, where $\Sigma_4$ is defined in \eqref{Sigma4}. Since $(\nabla\mcaM)(\xi_0)\neq0$, there exist a neighborhood $U(\xi_0)$ of $\xi_0$ and a constant $c(\xi_0)>0$, such that $|(\nabla\mcaM)(\xi)|\geq c(\xi_0)>0$ for all $\xi\in U(\xi_0)$. 

Choose an open neighborhood $U_1$ of $\xi_0$ such that $\overline{U_1}\subset U(\xi_0)$, and take $0\leq \varphi\in C_0^\infty(U(\xi_0))$ satisfying $\varphi\equiv1\  \text{on } U_1.$  Let $0\leq \rho\in C^{\infty}_0((-1,1))$ with $\int_{\R}\rho(s)ds=1$ and 
$\rho_{\varepsilon}(s)=\varepsilon^{-1}\rho(s\varepsilon^{-1})$. By Lemma \ref{coarea formula lemma},
\begin{equation*}
\int_{\T^3}\varphi(\xi)|f(\xi)|^2\rho_{\varepsilon}(\mcaM(\xi)-4)d\xi=\int_{0}^{12}\rho_{\varepsilon}(\lambda-4)\int_{\mcaM_{\lambda}}\frac{\varphi(\xi)|f(\xi)|^2d\sigma_{\lambda}(\xi)}{|(\nabla\mcaM)(\xi)|}d\lambda\rightarrow\int_{\mcaM_{4}}\frac{\varphi(\xi)|f(\xi)|^2d\sigma_{4}(\xi)}{|(\nabla\mcaM)(\xi)|}
\end{equation*}
as $\varepsilon\rightarrow0$. On the other hand, by \eqref{m-4 and f} and the uniform boundedness of $s^2\rho(s)$, we have
\begin{equation*}
\Big|\int_{\T^3}\varphi(\xi)|f(\xi)|^2\rho_{\varepsilon}(\mcaM(\xi)-4)d\xi\Big|=\Big|\int_{\T^3}\varphi(\xi)(\mcaM(\xi)-4)^2|(\mcaF\phi)(\xi)|^2\rho_{\varepsilon}(\mcaM(\xi)-4)d\xi\Big|\lesssim  \varepsilon\rightarrow0.
\end{equation*}
Consequently, the continuity of $f$ implies that $f(\xi)=0$ on $\mcaM_4\cap \{\varphi(\xi)>0\}$. In particular, since $\varphi(\xi_0)=1>0$, it follows that $f(\xi_0)=0$. It remains to consider $\xi\in\Sigma_4$. By the Morse lemma, there exist neighborhoods $U(\xi)$ of $
\xi$ and $V_0$ of $0$, and a $C^{\infty}$ diffeomorphism $G: U(\xi)\longrightarrow V_0$ such that
\begin{equation*}
    G(\xi)=0\quad {\rm and}\quad \mcaM(G^{-1}(\eta))-4=-|\eta'|^2+|\eta''|^2=:F(\eta), \quad \eta=(\eta',\eta''),\quad \eta'\in \R,\quad \eta''\in \R^2.
    \end{equation*}
Choose unit vector $e'\in\R$ and $e''\in\R^2$, and set $\eta_r=(re',re'')$ with $r>0$. Then   
$$F(\eta_r)=0,\quad (\nabla F)(\eta_r)=(-2re',2re'')\neq 0,\quad  \lim\limits_{r\rightarrow0}\eta_r=0.$$
Define $\xi_r:=G^{-1}({\eta_r})$. Then $\xi_r\in \mcaM_4\setminus\Sigma_4$ by  $\mcaM(\xi_r)=4+F(\eta_r)=4$ and $(\nabla\mcaM)(\xi_r)\neq0$. Moreover, 
\begin{equation*}
\quad \xi=G^{-1}(0)=\lim\limits_{r\rightarrow0}G^{-1}(\eta_r)=\lim\limits_{r\rightarrow0}\xi_r.
\end{equation*}
By the result above and the continuity of $f$, we have $f(\xi)=\lim\limits_{r\rightarrow0}f(\xi_r)=0$. Therefore, $f|_{\mcaM_4}=0$.
\end{proof}
\section{Kernel Coefficients in the Free Resolvent Expansions}\label{further infor of coeffi}
This appendix provides further details and explicit computations related to the free resolvent expansions discussed in Remark \ref{remark of free expan}. We compute the required coefficients $A_{2k+1}$ near threshold 0 using the exact discrete symbol, and at the threshold 4, we isolate the half-integer contribution and derive the kernel $B^\pm_0(n)$.

\begin{proposition}\label{prop of Ak calcul}
The expansion \eqref{expan of R0} holds. Moreover, $(-\Delta)A_0=\delta_0$, and the kernel $A_k(n-m)$ of $A_k$ for odd $k$ is given by
\begin{equation}\label{Ak for odd k}
 A_k(n)=\frac{i(\frac{k-1}{2})!}{4\pi k!}\sum\limits_{\substack{k_1+k_2+k_3=\frac{k-1}{2}\\  k_1,k_2,k_3\in \N}} \prod\limits_{j=1}^{3}\frac{(\frac{1}{2}-n_j)_{k_j}(\frac{1}{2}+n_j)_{k_j}}{k_j!},\quad n=(n_1,n_2,n_3)\in\Z^3
\end{equation}
with $(a)_{\ell}:=a(a+1)\cdots(a+\ell-1)$ for integer $\ell\geq1$ and $ (a)_0=1$. As a consequence, $A_1(n)=\frac{i}{4\pi}$.
\end{proposition} 

\begin{proof}
Fix $0<\delta<1$ and let $0<|\omega|<\delta/2$ with $0<\arg\omega<2\pi$. We decompose $R_0(\omega,n)$ as 
\begin{equation}\label{decom of R0}
R_0(\omega,n)=\frac{1}{8\pi^3}\Big(\int_{\{\mcaM(\theta)<\delta\}}+\int_{\{\mcaM(\theta)\geq\delta\}}\Big)\frac{e^{-i\theta\cdot n}}{\mcaM(\theta)-\omega}d\theta=:R_{01}(\omega,n)+R_{02}(\omega,n).
\end{equation}
Clearly, $R_{02}(\omega,n)$ is analytic in $0<|\omega|<\delta/2$. Moreover, for every $k\in\N$, the derivative $(\partial^{k}_\omega R_{02})(\omega,n)$ is uniformly bounded in $n$ and $0<|\omega|<\delta/2$. It remains to analyze $R_{01}(\omega,n)$. Making the change of variables 
\begin{equation*}
s_j=2\sin\frac{\theta_j}{2}, \qquad g(s_j)=2\arcsin{\frac{s_j}{2}},\qquad s=(s_1,s_2,s_3),
\end{equation*}
and using the polar coordinates $s=rs'$, $s'\in\mathbb{S}^2$, we have 
\begin{equation}\label{expre of R01}
    R_{01}(\omega,n)=\frac{1}{8\pi^3}\int_{0}^{\sqrt{\delta}}\frac{r^2}{r^2-\omega}F(r,n)dr,\qquad F(r,n)=\int_{\mathbb{S}^2}\prod\limits_{j=1}^{3}e^{-ig(rs'_j)n_j}g'(rs'_j)d\sigma(s'),
\end{equation}
where $d\sigma(s')$ denotes the surface measure on the unit sphere $\mathbb{S}^2$.
By symmetry, we note that
\begin{equation*}
F(r,n)=\int_{\mathbb{S}^2}G(r,s,n)d\sigma(s),\quad G(r,s,n)=\prod\limits_{j=1}^{3}g_{n_j}(rs_j),\quad g_{m}(x)={\rm cos}(mg(x))g'(x).
\end{equation*}
For any integer $M$, applying Taylor expansion to $G(r,s,n)$ near $r=0$, we obtain
\begin{equation*}
G(r,s,n)=\sum\limits_{k=0}^{M}\frac{1}{k!}(\partial^k_rG)(0,s,n)r^k+G_{M}(r,s,n).
\end{equation*}
Since $G$ is an even function about $r$, it follows that $(\partial^k_rG)(0,s,n)=0$ if $k$ is odd. Moreover, 
\begin{equation*}
G_{M}(r,s,n)=\begin{cases}
    \frac{r^{M+1}}{M!}\int_{0}^{1}(1-t)^{M}(\partial^{M+1}_rG)(tr,s,n)dt,&M\ is\ odd,\\
    \frac{r^{M+2}}{(M+1)!}\int_{0}^{1}(1-t)^{M+1}(\partial^{M+2}_rG)(tr,s,n)dt,&M\ is\ even.
\end{cases} 
\end{equation*}
Consequently, 
\begin{equation}\label{expan of Frn}
F(r,n)=\sum\limits_{k=0}^{M}f_{k}(n)r^k+E_{M}(r,n), \quad 0<r<\sqrt{\delta}
\end{equation}
with 
\begin{equation*}
f_k(n)=\frac{1}{k!}\int_{\mathbb{S}^2}(\partial^k_rG)(0,s,n)d\sigma(s),\qquad  E_{M}(r,n)=\int_{\mathbb{S}^2}G_M(r,s,n)d\sigma(s).
\end{equation*}
It is clear that $f_k(n)=0$ if $k$ is odd. Next we claim that
\begin{equation}\label{fk for even k}
f_k(n)=\frac{4\pi(\frac{k}{2})!}{(k+1)!}\sum\limits_{\substack{k_1+k_2+k_3=\frac{k}{2}\\  k_1,k_2,k_3\in \N}} \prod\limits_{j=1}^{3}\frac{(\frac{1}{2}-n_j)_{k_j}(\frac{1}{2}+n_j)_{k_j}}{k_j!},\qquad \text{if $k$ is even}, 
\end{equation}
and $E_{M}(r,n)$ satisfies the following estimate
\begin{equation}\label{estimate of fN21}
    |E_{M}(r,n)|\lesssim \begin{cases}
r^{M+1}\left<n\right>^{M+1}, &M\ is\ odd,\\
r^{M+2}\left<n\right>^{M+2}, &M\ is\ even,
    \end{cases}\qquad 0<r<\sqrt{\delta}.
\end{equation}
Indeed, first, by Leibniz's rule, 
\begin{equation}\label{leibniz of G}
(\partial^k_rG)(r,s,n)=\sum\limits_{\substack{k_1+k_2+k_3=k\\  k_1,k_2,k_3\in \N}} \frac{k!}{k_1!k_2!k_3!}\prod\limits_{j=1}^{3}s_j^{k_j}g^{(k_j)}_{n_j}(rs_j).
\end{equation}
We note that $
(4-x^2)g''_m(x)-3xg'_m(x)+(4m^2-1)g_m(x)=0$ for $ |x|<2.
$
This together with 
\begin{equation*}
g_m(x)=\sum\limits_{j=0}^{+\infty}b_j(m)x^{2j},\quad b_j(m)=\frac{g^{(2j)}_m(0)}{(2j)!},
\end{equation*}
yields $
    4(2j+2)(2j+1)b_{j+1}(m)+(4m^2-(2j+1)^2)b_j(m)=0$ and thus 
    $$g^{(2(j+1))}_{m}(0)=((j+\frac{1}{2})^2-m^2)g^{(2j)}_{m}(0)\Longrightarrow g^{(2j)}_{m}(0)=(\frac{1}{2}-m)_{j}(\frac{1}{2}+m)_{j}.$$
Hence, if $k$ is even, using \eqref{leibniz of G} and $$\int_{\mathbb{S}^2}s^{2k_1}_1s^{2k_2}_2s^{2k_3}_3d\sigma(s)=\frac{4\pi}{(k+1)!!}\prod\limits_{j=1}^{3}(2k_j-1)!!,$$ where $k_1+k_2+k_3=\frac{k}{2}$, we derive \eqref{fk for even k}.
To estimate the remainder term, since
$$g_m(x)=\frac{1}{2}(e^{img(x)}+e^{-img(x)})g'(x),\qquad(e^{\pm img(x)})^{(k)}(x)=e^{\pm img(x)} P_{k,m}(x),$$
where $P_{k,m}(x)$ is a polynomial in $m$ of degree at most $k$, and each of its coefficients is a finite product and linear combination of $\{g^{(j)}(x)\}^k_{j=0}$. This together with \eqref{leibniz of G} implies \eqref{estimate of fN21}. Taking \eqref{expan of Frn} into \eqref{expre of R01} yields
\begin{equation}\label{expre2 of R01}
R_{01}(\omega,n)=\frac{1}{8\pi^3}\sum\limits_{k=0}^{M}f_k(n)\int_{0}^{\sqrt{\delta}}\frac{r^{k+2}}{r^2-\omega}dr+\frac{1}{8\pi^3}\int_{0}^{\sqrt{\delta}}\frac{r^{2}E_{M}(r,n)}{r^2-\omega}dr.
\end{equation}
If $k$ is even, using 
$r^{k+2}=r^{2(\frac{k}{2}+1)}-\omega^{\frac{k}{2}+1}+\omega^{\frac{k}{2}+1}$, we derive 
\begin{equation}\label{fk rk integral}
f_k(n)\int_{0}^{\sqrt{\delta}}\frac{r^{k+2}}{r^2-\omega}dr=\frac{\pi i}{2}\omega^{\frac{k+1}{2}}f_k(n)+\Big(\underbrace{\sum\limits_{j=0}^{\frac{k}{2}}a_j(k)\omega^j+\frac{1}{2}\omega^{\frac{k+1}{2}}\ln\frac{1-\sqrt{\frac{\omega}{\delta}}}{1+\sqrt{\frac{\omega}{\delta}}}}_{=:P_k(\omega)\ \text{is analytic in}\ \omega }\Big)f_k(n).
\end{equation}
Since $f_k(n)=0$ when $k$ is odd. Thus, this formula is valid for any $k\in\N$.
\vskip0.1cm
We now analyze the remainder integral in \eqref{expre2 of R01}. Fix a sufficiently small $0<\varepsilon\ll1$. We further restrict $\omega$ to $\varepsilon\leq \arg\omega\leq 2\pi-\varepsilon$. Let $q=[\frac{M+1}{2}]$. Using 
\begin{equation*}
    \frac{r^2}{r^2-\omega}=\sum\limits_{j=0}^{q}\omega^j r^{-2j}+\omega^{q+1}\frac{r^{-2q}}{r^2-\omega},
\end{equation*}
and the estimate \eqref{estimate of fN21}, we write
\begin{equation}\label{expre of remainder}
\int_{0}^{\sqrt{\delta}}\frac{r^{2}E_{M}(r,n)}{r^2-\omega}dr=Q_M(\omega,n)+\omega^{q+1}\int_{0}^{\sqrt{\delta}}\frac{r^{-2q}}{r^2-\omega}E_M(r,n)dr,
\end{equation}
where $Q_M(\omega,n)$ is a polynomial in $\omega$. For the second term, we claim that 
\begin{equation}\label{estimate for w qplus1}
\Big|\omega^{q+1}\int_{0}^{\sqrt{\delta}}\frac{r^{-2q}}{r^2-\omega}E_M(r,n)dr\Big|=O(|\omega|^{\frac{M+2}{2}}\left<n\right>^{M+1}). 
\end{equation}
Indeed, if $M$ is odd, then $q=\frac{M+1}{2}$. Since $\varepsilon\leq \arg\omega\leq 2\pi-\varepsilon$, we have the following estimate 
\begin{equation*}
    |r^2-\omega|\gtrsim r^2+|\omega|.
\end{equation*}
Combining this with \eqref{estimate of fN21}, we obtain 
\begin{equation*}
\Big|\omega^{q+1}\int_{0}^{\sqrt{\delta}}\frac{r^{-2q}}{r^2-\omega}E_M(r,n)dr\Big|\lesssim |\omega|^{q+1}\int_{0}^{\sqrt{\delta}}\frac{dr}{r^2+|\omega|}\left<n\right>^{M+1}\lesssim |\omega|^{\frac{M+2}{2}}\left<n\right>^{M+1}.
\end{equation*}
If $M$ is even, then $q=\frac{M}{2}$. Let $C_{n,\delta}=\min\{\left<n\right>^{-1},\sqrt{\delta}\}$. 
We divide the integral into two parts
\begin{equation*}
\int_{0}^{\sqrt{\delta}}\frac{r^{-2q}}{r^2-\omega}E_M(r,n)dr=\Big(\int_{0}^{C_{n,\delta}}+\int_{C_{n,\delta}}^{\sqrt{\delta}}\Big)\frac{r^{-M}}{r^2-\omega}E_M(r,n)dr=:(I_1+I_2)(M,\omega,n).
\end{equation*}
It is clear that 
\begin{equation*}
    |I_1(M,\omega,n)|\lesssim \left<n\right>^{M+2}\int_{0}^{C_{n,\delta}}\frac{r^2}{r^2+|\omega|}dr\lesssim \left<n\right>^{M+1}.
\end{equation*}
Since $F(r,n)$ is uniformly bounded in $r\in [0,\sqrt{\delta}]$ and $n\in\Z^3$, and $|f_{2j}(n)|\lesssim \left<n\right>^{2j}$,  when $r\left<n\right>\geq1$, we have
\begin{equation*}
    |E_{M}(r,n)|=\Big|F(r,n)-\sum\limits_{j=0}^{q}f_{2j}(n)r^{2j}\Big|\lesssim r^{M}\left<n\right>^{M}.
\end{equation*}
Then 
\begin{equation*}
|I_2(M,\omega,n)|\lesssim \left<n\right>^{M}\int_{C_{n,\delta}}^{\sqrt{\delta}}\frac{dr}{r^2}\leq \left<n\right>^{M}\int_{\left<n\right>^{-1}}^{+\infty}\frac{dr}{r^2}\leq \left<n\right>^{M+1}.
\end{equation*}
Therefore, we obtain \eqref{estimate for w qplus1}. Combining \eqref{fk rk integral}, \eqref{expre of remainder}, \eqref{estimate for w qplus1}, \eqref{expre2 of R01} and \eqref{decom of R0}, we obtain 
\begin{align*}
    R_{0}(\omega,n)=\sum\limits_{k=1}^{M+1}\frac{if_{k-1}(n)}{16\pi^2}\omega^{\frac{k}{2}}+\frac{1}{8\pi^3}\Big(\sum\limits_{k=0}^{M}P_k(\omega)f_k(n)+Q_M(\omega,n)\Big)+R_{02}(\omega,n)+O(|\omega|^{\frac{M+2}{2}}\left<n\right>^{M+1}),
\end{align*}
which gives \eqref{expan of R0} by letting $N=M+1$. Moreover, if $k$ is odd,  $A_k(n)=\frac{if_{k-1}(n)}{16\pi^2}$. This combined with \eqref{fk for even k} yields \eqref{Ak for odd k}. 
As a consequence,  $R_0(\omega,n)=A_0(n)+\widetilde{A}(\omega,n)$ with $|\widetilde{A}(\omega,n)|=O(|\omega|^{1/2})$. We first fix $\omega=a+ib~(b\neq0)$ and $|\omega|<\delta/2$. Since for any $n\in\Z^3$, $\big((-\Delta-\omega)R_0(\omega,\cdot)\big)(n)=\delta_0(n)$. Hence, 
\begin{equation*}
\delta_0(n)+\omega R_0(\omega,n)=\big((-\Delta)R_0(\omega,\cdot)\big)(n)=\big((-\Delta)A_0(\cdot)\big)(n)+\big((-\Delta)\widetilde{A}(\omega,\cdot))(n).  
\end{equation*}
Letting $\omega\rightarrow0$, we then derive $\big((-\Delta)A_0(\cdot)\big)(n)=\delta_0(n)$ and prove that $(-\Delta)A_0=\delta_0$.
\end{proof}
\vskip0.2cm
\begin{proposition} \label{prop of B0 calcu}
Let $B^\pm_0$ and $E^\pm_0$ be defined in \eqref{expan at 4}. Then we have
\begin{equation*}
B^\pm_0(n)=\frac{\mp}{4\pi}\big((-1)^{n_1}+(-1)^{n_2}+(-1)^{n_3}\big),\quad (-\Delta-4)E^\pm_0=\delta_0.
\end{equation*}
\end{proposition}
\begin{proof} 
Fix $0<\delta<1$, and let $0<|\omega|<\min\{\frac{\delta}{2},\mu\}$, where $\mu$ will be specified later. By the conjugation symmetry, it suffices to consider the case $\operatorname{Im}\omega>0.$  Choose a partition of unity $$\sum\limits_{j=1}^{3}\xi_j(\theta)+\chi(\theta)+\eta(\theta)=1,$$
where the functions  $\xi_j(\theta)$ are later determined so that $\xi_j(\theta)=1$ in a neighborhood of $\pi e_j$ and
$$\operatorname{supp}\chi \subseteq \{|\mcaM(\theta)-4|<\delta\} \setminus \bigcup_{j=1}^{3}\operatorname{supp}\xi_j, \qquad \operatorname{supp}\eta \subseteq \{|\mcaM(\theta)-4|\geq\delta\}.$$
Accordingly, we decompose
\begin{equation*}
R_0(4+\omega,n)=\sum\limits_{j=1}^{3}\underbrace{\int_{\T^3}\frac{e^{-i\theta\cdot n}\xi_j(\theta)}{8\pi^3(\mcaM(\theta)-4-\omega)}d\theta}_{:=K_j(\omega,n)}+\underbrace{\int_{\T^3}\frac{e^{-i\theta\cdot n}\chi(\theta)}{8\pi^3(\mcaM(\theta)-4-\omega)}d\theta}_{:=K_r(\omega,n)}+\underbrace{\int_{\T^3}\frac{e^{-i\theta\cdot n}\eta(\theta)}{8\pi^3(\mcaM(\theta)-4-\omega)}d\theta}_{:=K_a(\omega,n)}.
\end{equation*}

Clearly, $K_a(\omega,n)$ is analytic in $\omega$ and hence contributes only integer powers of $\omega$. We next show that the same holds for $K_r(\omega,n)$. Consequently, all half-integer contributions arise from the terms $K_j(\omega,n)$. Following the argument in the proof of \cite[Proposition 2.6]{Kop09}, we then compute the leading half-integer coefficient $B_0^+(n)$.
\vskip0.2cm
\underline{ \bf (1) The analyis of $K_r$.} Choose a finite partition of unity 
$$\chi(\theta)=\sum_{j=1}^{N_0}\chi_j(\theta), \qquad \chi_j\in C^\infty_0(\T^3),$$
with each $\operatorname{supp}\chi_j$ sufficiently small and contained in $\operatorname{supp}\chi$. Then
\begin{equation*}
    K_r(\omega,n)=\frac{1}{8\pi^3}\sum\limits_{j=1}^{N_0}\int_{\T^3}\frac{e^{-i\theta\cdot n}\chi_j(\theta)}{\mcaM(\theta)-4-\omega}d\theta=:\frac{1}{8\pi^3}\sum\limits_{j=1}^{N_0}K_{rj}(\omega,n),\quad \chi_j\in C^{\infty}_0({\rm supp}\chi).
\end{equation*}
Let $S_j=\operatorname{supp}\chi_j\cap\mcaM_4$. If $S_j=\varnothing$, then $K_{rj}(\omega,n)$ is clearly analytic in $\omega$. Suppose therefore that $S_j\neq\varnothing$. Then any $\theta\in\operatorname{supp}\chi_j$ can be represented uniquely as $$\theta=s+t\mathbf{n}(s),\qquad  s\in S_j,$$
where $\mathbf{n}(s):=(\nabla\mcaM)(s)/|(\nabla\mcaM)(s)|$ is the outward unit normal vector at $s$. By the mean value theorem, 
\begin{equation*}
\mcaM(\theta)-4=\mcaM(s+t\mathbf{n}(s))-\mcaM(s)=t\int_{0}^{1}\mathbf{n}(s)\cdot(\nabla\mcaM)(s+t\rho\mathbf{n}(s)) d\rho=:u(s,t).
\end{equation*}
Denote by $d\sigma(s)$ the surface measure on $S_j$. Then 
\begin{equation*}
    K_{rj}(\omega,n)=\int_{S_j}d\sigma(s)\int_{-a(s)}^{b(s)}\frac{e^{-i(s+t\mathbf{n}(s))\cdot n}\chi_j(s+t\mathbf{n}(s))J(s,t)}{u(s,t)-\omega}dt,
\end{equation*}
where $J(s,t)$ is the Jacobian and $a(s),b(s)>0$. Since $(\partial_tu)(s,0)=|(\nabla\mcaM)(s)|>0$ for any $s\in S_j$. Define 
$$\gamma_j:=\frac{1}{2}\min\limits_{s\in S_j}|(\nabla\mcaM)(s)|.$$
By the compactness of $S_j$, there exists $\tau_j>0$, independent of $s$, such that $[-\tau_j,\tau_j]\subseteq(-a(s),b(s))$ and 
$$(\partial_tu)(s,t)\geq \gamma_j>0,\qquad s\in S_j,\quad |t|\le\tau_j.$$
Denote 
$$\widetilde{\chi}_j(s,t,n)=e^{-i(s+t\mathbf{n}(s))\cdot n}\chi_j(s+t\mathbf{n}(s))J(s,t).$$
Let
$$\mathcal I_j := \bigl\{(s,t):s\in S_j,\, t\in(-a(s),b(s))\setminus[-\tau_j,\tau_j]\bigr\} \cap\operatorname{supp}\widetilde\chi_j.$$
Since $u(s,t)\neq0$ on $\mathcal I_j$, compactness yields
$$m_j:=\min\limits_{(s,t)\in\mcaI_j}|u(s,t)|>0.$$ Choose 
$$0<a_j<\frac12\min\{m_j,\gamma_j\tau_j\},$$
and $\zeta\in C^{\infty}_0((-a_j,a_j))$, $0\leq \zeta\leq1$ and $\zeta(x)=1$ for $|x|\leq \frac{a_j}{2}$. Let $0<|\omega|<\min\limits_{1\leq j\leq N_0}\{\frac{a_j}{4}\}=:\mu$. We decompose 
\begin{align*}
    K_{rj}(\omega,n)&=\int_{S_j}d\sigma(s)\int_{-a(s)}^{b(s)}\frac{\zeta(u(s,t))\widetilde{\chi}_j(s,t,n)}{u(s,t)-\omega}dt+\int_{S_j}d\sigma(s)\int_{-a(s)}^{b(s)}\frac{(1-\zeta(u(s,t)))\widetilde{\chi}_j(s,t,n)}{u(s,t)-\omega}dt.
\end{align*}
Since supp $(1-\zeta(x))\subseteq\{|x|\geq \frac{a_j}{2}\}$, and therefore the second term is analytic in $\omega$. It remains to consider
\begin{equation*}
    \widetilde{K}_{rj}(\omega,n):=\int_{S_j}d\sigma(s)\int_{-a(s)}^{b(s)}\frac{\zeta(u(s,t))\widetilde{\chi}_j(s,t,n)}{u(s,t)-\omega}dt
\end{equation*}
By the choice $a_j<m_j/2$, the condition
$\zeta(u(s,t))\neq0$ implies $|t|<\tau_j$. Moreover, our choice $a_j<\gamma_j\tau_j/2$ guarantees that the inverse
$t=t(s,u)$ is well-defined for every $s\in S_j$ and $|u|<a_j.$ Therefore, by the change of variables $u=u(s,t)$, we obtain
\begin{equation*}
    \widetilde{K}_{rj}(\omega,n)=\int_{-a_j}^{a_j}\frac{1}{u-\omega}\int_{S_j}\zeta(u)\widetilde{\chi}_j(s,t(s,u),n)(\partial_ut)(s,u)d\sigma(s)du=:\int_{-a_j}^{a_j}\frac{A_j(u,n)}{u-\omega}du.
\end{equation*}
By definition, for each $n$, 
 $A_j(\cdot,n)\in C^{\infty}_0((-a_j,a_j))$, and for every $k\in\N$, 
 \begin{equation}\label{estima of parkuA_j}
  \sup\limits_{|u|\leq  a_j}|(\partial^k_uA_j)(u,n)|\lesssim\left<n\right>^{k}.   
 \end{equation}
 Since Im$\omega>0$, 
 $$\partial_\omega^k\frac1{u-\omega} = \frac{k!}{(u-\omega)^{k+1}}, \qquad \partial_u^k\frac1{u-\omega} = (-1)^k\frac{k!}{(u-\omega)^{k+1}}.$$
 Hence, integrating by parts $k$ times and using
$A_j(\cdot,n)\in C^\infty_0((-a_j,a_j))$, we obtain
\begin{equation*}
(\partial^k_{\omega}\widetilde{K}_{rj})(\omega,n)=k!\int_{-a_j}^{a_j}\frac{A_j(u,n)}{(u-\omega)^{k+1}}du=(-1)^k\int_{-a_j}^{a_j}A_j(u,n)\partial^k_u(\frac{1}{u-\omega})du=\int_{-a_j}^{a_j}\frac{(\partial^k_uA_j)(u,n)}{u-\omega}du.
\end{equation*}
For $0<\varepsilon\leq1$, let
$K(\varepsilon):=\widetilde{K}_{rj}(\varepsilon\omega,n)$. Applying Taylor's expansion to $K$ at $\varepsilon$, evaluated at 1, gives
\begin{align*}
  \widetilde{K}_{rj}(\omega,n)= \sum\limits_{k=0}^{N}\frac{\omega^k}{k!}(\partial^k_{\omega}\widetilde{K}_{rj})(\varepsilon\omega,n)(1-\varepsilon)^k+\frac{(\omega(1-\varepsilon))^{N+1}}{N!}\int_{0}^{1}(1-t)^{N}(\partial^{N+1}_{\omega}\widetilde{K}_{rj})((\varepsilon+t(1-\varepsilon))\omega,n)dt.
\end{align*}
Applying \eqref{estima of parkuA_j} and Lemma \ref{cauthy integral lemma} with $h=\min\{\frac{a_j}{4},\left<n\right>^{-1}\}$,we obtain
\begin{equation*}
    \sup\limits_{0<|\omega|<\mu,\ \operatorname{Im}\omega>0}|(\partial^k_{\omega}\widetilde{K}_{rj})(\omega,n)|\lesssim\left<n\right>^{k}(1+\ln(1+\left<n\right>)).
\end{equation*}
Applying this estimate to the remainder and letting $\varepsilon\rightarrow0$, we derive
\begin{align*}
    \widetilde{K}_{rj}(\omega,n)&=\sum\limits_{k=0}^{N}\frac{\omega^k}{k!}\Big((\partial^{k}_uA_j)(0,n)\pi i+\int_{-a_j}^{a_j}\frac{(\partial^{k}_uA_j)(u,n)-(\partial^{k}_uA_j)(0,n)}{u}du\Big)\\
&\quad+O\Big(|\omega|^{N+1}\left<n\right>^{N+1}(1+\ln(1+\left<n\right>))\Big).
\end{align*}
Therefore, we prove that $K_r(\omega,n)$ contributes only integer powers of $\omega$. We now turn to the contributions of $K_j$ and compute $B^+_0(n)$.
\vskip0.2cm
\underline{\bf(2) Calculation of $\bm {B^+_0(n)}$.} To this end, we briefly review and reorganize Kopylova's argument.
Let $$\theta'=\theta-\pi e_j,\qquad\eta_j(\theta')=\xi_j(\theta'+\pi e_j)$$ and choose $\eta_j$ to be even in each variable. Then
\begin{equation*}
K_j(\omega,n)=\frac{(-1)^{n_j}}{8\pi^3}\int_{\T^3}\frac{e^{-i\theta'\cdot n}\eta_j(\theta')d\theta'}{\mcaM(\theta'+\pi e_j)-4-\omega}=\frac{(-1)^{n_j}}{\pi^3}\int_{(0,+\infty)^3}\frac{\prod\limits_{k=1}^{3}{\rm cos}(n_k\theta'_k)\eta_j(\theta')}{\sum\limits_{k=1}^{3}4{\rm sin}^2\frac{\theta'_k}{2}-8{\rm sin}^2\frac{\theta'_j}{2}-\omega}d\theta'.
\end{equation*}
Since the three cases are analogous, it suffices to consider $K_1(\omega,n)$. 
\vskip0.1cm
{\bf Step 1}: We apply the following three changes of variables successively:
\begin{align*}
  &(1)~\theta'=(\theta'_1,\theta'_2,\theta'_3)\rightarrow s=(s_1,s_2,s_3): s_1=2{\rm sin}\frac{\theta'_2}{2},\quad s_2=2{\rm sin}\frac{\theta'_3}{2},\quad s_3=2{\rm sin}\frac{\theta'_1}{2},\\
  &(2)~~s=(s_1,s_2,s_3)\rightarrow (\tau,s_3):~~~~~~~~ s_1=\tau\cos\varphi,\quad s_2=\tau\sin\varphi,\quad s_3=s_3,\\
  &(3)~(\tau,s_3)\rightarrow \rho=(\rho_1,\rho_2):
  \rho_1=\tau^2-s^2_3,\quad \rho_2=2\tau s_3.
\end{align*}
We now specify the cutoff $\eta_1$ so that, under the first change of
variables, 
$$\eta_1(\theta')=\widetilde\eta_1(|s|^2), \qquad |s|^2=s_1^2+s_2^2+s_3^2,$$
where $\widetilde\eta_1\in C^\infty_0([0,\infty))$ is chosen with
sufficiently small support such that
$$\operatorname{supp}\bigl(\widetilde\eta_1(|\rho|)\bigr) \cap\{\rho_2\ge0\} \subset D_\delta := \{(\rho_1,\rho_2):|\rho_1|\le\delta,\ 0\le\rho_2\le\delta\},\qquad |\rho|=\sqrt{\rho^2_1+\rho^2_2}.$$
Then
\begin{align}
K_1(\omega,n)&=\int_{(0,+\infty)^3}\frac{\frac{(-1)^{n_1}}{\pi^3}\prod\limits_{k=1}^{3}{\rm cos}(n_k\theta'_k)\eta_1(\theta')d\theta'}{4\big({\rm sin}^2\frac{\theta'_2}{2}+{\rm sin}^2\frac{\theta'_3}{2}-{\rm sin}^2\frac{\theta'_1}{2}\big)-\omega}\xlongequal[]{(1)}\frac{(-1)^{n_1}}{\pi^3}\int_{(0,+\infty)^3}\frac{F_1(n,s^2_1,s^2_2,s^2_3)\widetilde{\eta}_1(|s|^2)}{s^2_1+s^2_2-s^2_3-\omega}ds\notag\\
&\xlongequal[]{(2)}\frac{(-1)^{n_1}}{\pi^3}\int_{(0,+\infty)^2}\frac{\tau\widetilde{\eta}_1(\tau^2+s^2_3)}{\tau^2-s^2_3-\omega}\mathbb{F}_1(n,\tau^2,s^2_3)d\tau ds_3,\notag\\
&\xlongequal[]{(3)}\frac{(-1)^{n_1}}{\pi^3}\int_{0}^{\delta}\int_{-\delta}^{\delta}\frac{\sqrt{r+\rho_1}}{(\rho_1-\omega)r}H_1(r,\rho_1,n)d\rho_1 d\rho_2,\quad r=|\rho|,\label{expre of K1 after transform}
\end{align}
where $g(x)=2{\arcsin}\frac{x}{2}$, $G(s)=\Big((1-\frac{s^2_1}{4})(1-\frac{s^2_2}{4})(1-\frac{s^2_3}{4})\Big)^{-1/2}$ and 
\begin{align*}
&F_1(n,s^2_1,s^2_2,s^2_3)={\rm cos}(n_1g(s_3)){\rm cos}(n_2g(s_1)){\rm cos}(n_3g(s_2))G(s),\\
&\mathbb{F}_1(n,\tau^2,s^2_3)=\int_{0}^{\frac{\pi}{2}}F_1(n,\tau^2\cos^2\varphi,\tau^2\sin^2\varphi,s^2_3)d\varphi,\quad
H_1(r,\rho_1,n)=\frac{\widetilde{\eta}_1(r)}{4\sqrt{2}}\mathbb{F}_1(n,\frac{r+\rho_1}{2},\frac{r-\rho_1}{2}).
\end{align*}
\vskip0.1cm
{\bf Step 2}: Let $N$ be an integer. Applying Taylor expansion to $H_1(r,\rho_1,n)$ near $\rho_1=0$, we derive
\begin{equation*}
H_1(r,\rho_1,n)=\sum\limits_{k=0}^{N}\frac{(\partial^{k}_{\rho_1}H_1)(r,0,n)}{k!}\rho^k_1+h_N(r,\rho_1,n)\rho^N_1,\quad |\rho_1|\leq \delta,
\end{equation*}
where $(\partial^{k}_{\rho_1}H_1)(r,0,n)$ is a polynomial in $n$ of order $2k$ and $$|(\partial^j_{\rho_1}h_N)(r,\rho_1,n)| \lesssim  \left<n\right>^{2N+2j}\ {\rm for\ } j=0,1.$$ 
Substituting this expansion into \eqref{expre of K1 after transform} and using \cite[(2.19)--(2.21), (2.29) and Lemma 2.8]{Kop09}, we obtain

\begin{align*}
K_1(\omega,n)&=\frac{(-1)^{n_1}}{\pi^3}\sum\limits_{k=0}^{N}\Big(\sum\limits_{j=0}^{k-1}a_{k,j}(n)\omega^j+\omega^k\int_{0}^{\delta}\Big(\pi i-\ln\frac{1-\sqrt{\frac{r-\omega}{2r}}}{1+\sqrt{\frac{r-\omega}{2r}}} \Big)\frac{(\partial^{k}_{\rho_1}H_1)(r,0,n)}{k!\sqrt{r-\omega}}dr\Big)\\
&\quad+\sum\limits_{j=0}^{N-1}b_j(N,n)\omega^j+O(\left<n\right>^{2N}(1+\ln\left<n\right> )^2|\omega|^N).
\end{align*}
Applying Taylor's expansion once more to each $(\partial_{\rho_1}^k H_1)(r,0,n)$ near $r=0$, up to order $N-k$, and using \cite[Lemma 2.7 and (2.24)--(2.28)]{Kop09}, we further obtain
\begin{align*}
K_1(\omega,n)&=\frac{(-1)^{n_1}}{\pi^3}\sum\limits_{\ell=0}^{N}\Big(\sum\limits_{k=0}^{\ell}\frac{(\partial^{\ell-k}_{r}\partial^{k}_{\rho_1}H_1)(0,0,n)}{k!(\ell-k)!}c_{\ell-k}\Big)\omega^{\ell+\frac{1}{2}}+\sum\limits_{k=0}^{N}\sum\limits_{j=0}^{N}d(k,j,N,n)\omega^j\\
&\quad +O(\left<n\right>^{2N}(1+\ln\left<n\right> )^2|\omega|^N),
\end{align*}
where $(\partial^{\ell}_{r}\partial^{k}_{\rho_1}H_1)(0,0,n)$ is a polynomial in $n$ of order $2(k+\ell)$ and $c_{\ell}$ denotes the coefficient of $\omega^{\ell+\frac{1}{2}}$ in the expansion of the integral below (see also \cite[Lemma 2.7]{Kop09})
\begin{equation*}
    \int_{0}^{\delta}\Big(\pi i-\ln \frac{1-\sqrt{\frac{r-\omega}{2r}}}{1+\sqrt{\frac{r-\omega}{2r}}}\Big)\frac{r^{\ell}dr}{\sqrt{r-\omega}},\qquad 0<|\omega|<\frac{\delta}{2},\qquad \operatorname{Im}\omega>0.
\end{equation*}
In particular, $c_0=-2\sqrt{2}\pi$ (this value corrects the miscalculation in the original \cite[Appendix A]{Kop09}). Applying this method to $j=2,3$, a similar expansion of $K_j(\omega,n)$ is obtained by replacing $H_1$ with $H_j$, where 
 \begin{equation*}
H_j(r,\rho_1,n)=\frac{\widetilde{\eta}_j(r)}{4\sqrt{2}}\mathbb{F}_j(n,\frac{r+\rho_1}{2},\frac{r-\rho_1}{2}),\quad \mathbb{F}_j(n,\tau^2,s^2_3)=\int_{0}^{\frac{\pi}{2}}F_j(n,\tau^2\cos^2\varphi,\tau^2\sin^2\varphi,s^2_3)d\varphi
 \end{equation*}
with $\widetilde{\eta}_j(r)$ satisfying the same property as $\widetilde{\eta}_1(r)$ and
\begin{equation*}
F_j(n,s^2_1,s^2_2,s^2_3)=G(s)
\begin{cases}
        {\rm cos}(n_1g(s_1)){\rm cos}(n_2g(s_3)){\rm cos}(n_3g(s_2)),&j=2\\
        {\rm cos}(n_1g(s_1)){\rm cos}(n_2g(s_2)){\rm cos}(n_3g(s_3)),&j=3
 \end{cases}\ .
\end{equation*}
Hence, we obtain 
\begin{equation*}
B^+_{\ell}(n)=\sum\limits_{j=1}^{3}\frac{(-1)^{n_j}}{\pi^3}\Big(\sum\limits_{k=0}^{\ell}\frac{(\partial^{\ell-k}_{r}\partial^{k}_{\rho_1}H_j)(0,0,n)}{k!(\ell-k)!}c_{\ell-k}\Big)
\end{equation*}
In particular, if $\ell=0$, we can compute all $H_j(0,0,n)=\frac{\pi}{8\sqrt{2}}$ and this together with $c_0=-2\sqrt{2}\pi$ gives 
 \begin{equation*}
B^+_0(n)=c_0\sum\limits_{j=1}^{3}\frac{(-1)^{n_j}}{\pi^3}H_j(0,0,n)=-\frac{1}{4\pi}\big((-1)^{n_1}+(-1)^{n_2}+(-1)^{n_3}\big).
\end{equation*}
Moreover, by the method used for $A_0$ in Proposition \ref{prop of Ak calcul}, we can also verify $(-\Delta-4)E^+_0=\delta_0.$
\end{proof}

\begin{lemma}\label{cauthy integral lemma}
Given $a>0$ and $B\in C^{\infty}_0((-a,a))$. Define 
\begin{equation*}
    F(\omega)=\int_{-a}^{a}\frac{B(u)}{u-\omega}du,\quad \omega =x+iy,\quad x\in \R,\quad y>0.
\end{equation*}
Assume that $|\omega|<\frac{a}{2}$. Then for any $0<h\leq \frac{a}{4}$, we have
\begin{equation}\label{esti of Fw}
\sup\limits_{0<|\omega|<\frac{a}{2},\ \operatorname{Im}\omega>0}|F(\omega)|\leq C_a\|B\|_{L^{\infty}}+2h\|B'\|_{L^{\infty}}+4\|B\|_{L^{\infty}}\ln\frac{3a}{2h}.
\end{equation}
Moreover, $$\lim\limits_{\omega\rightarrow0,\ \operatorname{Im}\omega>0}F(\omega)=\int_{-a}^{a}\frac{B(u)-B(0)}{u}du+B(0)i\pi.$$
\end{lemma}
\begin{proof}
 Let $\omega=x+iy$ with $y>0$ and $|\omega|<\frac{a}{2}$. We write
    \begin{equation}\label{decom of Fw}
    F(\omega)=B(x)\int_{-a}^{a}\frac{du}{u-\omega}+\int_{-a}^{a}\frac{B(u)-B(x)}{u-\omega}du=:F_1(\omega)+F_2(\omega).
    \end{equation}
 Since $|\omega|\leq a/2$, then by triangle inequality, $a/2\leq |a\pm \omega|\leq 3a/2$. We can compute 
\begin{equation*}
|F_1(\omega)|=\Big|B(x)\Big(\pi i+\ln \frac{a-\omega}{a+\omega}\Big)\Big|\leq C_a\|B\|_{L^{\infty}}.
\end{equation*}
For $F_2(\omega)$, we decompose
\begin{equation*}
    F_2(\omega)=\int_{-a}^{a}\frac{B(u)-B(x)}{u-\omega}du=\Big(\int\limits_{\substack{|u-x|<h\\ |u|<a}}+\int\limits_{\substack{|u-x|\geq h\\ |u|<a}}\Big)\frac{B(u)-B(x)}{u-\omega}du=:F_{21}(\omega)+F_{22}(\omega).
\end{equation*}
By the fundamental theorem of calculus, 
\begin{equation*}
    |F_{21}(\omega)|\leq \|B'\|_{L^{\infty}}\int\limits_{\substack{|u-x|<h\\ |u|<a}}\frac{|u-x|}{|u-x-iy|}du\leq 2h\|B'\|_{L^{\infty}}.
\end{equation*}
On the other hand, 
\begin{equation*}
    |F_{22}(\omega)|\leq 2\|B\|_{L^{\infty}}\int\limits_{\substack{|u-x|\geq h\\ |u|<a}}\frac{1}{|u-x|}du\leq 2\|B\|_{L^{\infty}}\int\limits_{h\leq|u|\leq \frac{3a}{2}}\frac{1}{|u|}du=4\|B\|_{L^{\infty}}\ln\frac{3a}{2h}.
\end{equation*}
This proves \eqref{esti of Fw}. By \eqref{decom of Fw} and the dominated convergence theorem, the limit follows.
\end{proof}

\section{Limiting Absorption Principle and Representation of \texorpdfstring{$e^{-itH}P_{ac}(H)$}{e-itH} }\label{LAP and repre of eitHpacH}
\subsection{Limiting Absorption Principle for \texorpdfstring{$-\Delta_{\Z^d}+V$}{Delta+V}}\label{LAP general d}
\ \

\vskip0.1cm
Let $d\geq1$. Recall  that
\begin{equation*}
(\Delta_{\Z^d}\phi)(n):=\sum_{j=1}^d(\phi(n+e_j)+\phi(n-e_j)-2\phi(n)),\quad e_j=(0,\cdots,0,\underset{j}{1},0,\cdots,0).
\end{equation*}
The limiting absorption principle (LAP) for discrete Schr\"odinger operators $-\Delta_{\Z^d}+V$ has received considerable attention over the past few decades. A classical result is due to Boutet de Monvel and Sahbani \cite{BS99}, who established the LAP in a Besov-space framework for potentials with arbitrarily slow decay at infinity, via Mourre theory. Subsequent works have treated the LAP under various assumptions on the potential and in different functional settings; see, for instance, \cite{IK12,Man17,GM21}.
  
By contrast, fewer results seem to explicitly address higher-order derivatives of the boundary resolvent in weighted $\ell^2$ spaces in the discrete setting. Motivated by \cite{HY25}, we establish the following result for the derivatives of the resolvent in $\mcaB(s,-s)$.
\begin{theorem}\label{general LAP}
    Let $H_d=-\Delta_{\Z^d}+V$ with a real-valued potential $V$ satisfying $|V(n)|\lesssim \left<n\right>^{-\beta}$ for some $\beta>2$ and $\mcaI=(0,4d)\setminus\{4j: 1\leq j\leq d-1,\ j\in \N^+\}$. Let $[\beta]$ denote the greatest integer not exceeding $\beta$. Then the following statements hold.
\begin{itemize}
\item [\rm(i)] The point spectrum $\sigma_p(H_d)\cap\mcaI$ is discrete, each eigenvalue has finite multiplicity, and the singular continuous spectrum $\sigma_{sc}(H_d)=\varnothing$.
\vskip0.2cm
\item [\rm(ii)] Let $j\in\left\{0,\cdots,[\beta]-2\right\}$ and $s>j+\frac{1}{2}$. Then 
the norm limits
    \begin{equation*}
   \frac{d^{j}}{d\lambda^j}(\mcaR^{\pm}_V(\lambda))=\lim\limits_{\varepsilon\downarrow0}R^{(j)}_V(\lambda\pm i\varepsilon) \quad {\rm in} \quad \mcaB(s,-s)
    \end{equation*}
exist and are norm-continuous on $\mcaI\setminus\sigma_p(H_d)$. Here $R_{V}(z)=(H_d-z)^{-1}$ is the resolvent of $H_d$ and $R^{(j)}_{V}(z)$ denotes the $j$-th derivative of $R_{V}(z)$.
\end{itemize}
\end{theorem}

Following the approach of \cite[Section 2]{HY25}, to derive this theorem, we introduce the conjugate operator 
$$A:=-i\sum\limits_{j=1}^{d}(N_jP_j-P_j^*N _j),$$
where 
$
(N_j\phi)(n):=n_j\phi(n)$ with $\phi\in D(N_j):=\big\{\phi\in\ell^2(\Z^d)\mid n_j\phi(n)\in\ell^2(\Z^d)\big\}
$
and 
\begin{equation*}
(P_j\phi)(n):=\phi(n+e_j)-\phi(n),\quad (P^*_j\phi)(n):=\phi(n-e_j)-\phi(n).
\end{equation*}
It then suffices to verify that $H_d$ satisfies the desired regularity with respect to $A$ (see Definition \ref{def of regularity} below) and the Mourre estimate \eqref{mourre estimate}. 
Recall that
\begin{definition}\label{def of regularity}
{\rm
Given an integer $k\in\N$, we say that $H_d$ belongs to the class $C^{k}(A)$, and write $H_d\in C^{k}(A)$, if the sesquilinear form $ad^{k}_{A}(H_d)$, initially defined on $\mcaD(A^k)\times\mcaD(A^k)$, extends continuously to $\ell^2(\Z^d)\times\ell^2(\Z^d)$. We identify this extension with the corresponding bounded operator on $\ell^2(\Z^d)$ and denote it by the same symbol. The iterated commutators $ad^{k}_{A}(H_d)$ are defined recursively by
\begin{align}\label{adk}
\begin{split}
ad^{0}_A(H_d)&=H_d,\quad
ad^{1}_{A}(H_d)=[H_d,A]=H_dA-AH_d,\\
ad^{k}_{A}(H_d)&=ad^{1}_{A}\left(ad^{k-1}_{A}(H_d)\right)=\sum\limits_{i,j\geq0,i+j=k}^{}\frac{k!}{i!j!}(-1)^{i}A^{i}H_dA ^{j}.
\end{split}
\end{align}
}
\end{definition}
First, the required regularity follows from \cite[Lemma 4.2, Proposition 5.1]{BS99} and \cite[Lemma A.5]{HY25}.
\begin{lemma}\label{regularity}
 Let $H_d=-\Delta_{\Z^d}+V$ with $|V(n)|\lesssim \left<n\right>^{-\beta}$ for some $\beta>0$. Then $H_d\in C^{[\beta]}(A)$.
\end{lemma}
We next establish the following Mourre estimate.
\begin{lemma}\label{mourre lemma}
 Let $H_d=-\Delta_{\Z^d}+V$ with $|V(n)|\lesssim \left<n\right>^{-\beta}$ for some $\beta>1$. For any $\lambda\in(0,4d)\setminus\{4j\}^{d-1}_{j=1}$, there exist constants $\delta,C>0$ and a compact operator $K$ such that
\begin{equation}\label{mourre estimate}
    E_{H_d}(J)[H_d,iA]E_{H_d}(J)\geq C E_{H_d}(J)+K,\quad J=(\lambda-\delta,\lambda+\delta),
\end{equation}
where $E_{H_d}(J)$ denotes the spectral projection of $H_d$ onto the interval $J$.
\end{lemma}
\begin{proof}
Under the Fourier transform $\mcaF_d$ on $\Z^d$:
\begin{equation*}
    (\mcaF_d\phi)(x):=(2\pi)^{-d/2}\sum_{n\in\Z^d}^{}e^{-in\cdot x}\phi(n), \quad \phi\in\ell^2(\Z^d),
\end{equation*}
it suffices to prove the corresponding Mourre estimate on
$L^2((-\pi,\pi]^d)$. Set
$$\mcaH:=\mcaH_0+\mcaV,\qquad \mcaV:=\mcaF_d V\mcaF_d^{-1},\qquad \mcaA:=\mcaF_d A\mcaF^{-1}_d,$$ where
\begin{equation*}
\mcaH_0:=\mcaF_d (-\Delta_{\Z^d})\mcaF^{-1}_d=2\big(d-\sum\limits_{j=1}^{d}{\rm cos}x_j\big)=:h(x),
\end{equation*}
and 
\begin{equation*}
  \mcaA=\sum\limits_{j=1}^{d}(2i{\rm sin}x_j\partial_{x_j}+i{\rm cos}x_j-{\rm sin}x_j).
\end{equation*}
We first establish the Mourre estimate for $\mcaH_0$. Let  $\lambda\in(0,4d)\setminus\{4j\}^{d-1}_{j=1}$. A direct computation gives
\begin{equation*}
[\mcaH_0,i\mcaA]=4\sum\limits_{j=1}^{d}({\rm sin}x_j)^2=|\nabla h(x)|^2.
 \end{equation*}
Let 
 $M_{\lambda}=\{x\in(-\pi,\pi]^d:h(x)=\lambda\}.$ Since $\lambda$ is not a threshold, $\nabla h(x)\neq0$ on
$M_\lambda$. Moreover, $M_\lambda$ is compact, and therefore
$$
C(\lambda)
:=
\inf_{x\in M_\lambda}|\nabla h(x)|^2
>0.
$$
Choose $0<\eta<C(\lambda)$. By continuity, there exists
$\varepsilon>0$ such that
$$
|\nabla h(x)|^2
\geq C(\lambda)-\eta,
\qquad
x\in h^{-1}([\lambda-\varepsilon,\lambda+\varepsilon]).
$$
Hence, with $I=(\lambda-\varepsilon,\lambda+\varepsilon)$, we obtain
\begin{equation}\label{mourre estimate free}
E_{\mcaH_0}(I)
[\mcaH_0,i\mcaA]
E_{\mcaH_0}(I)
\geq
\big(C(\lambda)-\eta\big)E_{\mcaH_0}(I).
\end{equation}
We next pass from $\mcaH_0$ to $\mcaH$. Since
$V(n)\to0$ as $|n|\to\infty$, $V$ is compact on $\ell^2(\Z^d)$,
and hence $\mcaV$ is compact on $L^2((-\pi,\pi]^d)$. It follows that
$E_{\mcaH}(I)-E_{\mcaH_0}(I)$ is compact. Moreover, let \begin{equation*}
\widetilde{A}=A-\frac{i}{2}\sum\limits_{j=1}^{d}(P_j-P^*_j)=\frac{i}{2}\sum\limits_{j=1}^{d}\big((P^*_j-P_j)N_j+N_j(P^*_j-P_j)\big),\quad \Phi(x)=-\sum\limits_{j=1}^{d}{\rm sin}x_j.
\end{equation*}
Denote $\widetilde{\mcaA}=\mcaF_d\widetilde{A}\mcaF^{-1}_d$. Since 
$[V,i\widetilde{A}]$ is compact on $\ell^2(\Z^d)$ (cf. \cite[Page 17]{Man17}), it follows that
$$[\mcaV,i\mcaA]=[\mcaV,i\widetilde{\mcaA}]+[\mcaV,i\Phi]=\mcaF_d[V,i\widetilde{A}]\mcaF^{-1}_d+[\mcaV,i\Phi]$$
is compact on $L^2((-\pi,\pi]^d)$. Combining these facts with \eqref{mourre estimate free} and the decompositions
$$E_{\mcaH}(I)=E_{\mcaH_0}(I)+E_{\mcaH}(I)-E_{\mcaH_0}(I),\quad [\mcaH,i\mcaA]=[\mcaH_0,i\mcaA]+[\mcaV,i\mcaA],$$
we obtain the desired estimates.
\end{proof}
 Theorem \ref{general LAP} now follows from Lemmas \ref{regularity} and \ref{mourre lemma}, together with \cite[Proposition 2.6 and Lemma 2.8]{HY25}.

\subsection{Stationary Representation of \texorpdfstring {$e^{-itH}P_{ac}(H)$}{e-ith}} \label{Stationary repres of eitHpacH}
\ \

\vskip0.1cm
In this subsection, under the assumptions of Theorem \ref{decay theorem}, we prove the stationary representation
\begin{equation}\label{eitHpacH}
   e^{-itH}P_{ac}(H)=\frac{1}{2\pi i}\sum\limits_{k=1}^{3}\int_{4(k-1)}^{4k}e^{-it\lambda}(\mcaR_V^{+}-\mcaR_V^{-})(\lambda)d\lambda,\quad H=-\Delta_{\Z^3}+V.
\end{equation}
\begin{proof}
By the assumptions on $V$ and Weyl's theorem, the essential spectrum of $H$ satisfies 
$$\sigma_{ess}(H)=\sigma_{ess}(-\Delta_{\Z^3})=[0,12],$$
and the discrete spectrum $\sigma_{dis}(H)\subseteq\R\setminus[0,12]$.
Then, by Cauchy's integral formula, we have
\begin{align*}
e^{-itH}&=-\frac{1}{2\pi i}\ointctrclockwise_{|\omega|=C}e^{-it\omega}R_V(\omega)d\omega,\qquad C>\sup\{{12,\{|\lambda_j|\}_{\lambda_j\in\sigma_{dis}(H)}}\},\\
&=\sum\limits_{\lambda_j\in\sigma_{dis}(H)}e^{-it\lambda_j}\mcaP_{\lambda_j}-\frac{1}{2\pi i}\int_{\Gamma(\varepsilon)}e^{-it\omega}R_V(\omega)d\omega,\quad 0<\varepsilon<1,
\end{align*}
where $\mcaP_{\lambda_j}$ is the projection on the eigenspace corresponding to the eigenvalue $\lambda_j$ and the path $\Gamma(\varepsilon)$ is chosen as in Figure \ref{myplot} below. 
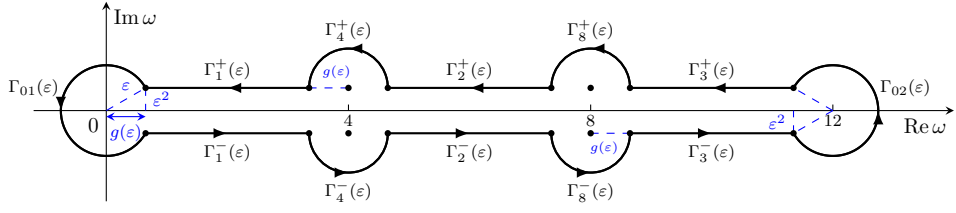
\begin{figure}[ht]
\centering
\captionsetup{labelformat=empty}
\centering
\begin{tikzpicture}[scale=0.8, transform shape, >=stealth]
\draw[->] (-1.2,0)--(14,0) node[below left]{$\operatorname{Re} \omega$};
\draw[->] (0,-1.5)--(0,1.8) node[below right]{$\operatorname{Im} \omega$};
\tikzset{ticklabel/.style={font=\footnotesize, anchor=north, inner sep=0pt, yshift=-1pt}}
\tikzset{major_tick/.style={font=\footnotesize, anchor=north, inner sep=0pt, yshift=-2pt}}
\foreach \x in {4,8,12} { 
    \draw (\x, 0) -- (\x, 0.1);
    \node[major_tick] at (\x, 0) {$\x$}; }
\node[below left] at (0,0) {$0$};
\draw[black, thick] (30:0.75) arc [start angle=30, end angle=330, radius=0.75] node[left, pos=0.4, font=\footnotesize] {$\Gamma_{01}(\varepsilon)$};
\draw[dashed, blue] (0,0) -- node[pos=0.5, above=1pt, font=\footnotesize] {$\varepsilon$} (30:0.75);
\coordinate (C1) at (30:0.75);
\coordinate (C2) at (330:0.75);
\draw[dashed, blue] (C1) -- node[midway, right, font=\footnotesize] {$\varepsilon^2$} (C1 |- 0,0);
\begin{scope} [shift={(12,0)}]
\draw[black, thick] (150:0.75) arc [start angle=150, end angle=-150, radius=0.75] node[right, pos=0.4, font=\footnotesize] {$\Gamma_{02}(\varepsilon)$};
 \draw[dashed, blue] (0,0) -- (150:0.75);
 \draw[dashed, blue] (0,0) -- (-150:0.75);
 \coordinate (D1) at (150:0.75);
\coordinate (D2) at (-150:0.75);
\draw[dashed, blue] (D2) -- node[midway, left, font=\footnotesize] {$\varepsilon^2$} (D2 |- 0,0);
\end{scope}
\begin{scope}[shift={(4,0.375)}]
\draw[black, thick] (0:0.65) arc [start angle=0, end angle=180, radius=0.65] node[above, pos=0.5, font=\footnotesize] {$\Gamma^{+}_{4}(\varepsilon)$};
    \draw[dashed, blue] (0,0)  -- node[pos=0.3, above=1pt, font=\tiny] {$g(\varepsilon)$}  (180:0.65);
    \coordinate (A1) at (0:0.65);
    \coordinate (A2) at (180:0.65);
    \coordinate (A5) at (0:0);
\end{scope}
\begin{scope}[shift={(4,-0.375)}]
\draw[black, thick] (180:0.65) arc [start angle=180, end angle=360, radius=0.65] node[below, pos=0.5, font=\footnotesize] {$\Gamma^{-}_{4}(\varepsilon)$};
    \coordinate (A3) at (180:0.65);
    \coordinate (A4) at (360:0.65);
    \coordinate (A6) at (0:0);
\end{scope}
\begin{scope}[shift={(8,-0.375)}]
\draw[black, thick] (180:0.65) arc [start angle=180, end angle=360, radius=0.65] node[below, pos=0.5, font=\footnotesize] {$\Gamma^{-}_{8}(\varepsilon)$};
\draw[dashed, blue] (0,0)  -- node[pos=0.3, below=0.5pt, font=\tiny] {$g(\varepsilon)$}  (360:0.75);
    \coordinate (B3) at (180:0.65);
    \coordinate (B4) at (360:0.65);
     \coordinate (B6) at (0:0);
\end{scope}
\begin{scope}[shift={(8,0.375)}]
\draw[black, thick] (0:0.65) arc [start angle=0, end angle=180, radius=0.65] node[above, pos=0.5, font=\footnotesize] {$\Gamma^{+}_{8}(\varepsilon)$};
    \coordinate (B1) at (0:0.65);
    \coordinate (B2) at (180:0.65);
    \coordinate (B5) at (0:0);
\end{scope}
\coordinate (proj) at ({0.75*cos(30)}, -0.1);
\draw[<->, >=stealth, blue!70!blue] (0,-0.1) -- node[midway, below, font=\footnotesize] {$g(\varepsilon)$} (proj);
  \tikzset{
    ccwarrow/.style={
      decoration={
        markings,
        mark=at position 0.5 with {\arrow{latex }}
      },
      postaction={decorate}
    }
  }
  \draw[black, thick, ccwarrow] (30:0.75) arc [start angle=30, end angle=330, radius=0.75];
\begin{scope}[shift={(12,0)}]
    \draw[black, thick, ccwarrow] (-150:0.75) arc [start angle=-150, end angle=150, radius=0.75];
\end{scope}
\begin{scope}[shift={(4,0.375)}]
    \draw[black, thick, ccwarrow] (0:0.65) arc [start angle=0, end angle=180, radius=0.65];
\end{scope}
\begin{scope}[shift={(4,-0.375)}]
    \draw[black, thick, ccwarrow] (180:0.65) arc [start angle=180, end angle=360, radius=0.65];
\end{scope}
\begin{scope}[shift={(8,0.375)}]
    \draw[black, thick, ccwarrow] (0:0.65) arc [start angle=0, end angle=180, radius=0.65];
\end{scope}
\begin{scope}[shift={(8,-0.375)}]
    \draw[black, thick, ccwarrow] (180:0.65) arc [start angle=180, end angle=360, radius=0.65];
\end{scope}
\begin{scope}[shift={(4,0)}]
    \fill (A1) circle (1.5pt); 
    \fill (A2) circle (1.5pt);
    \fill (A3) circle (1.5pt);
    \fill (A4) circle (1.5pt);
    \fill (A5) circle (1.5pt);
    \fill (A6) circle (1.5pt);
\end{scope}
\begin{scope}[shift={(8,0)}]
    \fill (B1) circle (1.5pt); 
    \fill (B2) circle (1.5pt);
    \fill (B3) circle (1.5pt);
    \fill (B4) circle (1.5pt);
    \fill (B5) circle (1.5pt);
    \fill (B6) circle (1.5pt);
\end{scope}    
\begin{scope}[shift={(12,0)}]
    \fill (D1) circle (1.5pt); 
    \fill (D2) circle (1.5pt);
\end{scope}
\fill (C1) circle (1.5pt); 
\fill (C2) circle (1.5pt);
\coordinate (E1) at (3,0.36);
\coordinate (E2) at (3,-0.36);
\coordinate (E3) at (5,0.36);
\coordinate (E4) at (5,-0.36);
\coordinate (F1) at (7,0.36);
\coordinate (F2) at (7,-0.36);
\coordinate (F3) at (9,0.36);
\coordinate (F4) at (9,-0.36);
\draw[black, thick, ccwarrow] (D1) -- (B1) node[above, pos=0.5, font=\footnotesize] {$\Gamma^+_{3}(\varepsilon)$};
\draw[black, thick, ccwarrow] (B2) -- (A1) node[above, pos=0.5, font=\footnotesize] {$\Gamma^+_{2}(\varepsilon)$};
\draw[black, thick,ccwarrow] (A2) -- (C1) node[above, pos=0.5, font=\footnotesize] {$\Gamma^+_{1}(\varepsilon)$};
\draw[black, thick, ccwarrow] (C2) -- (A3) node[below, pos=0.5, font=\footnotesize] {$\Gamma^-_{1}(\varepsilon)$};
\draw[black, thick, ccwarrow] (A4) -- (B3) node[below, pos=0.5, font=\footnotesize] {$\Gamma^-_{2}(\varepsilon)$};
\draw[black, thick, ccwarrow] (B4) -- (D2) node[below, pos=0.5, font=\footnotesize] {$\Gamma^-_{3}(\varepsilon)$};
\end{tikzpicture}
\caption{Figure 1: The path $\Gamma({\varepsilon})$.}
\label{myplot}
\end{figure}

\noindent Here the paths $\Gamma_{0\ell}(\varepsilon)$, $\Gamma^{\pm}_{j}(\varepsilon)$ and $\Gamma^{\pm}_{4k}(\varepsilon)$ are given as follows:
\begin{align}\label{expre of path Gammas}
\begin{split}
\Gamma_{01}(\varepsilon)&=\big\{|\omega|=\varepsilon~\big|~\arg\omega:\theta(\varepsilon)\to2\pi-\theta(\varepsilon)\big\},\quad {\rm sin}\theta(\varepsilon)=\varepsilon,\\
\Gamma^+_{j}(\varepsilon)&=\big\{\omega=\lambda+i\varepsilon^2~\big|~\lambda:4j-g(\varepsilon)\to4(j-1)+g(\varepsilon)\big\},\quad g(\varepsilon)=\varepsilon\sqrt{1-\varepsilon^2},\\
\Gamma^-_{j}(\varepsilon)&=\big\{\omega=\lambda-i\varepsilon^2~\big|~\lambda:4(j-1)+g(\varepsilon)\to4j-g(\varepsilon)\big\},\quad j=1,2,3,\\
\Gamma^{\pm}_{4k}(\varepsilon)&=\big\{|\omega-(4k\pm i\varepsilon^2)|=g(\varepsilon)~\big| \pm{\rm Im}\omega\geq\varepsilon^2\big\},\quad k=1,2,\\
\Gamma_{02}(\varepsilon)&=\big\{|\omega-12|=\varepsilon~\big|\arg(\omega-12):\pi+\theta(\varepsilon)\to3\pi-\theta(\varepsilon)\big\}.
\end{split}
\end{align}
Then
\begin{align}
&\qquad\quad e^{-itH}-   \sum\limits_{\lambda_j\in\sigma_{dis}(H)}e^{-it\lambda_j}\mcaP_{\lambda_j}\notag\\
&=-\frac{1}{2\pi i}\lim\limits_{\varepsilon\downarrow0}\left[\sum\limits_{j=1}^{3}\int_{\underset{\pm}{\bigcup}\Gamma^{\pm}_j(\varepsilon)}e^{-it\omega}R_V(\omega)d\omega+\sum\limits_{k=1}^{2}\int_{\underset{\pm}{\bigcup}\Gamma^{\pm}_{4k}(\varepsilon)}e^{-it\omega}R_V(\omega)d\omega+\sum\limits_{\ell=1}^{2}\int_{\Gamma_{0\ell}(\varepsilon)}e^{-it\omega}R_V(\omega)d\omega\right]\notag\\
&=:-\frac{1}{2\pi i}\lim\limits_{\varepsilon\downarrow0}\Big(\sum\limits_{j=1}^{3}I_j(\varepsilon)+\sum\limits_{k=1}^{2}I_{4k}(\varepsilon)+\sum\limits_{l=1}^{2}I_{0\ell}(\varepsilon)\Big).\label{eitH}
\end{align}
\vskip0.2cm
\underline{\bf Case (1): $I_j(\varepsilon)$}. We prove that, for each $j=1,2,3$, 
\begin{equation}\label{Ij123}
\lim\limits_{\varepsilon\downarrow0}I_j(\varepsilon)=\int_{4(j-1)}^{4j}e^{-it\lambda}(\mcaR^{-}_V-\mcaR^{+}_{V})(\lambda)d\lambda\quad {\rm in}\ \ \mcaB(\sigma,-\sigma),\quad \sigma>11/2. 
\end{equation}
Set $a_j=4(j-1)$. Fix $0<\mu_0\ll1$ sufficiently small and let  $0<\varepsilon<\mu_0$. Then  
\begin{align}\label{Ij}
\begin{split}
I_j(\varepsilon)&=\Big(\int_{\Gamma^+_j(\varepsilon)}+\int_{\Gamma^-_j(\varepsilon)}\Big)e^{-it\omega}R_V(\omega)d\omega\\
&=\Big(\int_{a_j+g(\varepsilon)}^{a_j+\mu_0}+\int_{a_j+\mu_0}^{a_{j+1}-\mu_0}+\int_{a_{j+1}-\mu_0}^{a_{j+1}-g(\varepsilon)}\Big)e^{-it\lambda}\Big(e^{-t\varepsilon^2}R_V(\lambda-i\varepsilon^2)-e^{t\varepsilon^2}R_V(\lambda+i\varepsilon^2)\Big)d\lambda\\
&=:I_{j1}(\varepsilon)+I_{j2}(\varepsilon)+I_{j3}(\varepsilon).
\end{split}
\end{align}
By Theorems \ref{Asymptotic expansion theorem_0} and \ref{Asymptotic expansion theorem_4}, the resolvent $R_V(a_j+\omega)$ near the threshold $a_j$ admits the unified expansion 
\begin{equation}\label{unified expan of Rv}
    R_V(a_j+\omega)=-\widetilde{\mcaP}_{a_j}\omega^{-1}+\mcaE_{a_j}(\omega), \quad 0<|\omega|\leq 2\mu_0,
\end{equation}
in the space $\mcaB(\sigma,-\sigma)$. Here $\widetilde{\mcaP}_{a_j}=\mcaP_{a_j}$ if $a_j$ is an eigenvalue of $H$, and $\widetilde{\mcaP}_{a_j}=0$ otherwise. Moreover, 
\begin{equation}\label{estim of mcaEaj}
    \|\mcaE_{a_j}(\omega)\|_{\mcaB(\sigma,-\sigma)}\lesssim|\omega|^{-1/2},\quad 0<|\omega|\leq2\mu_0.
\end{equation}
We first consider $I_{j1}(\varepsilon)$. Making the change of variables $x=\lambda-a_j$ and using \eqref{unified expan of Rv}, we obtain
\begin{align*}
    I_{j1}(\varepsilon)&=-\widetilde{\mcaP}_{a_j}\int_{g(\varepsilon)}^{\mu_0}e^{-it(x+a_j)}\Big(  \frac{e^{-t\varepsilon^2}}{x-i\varepsilon^2}-\frac{e^{t\varepsilon^2}}{x+i\varepsilon^2}\Big)dx\\
&\quad+\int_{g(\varepsilon)}^{\mu_0}e^{-it(x+a_j)}\Big(e^{-t\varepsilon^2}\mcaE_{a_j}(x-i\varepsilon^2)-e^{t\varepsilon^2}\mcaE_{a_j}(x+i\varepsilon^2)\Big)dx\\
&:=I^{(1)}_{j1}(\varepsilon)+I^{(2)}_{j1}(\varepsilon).
\end{align*}
For fixed $t$ and sufficiently small  $\varepsilon$, we have
\begin{equation*}
    \Big|\frac{e^{-t\varepsilon^2}}{x-i\varepsilon^2}-\frac{e^{t\varepsilon^2}}{x+i\varepsilon^2}\Big|=\Big|\frac{x(e^{-t\varepsilon^2}-e^{t\varepsilon^2})+i\varepsilon^2(e^{-t\varepsilon^2}+e^{t\varepsilon^2})}{x^2+\varepsilon^4}\Big|\lesssim_t\frac{x\varepsilon^2+\varepsilon^2}{x^2+\varepsilon^4}.
\end{equation*}
Consequently,
\begin{equation}\label{I1j1}
    \|I^{(1)}_{j1}(\varepsilon)\|_{\mcaB(\sigma,-\sigma)}\lesssim_t\varepsilon^2\Big(\int_{g(\varepsilon)}^{\mu_0}\frac{dx}{x}+\int_{g(\varepsilon)}^{\mu_0}\frac{dx}{x^2}\Big)\lesssim_t\varepsilon+\varepsilon^2|\ln g(\varepsilon)|\rightarrow0,\quad \varepsilon\rightarrow0.
\end{equation}
We next consider $I^{(2)}_{j1}(\varepsilon)$. By \eqref{estim of mcaEaj}, we note that 
\begin{equation*}
    \left\|e^{-it(x+a_j)}\Big(e^{-t\varepsilon^2}\mcaE_{a_j}(x-i\varepsilon^2)-e^{t\varepsilon^2}\mcaE_{a_j}(x+i\varepsilon^2)\Big)\right\|_{\mcaB(\sigma,-\sigma)}\lesssim_t x^{-1/2}\in L^{1}((0,\mu_0]).
\end{equation*}
For every fixed $x\in (0,\mu_0]$, the limiting absorption principle gives
\begin{equation*}
\lim\limits_{\varepsilon\rightarrow0}\mcaE_{a_j}(x\pm i\varepsilon^2)=\lim\limits_{\varepsilon\rightarrow0}\Big(R_V(a_j+x\pm i\varepsilon^2)+\widetilde{\mcaP}_{a_j}(x\pm i\varepsilon^2)^{-1}\Big)=\mcaR^{\pm}_V(a_j+x)+\widetilde{\mcaP}_{a_j}x^{-1}.
\end{equation*}
Therefore, by Lebesgue's dominated convergence theorem, we obtain
\begin{equation}\label{I2j1}
\lim\limits_{\varepsilon\rightarrow0}I^{(2)}_{j1}(\varepsilon)=\int_{0}^{\mu_0}e^{-it(x+a_j)}(\mcaR^-_V(x+a_j)-\mcaR^+_V(x+a_j))dx=\int_{a_j}^{a_j+\mu_0}e^{-it\lambda}(\mcaR^-_V(\lambda)-\mcaR^+_V(\lambda))d\lambda.
\end{equation}
Similarly, we can also prove 
\begin{equation}\label{Ij3}
\lim\limits_{\varepsilon\rightarrow0}{I}_{j3}(\varepsilon)= \int_{a_{j+1}-\mu_0}^{a_{j+1}}e^{-it\lambda}(\mcaR^-_V(\lambda)-\mcaR^+_V(\lambda))d\lambda.
\end{equation}
It remains to consider the interior part $I_{j2}(\varepsilon)$. Let $J=[a_j+\mu_0,a_{j+1}-\mu_0]$. By the limiting absorption principle, the absence of embedded eigenvalues in the open spectral band, and the compactness of $J$, we have (cf. \cite[Proposition 2.6]{HY25})
\begin{equation*}
    \sup\limits_{\lambda\in J,\ 0<\varepsilon<\mu_0}\|R_V(\lambda\pm i\varepsilon^2)\|_{\mcaB(\sigma,-\sigma)}<\infty.
\end{equation*}
Therefore, another application of the Lebesgue's dominated convergence theorem gives 
\begin{equation*}
\lim\limits_{\varepsilon\rightarrow0}{I}_{j2}(\varepsilon)= \int_{a_{j}+\mu_0}^{a_{j+1}-\mu_0}e^{-it\lambda}(\mcaR^-_V(\lambda)-\mcaR^+_V(\lambda))d\lambda.
\end{equation*}
Combining this limit with \eqref{I1j1}, \eqref{I2j1}, \eqref{Ij3}, and the decomposition \eqref{Ij}, we obtain \eqref{Ij123}.
\vskip0.2cm
\underline{\bf Case (2): $I_{4k}(\varepsilon)$}. We claim 
\begin{equation}\label{I4k}
\lim\limits_{\varepsilon\downarrow0}I_{4k}(\varepsilon)=\begin{cases}
0,&4k~\text{is a regular point of }H \\
-2\pi i e^{-4ikt} \mathcal{P}_{4k},&4k~\text{is an eigenvalue of }H
\end{cases}\quad {\rm in }\ \mathcal{B}(\sigma,-\sigma)
\end{equation}
with $\sigma>\frac{7}{2}$ in the regular case and $\sigma>\frac{11}{2}$ in the eigenvalue case. Indeed, for any $k=1,2$, recall 
\begin{equation}\label{I4keps}
    I_{4k}(\varepsilon)=\Big(\int_{\Gamma^+_{4k}(\varepsilon)}+\int_{\Gamma^-_{4k}(\varepsilon)}\Big)e^{-it\omega}R_V(\omega)d\omega=:I^+_{4k}(\varepsilon)+I^-_{4k}(\varepsilon).
\end{equation}
By \eqref{expre of path Gammas} and the change of variable $\omega=4k\pm i\varepsilon^2+g(\varepsilon) e^{i\theta}$, we have
\begin{equation*}
I^{\pm}_{4k}(\varepsilon)=ig(\varepsilon) e^{\pm t\varepsilon^2-4ikt}
\left\{
    \begin{aligned}
\int_{0}^{\pi}e^{i\theta-itg(\varepsilon) e^{i\theta}}R_V(4k+ i\varepsilon^2+g(\varepsilon )e^{i\theta})d\theta,&\quad ``+",\\
        \int_{\pi}^{2\pi}e^{i\theta-itg(\varepsilon) e^{i\theta}}R_V(4k- i\varepsilon^2+g(\varepsilon) e^{i\theta})d\theta,&\quad ``-".
    \end{aligned}
    \right.
\end{equation*}
(i) If $4k$ is a regular point of $H$, from \eqref{(t)asy expan of regular 4}, we have 
\begin{equation*}
\|R_V(4k+\omega)\|_{\mcaB(\sigma,-\sigma)}\lesssim 1,\quad |\omega|\ll1/2, \quad \sigma>\frac{7}{2}.
\end{equation*}
Hence, 
\begin{equation}\label{lim of regular I4keps}
  \|I_{4k}(\varepsilon)\|_{\mcaB(\sigma,-\sigma)} \leq\sum\limits_{\pm}\|I^{\pm}_{4k}(\varepsilon)\|_{\mcaB(\sigma,-\sigma)}\lesssim g(\varepsilon) e^{|t|(g(\varepsilon)+ \varepsilon^2)}\rightarrow0,\quad \varepsilon\rightarrow0,\quad \sigma>\frac{7}{2}.
\end{equation}
(ii) If $4k$ is an eigenvalue of $H$, by \eqref{(t)asy expan eigenvalue 4}, as $|\omega|\ll 1/2$, 
\begin{equation*}
    R_V(4k+\omega)=-\mcaP_{4k}\omega^{-1}+\widetilde{R}_V(4k+\omega),\quad {\rm in }~\mcaB(\sigma,-\sigma),\quad \sigma>\frac{11}{2},
\end{equation*}
where $ \|\widetilde{R}_V(4k+\omega)\|_{\mcaB(\sigma,-\sigma)}\lesssim |\omega|^{-\frac{1}{2}}$. Since $|g(\varepsilon) e^{i\theta}\pm i\varepsilon^2|^{-\frac{1}{2}}\lesssim \varepsilon^{-\frac{1}{2}}$, a similar argument to (i) yields
\begin{equation*}
\lim\limits_{\varepsilon\downarrow0}I^{\pm}_{4k}(\varepsilon)=-e^{-4ikt}\mcaP_{4k}\lim\limits_{\varepsilon\downarrow0}\left\{\begin{aligned}
    \int_{0}^{\pi}\frac{ e^{i\theta-itg(\varepsilon) e^{i\theta}}}{\varepsilon(1-\varepsilon^2)^{-1/2}-ie^{i\theta}}d\theta,\quad ``+"\\
    \int_{\pi}^{2\pi}\frac{e^{i\theta-itg(\varepsilon) e^{i\theta}}}{-\varepsilon(1-\varepsilon^2)^{-1/2}-ie^{i\theta}}d\theta,\quad ``-"
\end{aligned}\right.=-\pi ie^{-4ikt}\mcaP_{4k},
\end{equation*}
which together with \eqref{I4keps} and \eqref{lim of regular I4keps} gives \eqref{I4k}.

 \underline{\bf Case (3):  $I_{0\ell}(\varepsilon)$}. By Theorem \ref{Asymptotic expansion theorem_0} and an argument similar to $I_{4k}(\varepsilon)$, one can verify  
\begin{equation*}
\lim\limits_{\varepsilon\downarrow0}I_{01}(\varepsilon)~({\rm resp.~} I_{02}(\varepsilon))=\begin{cases}
    0,&0~({\rm resp.~12})~\text{is a regular point of}~ H\\
    0,&0~({\rm resp.~12})~\text{is a resonance of}~ H\\
    -2\pi i\mcaP_0~ ({\rm resp.~}-2\pi ie^{-12it}\mcaP_{12}),&0~({\rm resp.~12})~\text{is an eigenvalue of}~ H
\end{cases}
\end{equation*}
in $\mathcal{B}(\sigma,-\sigma)$ with $\sigma>5/2$ for the regular and resonance cases, and $\sigma>7/2$ for the eigenvalue case. Therefore, combining 
this formula and \eqref{eitH},\eqref{Ij123} and \eqref{I4k}, we obtain  \eqref{eitHpacH}.
\end{proof}
\vskip0.4cm
\subsection*{Acknowledgements}    
The authors are partially supported by NSFC (grants No. 12671120 and 12531005). We would like to express their thanks to Professor Avy Soffer for his interests and insightful discussions.
	
\subsection*{Data Availability}
Data sharing is not applicable to this article as no new data were created or analyzed in this study.
\subsection*{Competing Interests} The authors have no relevant financial or non-financial interests to disclose.
\normalem

\end{document}